\documentclass[12pt, righttag]{article}
\usepackage{amsmath}
\usepackage{amscd}
\usepackage{graphicx}
\usepackage{amsthm}
\usepackage{amsfonts}
\usepackage{amssymb}
\newtheorem{theorem}{Theorem}[section]

\newtheorem{lemma}[theorem]{Lemma}

\newtheorem{proposition}[theorem]{Proposition}
\newtheorem{remark}[theorem]{Remark}

\numberwithin{equation}{section}

\begin{document}

\title{The homotopy principle in the existence level for maps with only singularities
of types $A$, $D$ and $E$}
\author{YOSHIFUMI ANDO}
\date{}
\maketitle
\begin{abstract}
Let $N$ and $P$ be smooth manifolds of dimensions $n$ and $p$ ($n\geq p\geq2$)
respectively. Let $\Omega(N,P)$ denote an open subspace of $J^{\infty}(N,P)$
which consists of all regular jets and jets with prescribed singularities of
types $A_{i}$, $D_{j}$ and $E_{k}$. An $\Omega$-regular map $f:N\rightarrow P$
refers to a smooth map having only singularities in $\Omega(N,P)$ and
satisfying the transversality condition. We will prove what is called the
homotopy principle for $\Omega$-regular maps in the existence level. Namely, a
continuous section $s$ of $\Omega(N,P)$ over $N$ has an $\Omega$-regular map
$f$ such that $s$ and $j^{\infty}f$ are homotopic as sections.
\end{abstract}

\section*{Introduction}

Let $N$ and $P$ be smooth ($C^{\infty}$) manifolds of dimensions $n$ and $p$
respectively with $n\geq p\geq2$. Let $A_{k}$, $D_{k}$ and $E_{k}$ denote the
types of the singularities of function germs studied in [Ar]. We say that a
$C^{\infty}$ stable map germ $f:(N,x_{0})\rightarrow(P,y_{0})$ has a
singularity of type $A_{k}$, $D_{k}$ or $E_{k}$, when $f$ is $C^{\infty}$
equivalent to $g:(\mathbf{R}^{n},\mathbf{0})\rightarrow(\mathbf{R}%
^{p},\mathbf{0})$, which is a versal unfolding of the function germ with
respective singularities as follows. Here, we assume $n>p\geq2$ only when we
deal with the singularity of types $D_{k}$ and $E_{k}$. We take the
coordinates%
\[%
\begin{array}
[c]{ll}%
x=(x_{1},\ldots,x_{p-k},x_{p},\ldots,x_{n-1},t_{0},\ldots,t_{k-2},u) &
\text{for }A_{k}\text{,}\\
x=(x_{1},\ldots,x_{p-k},x_{p},\ldots,x_{n-2},t_{0},\ldots,t_{k-2},u,\ell) &
\text{for }D_{k}\text{ and }E_{k},\text{ }%
\end{array}
\]
and $(y_{1},\ldots,y_{p})$ of $\mathbf{R}^{p}$. Then $g$ is expressed by%
\[%
\begin{array}
[c]{ll}%
y_{i}\circ g(x)=x_{i} & \text{for }1\leq i\leq p-k,\\
y_{i}\circ g(x)=t_{i-p+k-1} & \text{for }p-k<i\leq p-1,
\end{array}
\]
and by $y_{p}\circ g(x)$, which is written as

\begin{description}
\item $(A_{k})$ $\pm u^{k+1}+\sum_{i=1}^{k-1}t_{i-1}u^{i}\pm x_{p}^{2}%
\pm\cdots\pm x_{n-1}^{2}$ ($k\geq1),$

\item $(D_{k})$ $u^{2}\ell\pm\ell^{k-1}+t_{0}u+t_{1}\ell+\sum_{i=2}^{k-2}%
t_{i}\ell^{i}+Q\ (k\geq4),$

\item $(E_{6})$ $u^{3}\pm\ell^{4}+t_{0}u+t_{1}\ell+t_{2}u\ell+t_{3}\ell
^{2}+t_{4}u\ell^{2}+Q,$

\item $(E_{7})$ $u^{3}+u\ell^{3}+t_{0}u+t_{1}\ell+t_{2}u\ell+t_{3}\ell
^{2}+t_{4}\ell^{3}+t_{5}\ell^{4}+Q,$

\item $(E_{8})$ $u^{3}+\ell^{5}+t_{0}u+t_{1}\ell+t_{2}u\ell+t_{3}\ell
^{2}+t_{4}u\ell^{2}+t_{5}\ell^{3}+t_{6}u\ell^{3}+Q,$
\end{description}

\noindent where $Q=\pm x_{p}^{2}\pm\cdots\pm x_{n-2}^{2}$ ([Ar], [Math1,
(5.8)] and [Mo]).

In [B] there has been defined what is called the Boardman-Thom manifold
$\Sigma^{I}(N,P)$ in $J^{\infty}(N,P)$. For the symbol $I=(n-p+1,\overbrace
{1,\cdots,1}^{k-1},0),$ a smooth map germ $f:(N,x)\rightarrow(P,y)$ has $x$ as
a singularity of type $A_{k}$ if and only if $j_{x}^{\infty}f\in\Sigma
^{I}(N,P)$ and $j^{\infty}f$ is transverse to $\Sigma^{I}(N,P)$ (see [Mo]). By
using the method developed in [B], we have constructed in [An3] the
submanifolds $\Sigma D_{k}(N,P)$ ($k\geq4$) and $\Sigma E_{k}(N,P)$ ($k=6$,
$7$ or $8$) in $J^{\infty}(N,P)$, which play the similar role for the
singularities $D_{k}$ and $E_{k}$ respectively as $\Sigma^{I}(N,P)$ does for
$A_{k}$. Let $\Omega(N,P)$ denote an open subset of $J^{\infty}(N,P)$, which
consists of all regular jets and a number of prescribed submanifolds $\Sigma
A_{i}(N,P)$, $\Sigma D_{j}(N,P)$ and $\Sigma E_{k}(N,P)$. It is known when
$\Omega(N,P)$ becomes open by the adjacency relations of these singularities
given in [Ar, Corollary 8.7] (Lemma 3.5). It is known that $\Omega(N,P)$ is
the subbundle of $J^{\infty}(N,P)$ with the projection $\pi_{N}^{\infty}%
\times\pi_{P}^{\infty}$, whose fiber is denoted by $\Omega(n,p)$. A smooth map
$f:N\rightarrow P$ is called an $\Omega$-\textit{regular map} if and only if
(i) $j^{\infty}f(N)\subset\Omega(N,P)$ and (ii) $j^{\infty}f$ is transverse to
all $\Sigma A_{i}(N,P)$, $\Sigma D_{j}(N,P)$ and $\Sigma E_{k}(N,P)$ in
$\Omega(N,P)$.

We will study when a given continuous map is homotopic to an $\Omega$-regular
map. Let $C_{\Omega}^{\infty}(N,P)$ denote the space consisting of all
$\Omega$-regular maps equipped with the $C^{\infty}$-topology. Let
$\Gamma_{\Omega}(N,P)$ (resp. $\Gamma_{\Omega}^{tr}(N,P))$ denote the space
consisting of all continuous sections (resp. sections transverse to all
$\Sigma X_{i}(N,P)$ for $X_{i}=A_{i}$, $D_{i}$ and $E_{i}$) of the fiber
bundle $\pi_{N}^{\infty}|\Omega(N,P):\Omega(N,P)\rightarrow N$ equipped with
the compact-open topology. Then there exists a continuous map
\[
j_{\Omega}:C_{\Omega}^{\infty}(N,P)\rightarrow\Gamma_{\Omega}(N,P)
\]
defined by $j_{\Omega}(f)=j^{\infty}f$. It follows from the well-known theorem
due to Gromov[G1] that if $N$ is a connected open manifold, then $j_{\Omega}$
is a weak homotopy equivalence. This property is called the homotopy principle
(the terminology used in [G2]). If $N$ is a closed manifold, then it becomes a
hard problem for us to prove the homotopy principle. As the investigation
preceding [G1], we should refer to the Smale-Hirsch Immersion Theorem ([H1]),
the $k$-mersion Theorem due to [F] and the Phillips Submersion Theorem for
open manifolds ([P]). Du Plessis has introduced the extensibility condition
under which the homotopy principle holds for maps with certain Thom-Boardman
singularities in [duP1], or with prescribed $C^{\infty}$ simple singularities
in [duP2]. \`{E}lia\v{s}berg has proved in [E1] and [E2] the homotopy
principle in the $1$-jet level for fold-maps. In [An9] we have proved the
homotopy principle in the existence level for maps with prescribed
Thom-Boardman singularities in dimensions $n\geq p\geq2$.

We prove the following theorem for closed manifolds.

\begin{theorem}
Let $n\geq p\geq2$. Let $N$ and $P$ be connected manifolds of dimensions $n$
and $p$ respectively with $\partial N=\emptyset$. Assume that $\Omega(N,P)$
contains $\Sigma A_{1}(N,P)$ at least. Let $C$ be a closed subset of $N$. Let
$s$ be a section of $\Gamma_{\Omega}(N,P)$ which has an $\Omega$-regular map
$g$ defined on a neighborhood of $C$ into $P$, where $j^{\infty}g=s$.

Then there exists an $\Omega$-regular map $f:N\rightarrow P$ such that
$j^{\infty}f$ is homotopic to $s$ relative to a neighborhood of $C$ by a
homotopy $s_{\lambda}$ in $\Gamma_{\Omega}(N,P)$ with $s_{0}=s$ and
$s_{1}=j^{\infty}f$.
\end{theorem}

This theorem has been stated in [An2, Theorems 1.1 and 2.1] and has
generalized [An1, Theorem 1]. However, the proof of Theorem 0.1 was not given,
because a more general theorem was possibly supposed to appear. As far as the
author knows, the proof \ of Theorem 0.1 has not been published until now.
Recently, it turns out that this kind of the homotopy principle has many
applications. Theorem 0.1 is very important even for fold-maps in [An6,
Theorem 1] and [An7, Theorems 0.2 and 0.3] and the homotopy type of
$\Omega(n,p)$ associated to fold-maps, which consists of all regular jets and
jets with $A_{1}$-singularities, are determined in [An5] and [An7]. The famous
theorem about the elimination of cusps in [L1] (see also [T]) is now a simple
consequence of these theorems. In [Sady1] Sadykov has applied [An1, Theorem 1]
to the elimination of higher $A_{k}$ singularities ($k\geq3$) for Morin maps
between orientable manifolds when $n-p$ is odd. We can describe the cobordism
group of maps with prescribed singularities of types $A$, $D$ and $E$ in terms
of certain stable homotopy groups by using the result in [An10, Theorem 4.2].

We refer the reader to [Saek] and [Sady2] as other applications of the
homotopy principle.

We give the following application of Theorem 0.1.

\begin{theorem}
Let $n>p\geq2$. If $n-p$ is even, then we have

$(1)$ a smooth map $f:N\rightarrow P$ admitting only singularities of types
$A_{i}$ $(i\geq1)$ and $D_{j}$ $(j\geq4)$ is homotopic to a smooth map
admitting only singularities of types $A_{i}$ $(i\geq1)$, $D_{4}$ and $D_{5}$,

$(2)$ a smooth map $f:N\rightarrow P$ admitting only singularities of types
$A_{i}$ $(i\geq1)$, $D_{j}$ $(j\geq4)$ and $E_{k}$\ $(8\geq k\geq6)$ is
homotopic to a smooth map admitting only singularities of types $A_{i}$
$(i\geq1)$, $D_{4}$, $D_{5}$ and $E_{6}.$
\end{theorem}

In Section 1 we explain notations which are used in this paper. In Sections 2
and 3\ we review the definitions and the fundamental properties of the
Boardman submanifolds, $\Sigma D_{k}(N,P)$ and $\Sigma E_{k}(N,P)$
respectively. In Section 4 we give Theorem 4.1 which is a simpler version of
Theorem 0.1. We will prove Theorem 0.1 by the induction, and prepare a certain
rotation of the tangent spaces defined around the singularities in $N$ to
deform $s$. In Section 5 we prepare several lemmas which are necessary in the
deformation of $s$. In Section 6 we prove Theorem 4.1. In Section 7 we prove
Theorem 0.2 by applying Theorem 0.1.

The author would like to thank the referee for his kind and helpful comments.

\section{Notations}

Throughout the paper all manifolds are Hausdorff, paracompact and smooth of
class $C^{\infty}$. Maps are basically continuous, but may be smooth (of class
$C^{\infty}$) if necessary. Given a fiber bundle $\pi:E\rightarrow X$ and a
subset $C$ in $X,$ we denote $\pi^{-1}(C)$ by $E|_{C}.$ Let $\pi^{\prime
}:F\rightarrow Y$ be another fiber bundle. A map $\tilde{b}:E\rightarrow F$ is
called a fiber map over a map $b:X\rightarrow Y$ if $\pi^{\prime}\circ
\tilde{b}=b\circ\pi$ holds. The restriction $\tilde{b}|(E|_{C}):E|_{C}%
\rightarrow F$ (or $F|_{b(C)}$) is denoted by $\tilde{b}|_{C}$. In particular,
for a point $x\in X,$ $E|_{x}$ and $\tilde{b}|_{x}$ are simply denoted by
$E_{x}$ and $\tilde{b}_{x}:E_{x}\rightarrow F_{b(x)}$ respectively. We denote,
by $b^{F}$, the induced fiber map $b^{\ast}(F)\rightarrow F$ covering $b$. For
a map $j:W\rightarrow X$, let $(b\circ j)^{\ast}(\tilde{b}):j^{\ast
}E\rightarrow(b\circ j)^{\ast}F$ over $W$ be the fiber map canonically induced
from $b$ and $j$. A fiberwise homomorphism $E\rightarrow F$ is simply called a
homomorphism. For a vector bundle $E$ with metric and a positive function
$\delta$ on $X$, let $D_{\delta}(E)$ be the associated disk bundle of $E$ with
radius $\delta$. If there is a canonical isomorphism between two vector
bundles $E$ and $F$ over $X=Y,$ then we write $E\cong F$.

When $E$ and $F$ are vector bundles over $X=Y$, Hom$(E,F)$ denotes the vector
bundle over $X$ with fiber Hom$(E_{x},F_{x})$, $x\in X$, which consists of all
homomorphisms $E_{x}\rightarrow F_{x}$.

Let $J^{k}(N,P)$ denote the $k$-jet space of manifolds $N$ and $P$ ($k$ may be
$\infty$). Let $\pi_{N}^{k}$ and $\pi_{P}^{k}$ be the projections mapping a
jet to its source and target respectively. The map $\pi_{N}^{k}\times\pi
_{P}^{k}:J^{k}(N,P)\rightarrow N\times P$ induces a structure of a fiber
bundle with structure group $L^{k}(p)\times L^{k}(n)$, where $L^{k}(m)$
denotes the group of all $k$-jets of local diffeomorphisms of $(\mathbf{R}%
^{m},0)$. The fiber $(\pi_{N}^{k}\times\pi_{P}^{k})^{-1}(x,y)$ is denoted by
$J_{x,y}^{k}(N,P)$.

Let $\pi_{N}$ and $\pi_{P}$ be the projections of $N\times P$ onto $N$ and $P$
respectively.\ We set
\begin{equation}
J^{k}(TN,TP)=\text{\textrm{Hom}}(\pi_{N}^{\ast}(TN)\oplus S^{2}(\pi_{N}^{\ast
}(TN))\oplus\cdots\oplus S^{k}(\pi_{N}^{\ast}(TN)),\pi_{P}^{\ast}(TP))
\end{equation}
over $N\times P$. Here, for a vector bundle $E$ over $X$, let $S^{i}(E)$ be
the vector bundle $\cup_{x\in X}S^{i}(E_{x})$ over $X$, where $S^{i}(E_{x})$
denotes the $i$-fold symmetric product of $E_{x}$. If we provide $N$ and $P$
with Riemannian metrics, then the Levi-Civita connections induce the
exponential maps $\exp_{N,x}:T_{x}N\rightarrow N$ and $\exp_{P,y}%
:T_{y}P\rightarrow P$. In dealing with the exponential maps we always consider
the convex neighborhoods ([K-N]). We define the smooth bundle map
\begin{equation}
J^{k}(N,P)\mathbf{\rightarrow}J^{k}(TN,TP)\text{ \ \ \ over }N\times P
\end{equation}
by sending $z=j_{x}^{k}f\in J_{x,y}^{k}(N,P)$ to the $k$-jet of $(\exp
_{P,y})^{-1}\circ f\circ\exp_{N,x}$ at $\mathbf{0}\in T_{x}N$, which is
regarded as an element of $J^{k}(T_{x}N,T_{y}P)(=J_{x,y}^{k}(TN,TP))$. The
structure group of $J^{k}(TN,TP)$ is reduced to $O(p)\times O(n)$.

Recall that $S^{i}(E)\ $has the inclusion $S^{i}(E)\rightarrow\otimes^{i}E$
and the canonical projection $\otimes^{i}E\rightarrow S^{i}(E)$ (see [B,
Section 4] and [Mats, Ch. III, Section 2]). Let $E_{j}$ be subbundles of $E$
$(j=1,\cdots,r)$. We define $E_{1}\bigcirc\cdots\bigcirc E_{i}=\bigcirc
_{j=1}^{i}E_{j}$ to be the image of $E_{1}\otimes\cdots\otimes E_{i}%
=\otimes_{j=1}^{i}E_{j}\rightarrow\otimes^{i}E\rightarrow S^{i}(E)$. When
$E_{j+1}=\cdots=E_{j+\ell}$, we often write $E_{1}\bigcirc\cdots\bigcirc
E_{j}\bigcirc^{\ell}E_{j+1}\bigcirc E_{j+\ell+1}\bigcirc\cdots\bigcirc E_{i}$
in place of $\bigcirc_{j=1}^{i}E_{j}.$

\section{Boardman manifolds}

We review well-known results about Boardman manifolds in $J^{\infty}%
(N,P)$\ ([B] and [L2]). Let $I=(i_{1},i_{2},\cdots,i_{k},\cdots)$ be a
Boardman symbol, which is a sequence of nonnegative integers with $i_{1}\geq
i_{2}\geq\cdots\geq i_{k}\geq\cdots$. Set $I_{k}=(i_{1},i_{2},\cdots,i_{k})$
and $(I_{k},0)=(i_{1},i_{2},\cdots,i_{k},0)$. In the infinite jet space
$J^{\infty}(N,P)$, there have been defined a sequence of the submanifolds
$\Sigma^{I_{1}}(N,P)\supseteq\cdots\supseteq\Sigma^{I_{k}}(N,P)\supseteq
\cdots$ with the following properties. In this paper we often write
$\Sigma^{I_{k}}$ for $\Sigma^{I_{k}}(N,P)$ if there is no confusion.

Let $\mathbf{P}=(\pi_{P}^{\infty})^{\ast}(TP)$ and $\mathbf{D}$ be the total
tangent bundle defined over $J^{\infty}(N,P)$.\ We explain an important
property of the total tangent bundle $\mathbf{D}$, which is often used in this
paper. Let $f:(N,x)\rightarrow(P,y)$ be a germ and $\digamma$ be a smooth
function in the sense of [B, Definition 1.4] defined on a neighborhood of
$j_{x}^{\infty}f$. Given a vector field $v$ defined on a neighborhood of $x$
in $N$, there is a total vector field $D$ defined on a neighborhood of
$j_{x}^{\infty}f$ such that $D\digamma\circ j^{\infty}f=v(\digamma\circ
j^{\infty}f)$. It follows that $d(j^{\infty}f)(v)(\digamma)=D\digamma
(j^{\infty}f)$ for $d(j^{\infty}f):TN\rightarrow T(J^{\infty}(N,P))$ around
$x$. This implies $d(j^{\infty}f)(v)=D.$ Hence, we have $\mathbf{D\cong(}%
\pi_{N}^{\infty})^{\ast}(TN)$.

First we have the first derivative $\mathbf{d}_{1}:\mathbf{D}\rightarrow
\mathbf{P}$\ over $J^{\infty}(N,P)$. We define $\Sigma^{I_{1}}(N,P)$ to be the
submanifold of $J^{\infty}(N,P)$ which consists of all jets $z$ such that the
kernel rank of $\mathbf{d}_{1,z}$ is $i_{1}$. Since $\mathbf{d}_{1}%
|_{\Sigma^{I_{1}}(N,P)}$ is of constant rank $n-i_{1}$, we set $\mathbf{K}%
_{1}=$Ker$(\mathbf{d}_{1})$ and $\mathbf{Q}_{1}=$Cok$(\mathbf{d}_{1})$, which
are vector bundles over $\Sigma^{I_{1}}(N,P)$. Set $\mathbf{K}_{0}=\mathbf{D}%
$, $\mathbf{P}_{0}=\mathbf{P}$ and $\Sigma^{I_{0}}(N,P)=J^{\infty}(N,P)$. We
can inductively define $\Sigma^{I_{k}}(N,P)$ and the bundles $\mathbf{K}_{k}$
and $\mathbf{P}_{k}$ over $\Sigma^{I_{k}}(N,P)$ ($k\geq1$) with the properties:

(1) $\mathbf{K}_{k-1}|_{\Sigma^{I_{k}}(N,P)}\supseteq\mathbf{K}_{k}$ over
$\Sigma^{I_{k}}(N,P)$.

(2) $\mathbf{K}_{k}$ is an $i_{k}$-dimensional subbundle of $T(\Sigma
^{I_{k-1}}(N,P))|_{\Sigma^{I_{k}}(N,P)}.$

(3) There exists the $(k+1)$-th intrinsic derivative $\mathbf{d}%
_{k+1}:T(\Sigma^{I_{k-1}}(N,P))|_{\Sigma^{I_{k}}(N,P)}\rightarrow
\mathbf{P}_{k}$ over $\Sigma^{I_{k}}(N,P),$ so that it induces the exact
sequence%
\[%
\begin{array}
[c]{l}%
\mathbf{0\rightarrow}T(\Sigma^{I_{k}}(N,P))\overset{\text{inclusion}%
}{\hookrightarrow}T(\Sigma^{I_{k-1}}(N,P))|_{\Sigma^{I_{k}}(N,P)}%
\overset{\mathbf{d}_{k+1}}{\longrightarrow}\mathbf{P}_{k}\text{ }%
\rightarrow\mathbf{0}\text{ \ \ over }\Sigma^{I_{k}}(N,P)\mathbf{.}%
\end{array}
\]
Namely, $\mathbf{d}_{k+1}$ induces the isomorphism of the normal bundle%
\[
\nu(I_{k}\subset I_{k-1})=(T(\Sigma^{I_{k-1}}(N,P))|_{\Sigma^{I_{k}}%
(N,P)})/T(\Sigma^{I_{k}}(N,P))
\]
of $\Sigma^{I_{k}}(N,P)$ in $\Sigma^{I_{k-1}}(N,P)$ onto $\mathbf{P}_{k}$.

(4) $\Sigma^{I_{k+1}}(N,P)$ is defined to be the submanifold of $\Sigma
^{I_{k}}(N,P)$, which consists of all jets $z$ with dim(Ker$(\mathbf{d}%
_{k+1,z}|\mathbf{K}_{k,z}))=i_{k+1}.$

(5) Set $\mathbf{K}_{k+1}=$Ker$(\mathbf{d}_{k+1}|\mathbf{K}_{k})$ and
$\mathbf{Q}_{k+1}\mathbf{=}$Cok$(\mathbf{d}_{k+1}|\mathbf{K}_{k})$ over
$\Sigma^{I_{k+1}}(N,P)$. Then it follows that $(\mathbf{K}_{k}|_{\Sigma
^{I_{k+1}}(N,P)})\cap T(\Sigma^{I_{k}}(N,P))|_{\Sigma^{I_{k+1}}(N,P)}%
=\mathbf{K}_{k+1}.$

(6) The intrinsic derivative%
\[
d(\mathbf{d}_{k+1}|\mathbf{K}_{k}):T(\Sigma^{I_{k}}(N,P))|_{\Sigma^{I_{k+1}%
}(N,P)}\rightarrow\text{Hom}(\mathbf{K}_{k+1},\mathbf{Q}_{k+1})\text{ \ \ over
}\Sigma^{I_{k+1}}(N,P)
\]
of $\mathbf{d}_{k+1}|\mathbf{K}_{k}$ is of constant rank $\dim(\Sigma^{I_{k}%
}(N,P))-\dim(\Sigma^{I_{k+1}}(N,P))$. We set $\mathbf{P}_{k+1}%
=\operatorname{Im}(d(\mathbf{d}_{k+1}|\mathbf{K}_{k}))$ and define
$\mathbf{d}_{k+2}$ to be%
\[
\mathbf{d}_{k+2}=d(\mathbf{d}_{k+1}|\mathbf{K}_{k}):T(\Sigma^{I_{k}%
}(N,P))|_{\Sigma^{I_{k+1}}(N,P)}\rightarrow\mathbf{P}_{k+1}%
\]
as the epimorphism.

(7) The codimension of $\Sigma^{I_{k}}(N,P)$ in $J^{k}(N,P)$ is described in
[B, Theorem 6.1].

(8) We define $\Sigma^{I_{k},0}(N,P)=\Sigma^{I_{k}}(N,P)\setminus
\Sigma^{I_{k+1}}(N,P).$

(9) The submanifold $\Sigma^{I_{k}}(N,P)$ is actually defined so that it
coincides with the inverse image of the submanifold in $J^{k}(N,P)$ by
$\pi_{k}^{\infty}$.

We apply the above construction to the two symbols $J=(n-p+1,1,\cdots
,1,\cdots)$ and $\frak{J}=(n-p+1,2,1,\cdots,1,\cdots)$. For both symbols, we
write $\mathbf{K}$ and $\mathbf{Q}$ for $\mathbf{K}_{1}$ and $\mathbf{Q}_{1}$
respectively in the rest of the paper. We have dim$\mathbf{K}=n-p+1$,
dim$\mathbf{Q}=1$ and
\begin{equation}
\mathbf{d}_{2}:T(J^{\infty}(N,P))|_{\Sigma^{n-p+1}(N,P)}\rightarrow
\mathbf{P}_{1}=\text{Hom}(\mathbf{K},\mathbf{Q)}\text{ \ \ over }%
\Sigma^{n-p+1}(N,P).
\end{equation}

\noindent When $k\geq2$, we usually write $\mathbf{K}_{k}^{I}$, $\mathbf{Q}%
_{k}^{I}$, $\mathbf{P}_{k}^{I}$ and $\mathbf{d}_{k}^{I}$ in place of
$\mathbf{K}_{k}$, $\mathbf{Q}_{k}$, $\mathbf{P}_{k}$ and $\mathbf{d}_{k}$
respectively for $I=J$ and $\frak{J}$.

We first deal with the symbol $J$. In this paper $\Sigma^{J_{k}}(N,P)$ and
$\Sigma^{J_{k},0}(N,P)$ are often denoted by $\Sigma{{\overline{A}_{k}}}(N,P)$
and $\Sigma{A_{k}}(N,P)$ respectively.

($J$-1) We have that $\dim\mathbf{K}_{2}^{J}=1$, $\mathbf{Q}_{k}^{J}%
=$Hom$(\bigcirc^{k-1}\mathbf{K}_{2}^{J},\mathbf{Q})|_{\Sigma{{\overline{A}%
_{k}}}}$ and the exact sequence
\begin{equation}
\mathbf{0\rightarrow}T(\Sigma{{\overline{A}_{k}}}(N,P))\hookrightarrow
T(\Sigma{{\overline{A}_{k-1}}}(N,P))|_{\Sigma{{\overline{A}_{k}}}}%
\overset{\mathbf{d}_{k+1}^{J}}{\longrightarrow}\mathbf{P}_{k}^{J}%
\cong\text{Hom}(\bigcirc^{k}\mathbf{K}_{2}^{J},\mathbf{Q})|_{\Sigma
{{\overline{A}_{k}}}}\rightarrow\mathbf{0}.
\end{equation}
Set $\nu({\overline{A}_{k}}\subset{\overline{A}_{k-1}})=\nu(J_{k}\subset
J_{k-1})$. Then $\mathbf{d}_{k+1}^{J}$induces the canonical isomorphism
$\nu({\overline{A}_{k}}\subset{\overline{A}_{k-1}})\rightarrow\mathbf{P}%
_{k}^{J}$.

($J$-2) $\Sigma^{J_{k}}(N,P)$ is of codimension $n-p+k$ in $J^{\infty}(N,P).$

($J$-3) A smooth map germ $f:(N,x)\rightarrow(P,y)$ has $x$ as a singularity
of type $A_{k}$ if and only if $j^{\infty}f$ is transverse to $\Sigma
^{J_{k},0}(N,P)$ at $x$ and $j_{x}^{\infty}f\in\Sigma^{J_{k},0}(N,P)$ (see [Mo]).

We next turn to the symbol $\frak{J}=(n-p+1,2,1,\cdots,1,\cdots)$. Then we
have the following by using [B, Corollary 7.10].

($\frak{J}$-1) We have that $\dim\mathbf{K}_{2}^{\frak{J}}=2$, $\mathbf{Q}%
_{2}^{\frak{J}}=\mathrm{Hom}(\mathbf{K}_{2}^{\frak{J}},\mathbf{Q}%
)|_{\Sigma^{\frak{J}_{2}}}$,%
\begin{equation}
\mathbf{d}_{3}^{\frak{J}}:T(\Sigma^{n-p+1}(N,P))|_{\Sigma^{\frak{J}_{2}}%
}\rightarrow\mathbf{P}_{2}^{\frak{J}}\cong\text{Hom}(\bigcirc^{2}%
\mathbf{K}_{2}^{\frak{J}},\mathbf{Q)}|_{\Sigma^{\frak{J}_{2}}},
\end{equation}
and the exact sequence%
\begin{equation}
\mathbf{0\rightarrow K}_{3}^{\frak{J}}\rightarrow\mathbf{K}_{2}^{\frak{J}%
}\overset{\mathbf{d}_{3}^{\frak{J}}|\mathbf{K}_{2}^{\frak{J}}}{\longrightarrow
}\text{Hom}(\bigcirc^{2}\mathbf{K}_{2}^{\frak{J}},\mathbf{Q})\rightarrow
\text{Hom}(\mathbf{K}_{3}^{\frak{J}}\bigcirc\mathbf{K}_{2}^{\frak{J}%
},\mathbf{Q)\rightarrow0}\text{ \ over }\Sigma^{\frak{J}_{3}}(N,P).
\end{equation}
Here, $\mathbf{K}_{3}^{\frak{J}}$ $(=$Ker$(\mathbf{d}_{3}^{\frak{J}%
}|\mathbf{K}_{2}^{\frak{J}}))$ and Im$(\mathbf{d}_{3}^{\frak{J}}%
|\mathbf{K}_{2}^{\frak{J}})\cong\mathrm{Hom}(\mathbf{K}_{2}^{\frak{J}%
}/\mathbf{K}_{3}^{\frak{J}}\bigcirc\mathbf{K}_{2}^{\frak{J}}/\mathbf{K}%
_{3}^{\frak{J}},\mathbf{Q})|_{\Sigma^{\frak{J}_{3}}}$ are of dimension
$1$\ and $\mathbf{Q}_{3}^{\frak{J}}\cong\mathrm{Hom}(\mathbf{K}_{3}^{\frak{J}%
}\bigcirc\mathbf{K}_{2}^{\frak{J}},\mathbf{Q})|_{\Sigma^{\frak{J}_{3}}}$ with
dim$\mathbf{Q}_{3}^{\frak{J}}=2\mathbf{.}$

($\frak{J}$E-1) Since $\mathbf{P}_{3}^{\frak{J}}=\mathrm{Hom}(\mathbf{K}%
_{3}^{\frak{J}},\mathbf{Q}_{3}^{\frak{J}})\cong\mathrm{Hom}((\bigcirc
^{2}\mathbf{K}_{3}^{\frak{J}})\bigcirc\mathbf{K}_{2}^{\frak{J}},\mathbf{Q}%
)|_{\Sigma^{\frak{J}_{3}}}$, we have the epimorphism%
\begin{equation}
\mathbf{d}_{4}^{\frak{J}}:T(\Sigma^{\frak{J}_{2}}(N,P))|_{\Sigma^{\frak{J}%
_{3}}}\rightarrow\text{Hom}((\bigcirc^{2}\mathbf{K}_{3}^{\frak{J}}%
)\bigcirc\mathbf{K}_{2}^{\frak{J}},\mathbf{Q})|_{\Sigma^{\frak{J}_{3}}}.
\end{equation}

($\frak{J}$E-2) By $\mathbf{K}_{4}^{\frak{J}}=\mathbf{K}_{3}^{\frak{J}%
}|_{\Sigma^{\frak{J}_{4}}}$, $\mathbf{Q}_{4}^{\frak{J}}\cong\mathrm{Hom}%
((\bigcirc^{2}\mathbf{K}_{4}^{\frak{J}})\bigcirc\mathbf{K}_{2}^{\frak{J}%
},\mathbf{Q})|_{\Sigma^{\frak{J}_{4}}}$ and $\mathbf{P}_{4}^{\frak{J}%
}=\mathrm{Hom}(\mathbf{K}_{4}^{\frak{J}},\mathbf{Q}_{4}^{\frak{J}}%
)\cong\mathrm{Hom}((\bigcirc^{3}\mathbf{K}_{4}^{\frak{J}})\bigcirc
\mathbf{K}_{2}^{\frak{J}},\mathbf{Q})|_{\Sigma^{\frak{J}_{4}}}$, we have the
epimorphism%
\begin{equation}
\mathbf{d}_{5}^{\frak{J}}:T(\Sigma^{\frak{J}_{3}}(N,P))|_{\Sigma^{\frak{J}%
_{4}}}\rightarrow\text{Hom}((\bigcirc^{3}\mathbf{K}_{4}^{\frak{J}}%
)\bigcirc\mathbf{K}_{2}^{\frak{J}},\mathbf{Q})|_{\Sigma^{\frak{J}_{4}}}.
\end{equation}

($\frak{J}$E-3) The homomorphisms $\mathbf{d}_{3}^{\frak{J}}$, $\mathbf{d}%
_{4}^{\frak{J}}$ and $\mathbf{d}_{5}^{\frak{J}}$ induce the homomorphisms%
\begin{equation}%
\begin{array}
[c]{ll}%
\widetilde{\mathbf{d}}_{3}^{\frak{J}}:\bigcirc^{3}\mathbf{K}_{2}^{\frak{J}%
}\rightarrow\mathbf{Q} & \text{over }\Sigma^{\frak{J}_{2}}(N,P),\\
\widetilde{\mathbf{d}}_{4}^{\frak{J}}:(\bigcirc^{3}\mathbf{K}_{3}^{\frak{J}%
})\bigcirc\mathbf{K}_{2}^{\frak{J}}\rightarrow\mathbf{Q} & \text{over }%
\Sigma^{\frak{J}_{3}}(N,P),\\
\widetilde{\mathbf{d}}_{5}^{\frak{J}}:(\bigcirc^{4}\mathbf{K}_{4}^{\frak{J}%
})\bigcirc\mathbf{K}_{2}^{\frak{J}}\rightarrow\mathbf{Q} & \text{over }%
\Sigma^{\frak{J}_{4}}(N,P)
\end{array}
\end{equation}
respectively.

\section{Submanifolds $\Sigma D_{k}(N,P)$ and $\Sigma E_{k}(N,P)$}

In this section we review the definition and the properties of the manifolds
$\Sigma{X{_{k}}}(N,P)$ and $\Sigma\overline{X}{{_{k}}}(N,P)$ for $X_{k}=D_{k}$
and $E_{k}$\ in [An3]. We identify, as usual, Hom$(\bigcirc^{3}\mathbf{R}%
^{2},\mathbf{R})$ with the set of all cubic forms with variables $u$ and $v$
on $\mathbf{R}^{2}$. By [Ar, Lemma 5.1] it is decomposed into the five
manifolds $S_{4}^{\pm}$, $S_{5}$, $S_{E}$ and $\mathbf{0}$ which are orbit
manifolds by $GL(2)$ through $u^{2}\ell\pm\ell^{3}$, $u^{2}\ell$, $u^{3}$ and
$0$ respectively. Let us recall the bundle Hom$(\bigcirc^{3}\mathbf{K}%
_{2}^{\frak{J}},\mathbf{Q)}$ over $\Sigma^{\frak{J}_{2}}(N,P)$. Then we obtain
the subbundles $S_{4}^{\pm}(\bigcirc^{3}\mathbf{K}_{2}^{\frak{J}},\mathbf{Q)}%
$, $S_{5}(\bigcirc^{3}\mathbf{K}_{2}^{\frak{J}},\mathbf{Q)}$, $S_{E}%
(\bigcirc^{3}\mathbf{K}_{2}^{\frak{J}},\mathbf{Q)}$ of Hom$(\bigcirc
^{3}\mathbf{K}_{2}^{\frak{J}},\mathbf{Q)}$ associated to $S_{4}^{\pm}$,
$S_{5}$, $S_{E}$ respectively. Let $\nu({{{{\overline{X}_{k+1}\subset
\overline{X}_{k})}}}}$ be the normal bundle $(T(\Sigma{\overline{X}_{k}%
}(N,P))|_{\Sigma{\overline{X}_{k+1}}})/T(\Sigma{\overline{X}_{k+1}}(N,P))$ and
set $\nu(X_{k+1}\subset\overline{X}_{k})=\nu(\overline{X}_{k+1}\subset
\overline{X}_{k})|_{\Sigma{X_{k+1}}}$.

Let us consider $\widetilde{\mathbf{d}}_{3}^{\frak{J}}|_{\Sigma^{\frak{J}%
_{2},0}}:\bigcirc^{3}\mathbf{K}_{2}^{\frak{J}}\rightarrow\mathbf{Q}$ over
$\Sigma^{\frak{J}_{2},0}(N,P)$. We define the\ submanifolds%
\begin{align*}
\Sigma D_{4}^{\pm}(N,P)  &  =\left(  \widetilde{\mathbf{d}}_{3}^{\frak{J}%
}|_{\Sigma^{\frak{J}_{2},0}}\right)  ^{-1}(S_{4}^{\pm}(\bigcirc^{3}%
\mathbf{K}_{2}^{\frak{J}},\mathbf{Q)),}\\
\Sigma{{\overline{D}_{5}}}(N,P)  &  =\left(  \widetilde{\mathbf{d}}%
_{3}^{\frak{J}}|_{\Sigma^{\frak{J}_{2},0}}\right)  ^{-1}(S_{5}(\bigcirc
^{3}\mathbf{K}_{2}^{\frak{J}},\mathbf{Q))}%
\end{align*}
in $\Sigma^{\frak{J}_{2},0}(N,P)$ respectively. Its transversality follows
from [An3, Proposition 3.5].

(D-i) We set $\Sigma D_{4}(N,P)=\Sigma{{D_{4}^{+}}}(N,P)\cup\Sigma{{D_{4}^{-}%
}}(N,P)$, and $\Sigma{{\overline{D}_{4}}}(N,P)=\Sigma^{\frak{J}_{2},0}(N,P).$
Since $\mathbf{d}_{3}^{\frak{J}}|\mathbf{K}_{2}^{\frak{J}}:\mathbf{K}%
_{2}^{\frak{J}}\rightarrow\mathrm{Hom}(\bigcirc^{2}\mathbf{K}_{2}^{\frak{J}%
},\mathbf{Q)}$ over $\Sigma^{\frak{J}_{2},0}(N,P)$ is injective, it follows
that $\mathbf{K}_{2}^{\frak{J}}=(T(\Sigma^{\frak{J}_{1}}(N,P))\cap
\mathbf{K)|}_{\Sigma^{\frak{J}_{2}}(N,P)}$ and $\mathbf{K}_{2}^{\frak{J}}\cap
T(\Sigma^{\frak{J}_{2},0}(N,P))=\{\mathbf{0}\}$.

(D-ii) There exist a small neighborhood $U(\Sigma{{\overline{D}_{5})}}$ of
$\Sigma{{\overline{D}_{5}}}(N,P)$ in $\Sigma^{\frak{J}_{2},0}(N,P)$ and the
line bundle $\mathbf{L}$\textbf{ }defined over $\Sigma{{\overline{D}_{5}}%
}(N,P)$, which is uniquely extended to the subbundle $\widetilde{\mathbf{L}}$
over $U(\Sigma{{\overline{D}_{5})}}$ of $\mathbf{K}_{2}^{\frak{J}}$ such that
for any $z\in\Sigma{{\overline{D}_{5}}}(N,P)$, $\mathbf{L}_{z}$ coincides with
$\left(  \mathbf{d}_{3,z}^{\frak{J}}|\mathbf{K}_{2,z}^{\frak{J}}\right)
^{-1}(H_{z})$, where $H_{z}$ is the subset of Hom$(\bigcirc^{2}\mathbf{K}%
_{2,z}^{\frak{J}},\mathbf{Q}_{z}\mathbf{)}$ of all quadratic forms of rank $1$
or $0$. Furthermore, we have Cok$(\mathbf{d}_{3}^{\frak{J}}|\mathbf{K}%
_{2}^{\frak{J}})=\mathrm{Hom}(\bigcirc^{2}\mathbf{L},\mathbf{Q)}$ over
$\Sigma\overline{D}_{5}(N,P).$

(D-iii) Let $z\in U(\Sigma{{\overline{D}_{5})}}$. Then $z$ lies in
$\Sigma{{\overline{D}_{5}(N,P)}}$ if and only if the restriction
$\widetilde{\mathbf{d}}_{3,z}^{\frak{J}}|\bigcirc^{3}\widetilde{\mathbf{L}%
}_{z}$ is a null homomorphism.

(D-iv) We successively construct the submanifolds $\Sigma{{\overline{D}_{k}}%
}(N,P)$ ($k\geq5$) with $\Sigma{\overline{D}_{4}}(N,P)\linebreak \supset
U(\Sigma{\overline{D}_{5}})\supset\Sigma{\overline{D}_{5}}(N,P)\supset
\cdots\supset\Sigma{\overline{D}_{k}}(N,P)\supset\cdots$ and the homomorphisms%
\begin{equation}
\mathbf{r}_{3}:\bigcirc^{3}\widetilde{\mathbf{L}}\mathbf{\rightarrow Q}\text{
\ over }U({{{{\overline{D}_{5}}})}}\text{ \ and \ }\mathbf{r}_{k-1}%
:\bigcirc^{k-1}\mathbf{L\rightarrow Q}\text{ \ \ over }{{\Sigma{{\overline
{D}_{k}}}(N,P)}}\text{ (}k\geq5\text{)}%
\end{equation}
with the following properties:

(D-iv-1) $\widetilde{\mathbf{d}}_{3}^{\frak{J}}|\bigcirc^{3}\widetilde
{\mathbf{L}}\mathbf{=r}_{3}.$

(D-iv-2) ${{\Sigma{{\overline{D}_{k}}}(N,P)}}$ is of codimension $n-p+k$ in
$J^{\infty}(N,P).$

(D-iv-3) An element $z\in U({{\Sigma{{\overline{D}_{5}}})}}$ lies in
${{\Sigma{{\overline{D}_{5}}}(N,P)}}$ if and only if $\mathbf{r}_{3,z}$
vanishes. An element $z\in{{\Sigma{{\overline{D}_{k}}}(N,P)}}$ lies in
${{\Sigma{{\overline{D}_{k+1}}}(N,P)}}$ if and only if $\mathbf{r}_{k-1,z}$ vanishes.

(D-iv-4) The intrinsic derivative of $\mathbf{r}_{k-1}$%
\begin{equation}
d(\mathbf{r}_{k-1}):T({{\Sigma{{\overline{D}_{k}}}(N,P))|}}_{{{\Sigma
{{\overline{D}_{k+1}}}}}}\rightarrow\text{Hom}(\bigcirc^{k-1}\mathbf{L}%
,\mathbf{Q)}{{|}}_{{{\Sigma{{\overline{D}_{k+1}}}}}}%
\end{equation}
is surjective. In other words, $\mathbf{r}_{k-1}$ is, as a section, transverse
to the zero-section of Hom$(\bigcirc^{k-1}\mathbf{L},\mathbf{Q)}|
_{\Sigma{\overline{D}_{k}}}$ on $\Sigma{\overline{D}_{k+1}}(N,P)$. Then,
$d(\mathbf{r}_{k-1})$ induces the isomorphism $\nu(\overline{D}_{k+1}%
\subset\overline{D}_{k})\rightarrow\mathrm{Hom}(\bigcirc^{k-1}\mathbf{L}%
,\mathbf{Q)}|_{\Sigma{\overline{D}_{k+1}}}.$

(D-v) We define ${{\Sigma{{D_{k}}}(N,P)=\Sigma{{\overline{D}_{k}}%
}(N,P)\setminus\Sigma{{\overline{D}_{k+1}}}(N,P)}}$. Then a smooth map germ
$f:(N,x)\rightarrow(P,y)$ has $x$ as a singularity of type $D_{k}$ if and only
if $j^{\infty}f$ is transverse to ${{\Sigma{{D_{k}}}(N,P)}}$ at $x$\ and
$j_{x}^{\infty}f\in{{\Sigma{{D_{k}}}(N,P).}}$

Next we turn to define ${{\Sigma}}E{{{{_{k}}}(N,P)}}$. We note that the
singularities $E{{{{_{6}}}}}$, $E{{{{_{7}}}}}$ and $E{{{{_{8}}}}}$ have the
Boardman symbols $(\frak{J}_{3},0)$, $(\frak{J}_{3},0)$ and $(\frak{J}_{4},0)$
respectively (the method for the calculation of Boardman symbols in [Math2]
may be convenient).

(E-i) We define $\Sigma E_{6}(N,P)$ to be the open submanifold of
$\Sigma^{\frak{J}_{3},0}(N,P)$ which consists of all jets $z$ such that
$\widetilde{\mathbf{d}}_{4,z}^{\frak{J}}|\bigcirc^{4}\mathbf{K}_{3,z}%
^{\frak{J}}$ does not vanish. We set $\Sigma E_{7}(N,P) =\Sigma^{\frak{J}%
_{3},0}(N,P)\setminus\Sigma E_{6}(N,P)$. Namely, for $z\in{{\Sigma}}E{{{{_{7}%
}}(N,P)}}$, $\widetilde{\mathbf{d}}_{4,z}^{\frak{J}}|\bigcirc^{3}%
\mathbf{K}_{3,z}^{\frak{J}}\bigcirc(\mathbf{K}_{2,z}^{\frak{J}}/\mathbf{K}%
_{3,z}^{\frak{J}})$ does not vanish.

(E-ii) Define ${{\Sigma{{\overline{E}_{7}}}(N,P)}}$ in ${{\Sigma}}%
^{\frak{J}_{3}}{{(N,P)}}$\ to be the inverse image of the zero-section of
$\widetilde{\mathbf{d}}_{4}^{\frak{J}}|\bigcirc^{4}\mathbf{K}_{3}^{\frak{J}}$
over ${{\Sigma}}^{\frak{J}_{3}}{{(N,P)}}$, which is regarded as the section
over ${{\Sigma}}^{\frak{J}_{3}}{{(N,P)}}$. By Lemma 3.1 (i)\ below,\ ${{\Sigma
{{\overline{E}_{7}}}(N,P)}}$ is a submanifold of ${{\Sigma}}^{\frak{J}_{3}%
}{{(N,P)}}$ and the intrinsic derivative $d(\widetilde{\mathbf{d}}%
_{4}^{\frak{J}}|\bigcirc^{4}\mathbf{K}_{3}^{\frak{J}}):T({{\Sigma}}%
^{\frak{J}_{3}}{{(N,P))|}}_{{{\Sigma{{\overline{E}_{7}}}}}}\rightarrow
\mathrm{Hom}(\bigcirc^{4}\mathbf{K}_{3}^{\frak{J}}\mathbf{,Q)}{{|}}%
_{{{\Sigma{{\overline{E}_{7}}}}}}$ of $\widetilde{\mathbf{d}}_{4}^{\frak{J}%
}|\bigcirc^{4}\mathbf{K}_{3}^{\frak{J}}$ is an epimorphism. Hence, we have the
exact sequence%
\begin{equation}
0\rightarrow T({{\Sigma{{\overline{E}_{7}(N,P))\rightarrow}}}}T({{\Sigma}%
}^{\frak{J}_{3}}{{(N,P))|}}_{{{\Sigma{{\overline{E}_{7}}}}}}\overset
{\underrightarrow{d(\widetilde{\mathbf{d}}_{4}^{\frak{J}}|\bigcirc
^{4}\mathbf{K}_{3}^{\frak{J}})~}}{}\text{Hom}(\bigcirc^{4}\mathbf{K}%
_{3}^{\frak{J}}\mathbf{,Q)}{{|}}_{{{\Sigma{{\overline{E}_{7}}}}}}\rightarrow0.
\end{equation}
Let $z\in{{\Sigma}}^{\frak{J}_{4}}{{(N,P)}}$. Since $\mathbf{K}_{4,z}%
^{\frak{J}}=\mathbf{K}_{3,z}^{\frak{J}}$, we have that $\mathbf{d}%
_{4,z}^{\frak{J}}{{{{|}}}}\mathbf{K}_{4,z}^{\frak{J}}$ vanishes. Namely,
$\widetilde{\mathbf{d}}_{4,z}^{\frak{J}}|\bigcirc^{4}\mathbf{K}_{3,z}%
^{\frak{J}}$ vanishes. This implies that ${{\Sigma}}^{\frak{J}_{4}%
}{{(N,P)\subset\Sigma{{\overline{E}_{7}}}(N,P)}}$.

(E-iii) We define ${{\Sigma}}E{{{{_{8}}}(N,P)}}$ to be the submanifold of
${{\Sigma{{\overline{E}_{7}}}(N,P)}}$ which consists of all jets $z$ such that

(i) $\widetilde{\mathbf{d}}_{4,z}^{\frak{J}}|\bigcirc^{3}\mathbf{K}%
_{4,z}^{\frak{J}}\bigcirc(\mathbf{K}_{2,z}^{\frak{J}}/\mathbf{K}%
_{4,z}^{\frak{J}})$ vanishes,

(ii) $\widetilde{\mathbf{d}}_{5,z}^{\frak{J}}|\bigcirc^{5}\mathbf{K}%
_{4,z}^{\frak{J}}$ does not vanish.

\noindent We have ${\Sigma}E_{8}(N,P)\subset{\Sigma}^{\frak{J}_{4}}(N,P)$ by
(i), ($\frak{J}$E-2) and (E-ii). By (5) in Section2 for $z\in{\Sigma
}^{\frak{J}_{4}}(N,P)$, $\mathbf{K}_{4,z}^{\frak{J}}\cap T_{z}({\Sigma
}^{\frak{J}_{4}}(N,P))=\{\mathbf{0}\}$ if and only if $z\in{{\Sigma}%
}^{\frak{J}_{4},0}{{(N,P)}}$. Hence, we have ${{\Sigma}}E{{{{_{8}}%
}(N,P)\subset\Sigma}}^{\frak{J}_{4},0}{{(N,P)}}$ by (ii). By Lemma 3.1 (ii)
below, the section $\widetilde{\mathbf{d}}_{4}^{\frak{J}}|\bigcirc
^{3}\mathbf{K}_{4}^{\frak{J}}\bigcirc(\mathbf{K}_{2}^{\frak{J}}/\mathbf{K}%
_{4}^{\frak{J}})$ of Hom($\bigcirc^{3}\mathbf{K}_{4}^{\frak{J}}\bigcirc
(\mathbf{K}_{2}^{\frak{J}}/\mathbf{K}_{4}^{\frak{J}}),\mathbf{Q})$ over
${{\Sigma{{\overline{E}_{7}}}(N,P)}}$ is transverse to the zero-section, whose
inverse image of this section coincides with ${{\Sigma}}^{\frak{J}_{4}%
}{{(N,P)}}$. Hence, ${{\Sigma}}^{\frak{J}_{4}}{{(N,P)}}$ is a submanifold of
${{\Sigma{{\overline{E}_{7}}}(N,P)}}$ and the intrinsic derivative
$d(\widetilde{\mathbf{d}}_{4}^{\frak{J}}|\bigcirc^{3}\mathbf{K}_{4}^{\frak{J}%
}\bigcirc(\mathbf{K}_{2}^{\frak{J}}/\mathbf{K}_{4}^{\frak{J}}))$ over
${{\Sigma}}^{\frak{J}_{4}}{{(N,P)}}$\ is an epimorphism. Namely, we have the
exact sequence%
\[
\mathbf{0}\rightarrow T({{\Sigma}}^{\frak{J}_{4}}{{(N,P)}})\rightarrow
T({{\Sigma{{\overline{E}_{7}}}(N,P))|}}_{{{\Sigma}}^{\frak{J}_{4}}}%
\rightarrow\text{Hom}(\bigcirc^{3}\mathbf{K}_{4}^{\frak{J}}\bigcirc
(\mathbf{K}_{2}^{\frak{J}}/\mathbf{K}_{4}^{\frak{J}}),\mathbf{Q})|_{{{\Sigma}%
}^{\frak{J}_{4}}}\rightarrow\mathbf{0},
\]
which yields%
\begin{equation}
\mathbf{0}\rightarrow T({{\Sigma}}E{{{{_{8}}}(N,P)}})\rightarrow
T({{\Sigma{{\overline{E}_{7}}}(N,P))|}}_{{{\Sigma}}E{{{{_{8}}}}}}%
\rightarrow\text{Hom}(\bigcirc^{3}\mathbf{K}_{4}^{\frak{J}}\bigcirc
(\mathbf{K}_{2}^{\frak{J}}/\mathbf{K}_{4}^{\frak{J}}),\mathbf{Q})|_{{{\Sigma}%
}E{{{{_{8}}}}}}\rightarrow\mathbf{0}.
\end{equation}

(E-iii-a) By (E-i), we have that%
\[
T({{{{\Sigma}}^{\frak{J}_{2}}{{(N,P))|}}_{{{\Sigma{{E_{6}}}}}}\supset}%
}\mathbf{K}_{3}^{\frak{J}}{{{{|}}_{{{\Sigma{{E_{6}}}}}}\neq}}\{\mathbf{0}%
\}\text{ \ and \ }\mathbf{K}_{3}^{\frak{J}}{{{{|}}_{{{\Sigma{{E_{6}}}}}}\cap}%
}T({{\Sigma{{E_{6}}}(N,P))=}}\{\mathbf{0}\}.
\]

(E-iii-b) Since $\mathbf{K}_{3}^{\frak{J}}\cap T({{\Sigma}}^{\frak{J}_{3}%
,0}{{(N,P))=}}\{\mathbf{0}\}$, we have $\mathbf{K}_{3}^{\frak{J}}{{{{|}%
}_{{{\Sigma{{E_{7}}}}}}}}\cap T({{\Sigma}}^{\frak{J}_{3},0}{{(N,P)){{{|}%
}_{{{\Sigma{{E_{7}}}}}}}=}}\{\mathbf{0}\}$.

(E-iii-c) Let $r:\mathrm{Hom}(\bigcirc^{3}\mathbf{K}_{4}^{\frak{J}}%
\bigcirc\mathbf{K}_{2}^{\frak{J}},\mathbf{Q}){{{{|}}_{{{\Sigma{{E_{8}}}}}}}%
}\rightarrow\mathrm{Hom}(\bigcirc^{4}\mathbf{K}_{4}^{\frak{J}},\mathbf{Q}%
){{{{|}}_{{{\Sigma{{E_{8}}}}}}}}$ be the map induced from the restriction.
From this definition, it follows that%
\[
r\circ\mathbf{d}_{5}^{\frak{J}}{{|(\mathbf{K}_{4}^{\frak{J}}{{{{|}}%
_{{{\Sigma{{E_{8}}}}}})}}=}}d(\widetilde{\mathbf{d}}_{4}^{\frak{J}}%
|\bigcirc^{4}\mathbf{K}_{3}^{\frak{J}})|(\mathbf{K}_{4}^{\frak{J}}{{{{|}%
}_{{{\Sigma{{E_{8}}}}}})}}:\mathbf{K}_{4}^{\frak{J}}{{{{|}}_{{{\Sigma{{E_{8}}%
}}}}}}\rightarrow\text{Hom}(\bigcirc^{4}\mathbf{K}_{4}^{\frak{J}}%
,\mathbf{Q}){{{{|}}_{{{\Sigma{{E_{8}}}}}}}}%
\]
is an isomorphism. From (3.3), it follows that $\mathbf{K}_{4}^{\frak{J}%
}{{{{|}}_{{{\Sigma{{E_{8}}}}}}}}\cap T({{\Sigma{{\overline{E}_{7}(N,P))|}%
}_{{{\Sigma{{E_{8}}}}}}=}}\{\mathbf{0}\}{{,}}$ and the inclusion
$\mathbf{K}_{4}^{\frak{J}}{{{{|}}_{{{\Sigma{{E_{8}}}}}}\hookrightarrow}%
}T({{\Sigma{{\overline{E}_{6}(N,P)){{{|}}_{{{\Sigma{{E_{8}}}}}}}}}}}$ induces
the isomorphism $\mathbf{K}_{4}^{\frak{J}}|_{\Sigma{E_{8}}}\rightarrow
\nu(\overline{E}_{7}\subset\overline{E}_{6})|_{\Sigma{E_{8}}}$.

(E-iv) ${{\Sigma}}E{{{{_{k}}}(N,P)}}$ is of codimension $n-p+k$ in $J^{\infty
}(N,P)$.

(E-v) A smooth map germ $f:(N,x)\rightarrow(P,y)$ has $x$ as a singularity at
$x$ of type $E_{k}$ if and only if $j^{\infty}f$ is transverse to ${{\Sigma}%
}E{{{{_{k}}}(N,P)}}$ at $x$\ and $j_{x}^{\infty}f\in{{\Sigma}}E{{{{_{k}}%
}(N,P).}}$

When we deal with only $\Omega$-regular maps, it is convenient in the proof of
Lemma 3.1 and also in the calculation of Thom polynomials of ${{\Sigma}%
}E{{{{_{k}}}(N,P)}}$ in [An4] to set \ $\Sigma{{\overline{E}_{6}}%
}(N,P)={{\Sigma}}^{\frak{J}_{3}}{{(N,P)}}$ and ${{\Sigma{{\overline{E}_{8}}%
}(N,P)={{\Sigma}}^{\frak{J}_{4}}{{(N,P)}}}}$, since we can readily deal with
${{{{\Sigma}}^{\frak{J}_{3}}{{(N,P)}}}}$ and ${{{{\Sigma}}^{\frak{J}_{4}%
}{{(N,P)}}}}$. This does not matter, since the subset of all jets
$z\in{{\Sigma}}^{\frak{J}_{4},0}{{(N,P),}}$ such that $\widetilde{\mathbf{d}%
}_{5,z}^{\frak{J}}|\bigcirc^{5}\mathbf{K}_{4,z}^{\frak{J}}$ vanishes and
$\widetilde{\mathbf{d}}_{5,z}^{\frak{J}}|\bigcirc^{4}\mathbf{K}_{4,z}%
^{\frak{J}}\bigcirc(\mathbf{K}_{2,z}^{\frak{J}}/\mathbf{K}_{4,z}^{\frak{J}})$
does not vanish, has modality ([Ar, Lemma 6.1, 7$^{\circ}$]). Namely, such a
jet is not simple and does not appear in $\Omega(N,P)$.

\begin{lemma}
$\mathrm{(i)}$ The section $\widetilde{\mathbf{d}}_{4}^{\frak{J}}|\bigcirc
^{4}\mathbf{K}_{3}^{\frak{J}}$ of $\mathrm{Hom}(\bigcirc^{4}\mathbf{K}%
_{3}^{\frak{J}}\mathbf{,Q)}$ over ${{{{\Sigma}}^{\frak{J}_{3}}{{(N,P)}}}}$\ is
transverse to the zero-section on ${{\Sigma{{\overline{E}_{7}}}(N,P)}}$.

$\mathrm{(ii)}$ The section $\widetilde{\mathbf{d}}_{4}^{\frak{J}}%
|\bigcirc^{3}\mathbf{K}_{4}^{\frak{J}}\bigcirc(\mathbf{K}_{2}^{\frak{J}%
}/\mathbf{K}_{4}^{\frak{J}})$ of $\mathrm{Hom}(\bigcirc^{3}\mathbf{K}%
_{4}^{\frak{J}}\bigcirc(\mathbf{K}_{2}^{\frak{J}}/\mathbf{K}_{4}^{\frak{J}%
}),\mathbf{Q})$ over ${{\Sigma{{\overline{E}_{7}}}(N,P)}}$ is transverse to
the zero-section on ${{\Sigma{{E_{8}}}(N,P)}}$.
\end{lemma}

\begin{proof}
The proof is very like the proof of [An3, Lemmas 3.4 and Proposition 3.5]. We
only give the proof of (i), since the proof of (ii) is similar.

Take an element $z\in{{\Sigma}}\overline{E}_{7}{{(N,P)}}$ over $(x_{0}%
,y_{0}).$ We choose local coordinates $x=(x_{1},\cdots,x_{n-2},u,\ell)$ around
$x_{0}$ and $y=(y_{1},\cdots,y_{p})$ around $y_{0}$ such that $z$ is expressed
as $j_{0}^{\infty}h^{z}$ with
\begin{align}
y_{i}\circ h^{z}(x)  &  =x_{i}\text{ \ \ (}1\leq i\leq p-1\text{),}\nonumber\\
y_{p}\circ h^{z}(x)  &  =\pm x_{p}^{2}\pm\cdots\pm x_{n-2}^{2}+\overline
{h^{z}}(x_{1},\cdots,x_{p-1},u,\ell),
\end{align}
where $\overline{h^{z}}\in(x_{1},\cdots,x_{p-1},u,\ell)^{2}$ and
$(\partial^{4}\overline{h^{z}}/\partial\ell^{4})|_{x=0}=0.$\ Let
$\delta_{x_{1}},\cdots,\delta_{x_{n-2}},$ $\delta_{u}$ and $\delta_{\ell}$
express the vector fields of the total tangent bundle $\mathbf{D}$, which
correspond to the coordinates $x_{1},\cdots,x_{n-2},u,\ell$ respectively
defined in [B, Definition 1.6] (we avoid to use the notation used there, which
are used for the singularities $D_{k}$). Then $\mathbf{K}_{z}$ is spanned by
$\delta_{x_{p}},\cdots,\delta_{x_{n-2}},$ $\delta_{u}$ and $\delta_{\ell}$,
$\mathbf{K}_{2,z}^{\frak{J}}$ spanned by $\delta_{u}$ and $\delta_{\ell}$,
$\mathbf{K}_{3,z}^{\frak{J}}$ spanned by $\delta_{\ell}$, and $\mathbf{Q}_{z}$
spanned by the image of $\partial/\partial y_{p}$, say $\mathbf{e}_{p}$. We
consider the submanifold $S$ in ${{\Sigma}}^{\frak{J}_{3}}{{(N,P)\cap J}%
}_{x_{0},y_{0}}^{\infty}(N,P)$, which consists of all jets $j_{0}^{\infty}h$
such that $y_{i}\circ h(x)=x_{i}$\ ($1\leq i\leq p-1$) and that $y_{p}\circ
h(x)=y_{p}\circ h^{z}(x)+\overline{h}(x_{1},\cdots,x_{p-1},u,\ell)$, where
$\overline{h}$ is a polynomial of degree $5$. Then the vector fields
$\delta_{\ell}$ and $\mathbf{e}_{p}$ determine the trivializations
$\mathbf{K}_{3}|_{S}\cong S\times\mathbf{R}$ and $\mathbf{Q}|_{S}\cong
S\times\mathbf{R}$ respectively. Under these trivializations, we calculate
$\widetilde{\mathbf{d}}_{4}^{\frak{J}}|_{S}:\bigcirc^{4}\mathbf{K}%
_{3}^{\frak{J}}|_{S}\rightarrow\mathbf{Q}|_{S}$ as follows. Let $d^{j}%
=a^{j}\delta_{\ell}$ ($1\leq j\leq4$) be the vector fields of $\mathbf{K}%
_{3}^{\frak{J}}|_{S}$. Then we have
\begin{align*}
\widetilde{\mathbf{d}}_{4}^{\frak{J}}(d^{1}\bigcirc d^{2}\bigcirc
d^{3}\bigcirc d^{4})  &  =(a^{1}\partial/\partial\ell)(a^{2}\partial
/\partial\ell)(a^{3}\partial/\partial\ell)(a^{4}\partial/\partial\ell
)y_{p}\circ h\\
&  =a^{1}a^{2}a^{3}a^{4}(\partial^{4}y_{p}\circ h/\partial\ell^{4}).
\end{align*}
This implies that $\widetilde{\mathbf{d}}_{4}^{\frak{J}}|_{S}$ is transverse
to the zero-section of Hom$(\bigcirc^{4}\mathbf{K}_{3}^{\frak{J}},\mathbf{Q})$
at $z\in S\cap{{\Sigma{{\overline{E}_{7}}}(N,P)}}$. Hence, (i) is proved.
\end{proof}

\begin{remark}
We can prove by the above definition of $\Sigma X_{k}(N,P)$\ that the normal
form for $X_{k}$ in Introduction has the origin as the singularity of the type
$X_{k}$ respectively.
\end{remark}

\begin{remark}
$\mathrm{(i)}$ The submanifold $\Sigma D_{k}(N,P)$ is actually defined so that
it coincides with the inverse image of the submanifold in $J^{k-1}(N,P)$ by
$\pi_{k-1}^{\infty}$.

$\mathrm{(ii)}$ The submanifold $\Sigma E_{6}(N,P)$ $($resp. $\Sigma
E_{7}(N,P)$ and $\Sigma E_{8}(N,P))$ is actually defined so that it coincides
with the inverse image of the submanifold in $J^{4}(N,P)$ by $\pi_{4}^{\infty
}$ $($resp. in $J^{5}(N,P)$\ by $\pi_{5}^{\infty})$.
\end{remark}

\begin{remark}
We can entirely do the arguments in the definition of Boardman manifolds,
$\Sigma D_{i}(N,P)$ and $\Sigma E_{j}(N,P)$ on $J^{\ell}(N,P)$ for a large
$\ell$\ in Sections 2 and 3. We provide $N$ and $P$ with Riemannian metrics,
which enable us to consider the exponential maps $TN\rightarrow N$ and
$TP\rightarrow P$ by the Levi-Civita connections$.$ For points $x\in N$, $y\in
P$ and orthonormal basis of $T_{x}N$\ and $T_{y}P$, we take the associated
convex neighborhoods $U\subset N$\ around $x$\ and $V\subset P$\ around $y$
with the normal coordinate systems $(x_{1},...,x_{n})$ and $(y_{1},...,y_{p})$
so that $TU\approx U\times T_{x}N$\ and $TV\approx V\times T_{y}P$ by the
connections respectively (see [K-N])$.$ Let us define the canonical embedding
$\mu_{\infty}^{\ell}:J^{\ell}(TN,TP)\rightarrow J^{\infty}(TN,TP)$ by putting
the null homomorphisms of $\mathrm{Hom}(S^{i}(\pi_{N}^{\ast}(TN)),\pi
_{P}^{\ast}(TP))$ as the $i$-th component for $i>\ell$. We regard $\mu
_{\infty}^{\ell}$\ as the map to $J^{\infty}(N,P)$\ under the identification
(1.2). Any element $z\in\mu_{\infty}^{\ell}(J^{\ell}(TN,TP))$ is represented
by a $C^{\infty}$ polynomial map germ $f:(N,x)\rightarrow(P,y)$ of degree
$\ell$ under these coordinates. It is clear that $\pi_{\ell}^{\infty}\circ
\mu_{\infty}^{\ell}=id_{J^{\ell}(TN,TP)}$ and $\mu_{\infty}^{\ell}\circ
\pi_{\ell}^{\infty}|(\mu_{\infty}^{\ell}(J^{\ell}(TN,TP)))=id_{\mu_{\infty
}^{\ell}(J^{\ell}(TN,TP))}$.

We can prove that $\mathbf{D}|_{\mu_{\infty}^{\ell}(J^{\ell}(TN,TP))}$ is
tangent to $\mu_{\infty}^{\ell}(J^{\ell}(TN,TP))$. Indeed, let $x\in U\subset
N$ and $y\in V\subset P$ be as above. For any poins $u\in U$ and $v\in V$, we
take the normal coordinate systems $(u_{1},...,u_{n})$ around $u$ in $U$ and
$(v_{1},...,v_{p})$\ around $v$ in $V$ associated to the orthonormal basis of
$T_{u}U$\ and $T_{v}V$, which are induced from the above orthonormal basis of
$T_{x}N$\ and $T_{y}P$ by the parallel transformation along the geodesics of
$x$ to $u$ and $y$ to $v$ respectively. We note that $u_{1},...,u_{n}$ are
smoooth functions of $x_{1},...,x_{n}$. Let $\sigma=(\sigma_{1},...,\sigma
_{n})$ with non-negative integers $\sigma_{i}$. We define the coordinate
system $X_{i}$, $Y_{j}$ and $W_{j,\sigma}$ of $J^{\infty}(TU,TV)$ by%
\[
X_{i}=x_{i}\circ\pi_{U}^{\infty}\text{, }Y_{j}=y_{j}\circ\pi_{V}^{\infty
}\text{, }W_{j,\sigma}(j_{u}^{\infty}f)=\frac{\partial^{|\sigma|}(v_{j}\circ
f)}{\partial u^{\sigma}}(u),
\]
where $|\sigma|\geq0$\ and $Y_{j}=W_{j,\mathbf{0}}$ for $|\sigma|=0$. In the definition,
we should note that the normal coordinate systems $(u_{1},...,u_{n})$ and
$(v_{1},...,v_{p})$\ vary depending on points $u$\ and $v$.

A smooth function $\Phi$ defined on an open subset of $\mu_{\infty}^{\ell
}(J^{\ell}(TN,TP))$ is written as $\Phi\circ\mu_{\infty}^{\ell}\circ\pi_{\ell
}^{\infty}$ on the same open subset. Hence, $\Phi$ is extended to $(\Phi
\circ\mu_{\infty}^{\ell})\circ\pi_{\ell}^{\infty}$ defined on an open subset
of $J^{\infty}(TN,TP)$, which is a smooth function in the sence of [B,
Definition (1.4)]. Let us consider the total tangent bundle $\mathbf{D}$
defined on $J^{\infty}(TN,TP)=J^{\infty}(N,P)$\ and its vector field $D_{i}$
corresponding to $\partial/\partial x_{i}$. By using [B, (1.8)], we have, for
$z=j_{u}^{\infty}f\in\mu_{\infty}^{\ell}(J^{\ell}(TN,TP))$,%
\begin{align*}
D_{i}(\Phi)(z)  &  =\frac{\partial(\Phi\circ j^{\infty}f)}{\partial x_{i}%
}(u)\\
&  =\frac{\partial\Phi}{\partial X_{i}}(z)\frac{\partial(X_{i}\circ j^{\infty
}f)}{\partial x_{i}}(u)+\underset{j,\sigma}{\sum}\frac{\partial\Phi}{\partial
W_{j,\sigma}}(z)\frac{\partial(W_{j,\sigma}\circ j^{\infty}f)}{\partial x_{i}%
}(u)\\
&  =\frac{\partial\Phi}{\partial X_{i}}(z)+\underset{j,\sigma}{\sum}%
\frac{\partial\Phi}{\partial W_{j,\sigma}}(z)\left(  \underset{}{\sum_{h}%
\frac{\partial(W_{j,\sigma}\circ j^{\infty}f)}{\partial u_{h}}(u)}%
\frac{\partial u_{h}}{\partial x_{i}}(u)\right)  ,
\end{align*}
where

(i) if we define $g_{i}^{h}(u)=\frac{\partial u_{h}}{\partial x_{i}}(u)$, then
they are smooth functions of $x_{1},...,x_{n}$,

(ii) $\frac{\partial W_{j,\sigma}}{\partial u_{h}}=W_{j,\sigma^{\prime}}$ with
$\sigma^{\prime}=(\sigma_{1},...,\sigma_{h-1},\sigma_{h}+1,\sigma
_{h+1},...,\sigma_{n}),$

(iii) $|\sigma|\leq\ell-1$, since $W_{j,\sigma^{\prime}}$ vanishes for
$|\sigma|=\ell.$

Since $\Phi$ is a function of variables $X_{i}$, $Y_{j}$ and $W_{j,\sigma}$
with $|\sigma|\leq\ell$, so is $D_{i}(\Phi)$. Therefore, $D_{i}$ is tangent to
$\mu_{\infty}^{\ell}(J^{\ell}(TN,TP))$ at $z$ for each $i$. Since $\mathbf{D}$
is generated by $D_{i}$, $\mathbf{D}_{z}$ is tangent to $\mu_{\infty}^{\ell
}(J^{\ell}(TN,TP))$ at $z$.
\end{remark}

For a symbol $I=(i_{1},i_{2},\cdots,i_{r})$ we define $\Omega^{I}(n,p)$ to be
a subset of $J^{\infty}(n,p)$ which consists of all Boardman manifolds
$\Sigma^{J}(n,p)$ with symbols $J$ of length $r$ satisfying $J\leq I$ in the
lexicographic order. The following lemma is well known.

\begin{lemma}
$\Omega^{I}(n,p)$ is an open subset of $J^{\infty}(n,p)$.
\end{lemma}

We next describe the adjacency relations between $\Sigma A_{i}(n,p)$, $\Sigma
D_{j}(n,p)$ and $\Sigma E_{k}(n,p)$ in $J^{\infty}(n,p)$. This adjacency
relations are quite parallel to the result given in [Ar, Corollary 8.7].
However, we need to explain them, since $\Sigma X_{j}(n,p)$ is not necessarily
an orbit through the jet of a germ with singularity $X_{j}$ given in
Introduction by $L^{\infty}(p)\times L^{\infty}(n)$.

\begin{lemma}
Let $n\geq p\geq2$. A subset $\Omega(n,p)$ consisting of all regular jets and
a number of $\Sigma A_{i}(n,p)$, $\Sigma D_{j}(n,p)$ and $\Sigma E_{k}(n,p)$
is an open subset of $J^{\infty}(n,p)$ if and only if the following three
conditions are satisfied.

$(\mathrm{i})$ If $\Sigma A_{i}(n,p)\subset\Omega(n,p)$, then $\Sigma A_{\ell
}(n,p)\subset\Omega(n,p)$ for all $\ell$ with $1\leq\ell<i$.

$(\mathrm{ii})$ If $\Sigma D_{i}(n,p)\subset\Omega(n,p)$, then $\Sigma
A_{\ell}(n,p)$ $(1\leq\ell<i)$ and $\Sigma D_{\ell}(n,p)$ $(4\leq\ell<i)$\ are
all contained in $\Omega(n,p)$.

$(\mathrm{iii})$ If $\Sigma E_{i}(n,p)\subset\Omega(n,p)$, then $\Sigma
A_{\ell}(n,p)$ $(1\leq\ell<i)$, $\Sigma D_{\ell}(n,p)$ $(4\leq\ell<i)$ and
$\Sigma E_{\ell}(n,p)$ $(6\leq\ell<i)$ are all contained in $\Omega(n,p)$.
\end{lemma}

\begin{proof}
Suppose that $\Omega(n,p)$ is open. Let $z=j_{0}^{\infty}f$ be a jet of
$\Sigma X_{i}(n,p)$. Even if $j^{\infty}f:(\mathbf{R}^{n},0)\rightarrow
(J^{\infty}(n,p),z)$ is not transverse to $\Sigma X_{i}(n,p)$, there exists a
sequence $\{z_{k}\}$ with $z_{k}=j_{0}^{\infty}f_{k}$ such that $\lim
_{k\rightarrow\infty}z_{k}=z$, $z_{k}\in\Sigma X_{i}(n,p)$ and $j^{\infty
}f_{k}$ is transverse to $\Sigma X_{i}(n,p)$ at $0$. Therefore, (i), (ii) and
(iii) follow from the adjacency relation in [Ar, Corollary 8.7].

We next assume (i), (ii) and (iii). Suppose to the contrary that $\Omega(n,p)$
is not open. Then there exist a jet $z\in$ $\Omega(n,p)$ and a sequence
$\{z_{k}\}$ with $z_{k}=j_{0}^{\infty}f_{k}$ such that $\lim_{k\rightarrow
\infty}z_{k}=z$ and $z_{k}\notin\Omega(n,p)$ for all $k$. If there exists
infinite $k$'s with $z_{k}\notin\Omega^{(\frak{I}_{4},0)}(n,p)$, then we have
$z\notin\Omega^{(\frak{I}_{4},0)}(n,p)$\ by Lemma 3.5. This is a
contradiction. Note that $\Sigma^{J_{2}}(n,p)\setminus(\cup_{i=1}^{\infty
}\Sigma A_{i}(n,p))$, $\Sigma^{(\frak{I}_{2},0)}(n,p)\setminus(\cup
_{i=4}^{\infty}\Sigma D_{i}(n,p))$ and $\Sigma^{\frak{I}_{3}}(n,p)\setminus
(\cup_{i=6}^{8}\Sigma E_{i}(n,p))$ are closed subsets of $\Sigma^{J_{2}}%
(n,p)$, $\Sigma^{(\frak{I}_{2},0)}(n,p)$ and $\Sigma^{\frak{I}_{3}}(n,p)$
respectively. Hence, we may assume that there exists a number $q$ and infinite
$k$'s with $z_{k}\in\Sigma\overline{X}_{q}(n,p)$. Let $q_{0}$ be the smallest
number among such number $q$'s with $z_{k}\in\Sigma\overline{X}_{q}(n,p)$.
From the adjacency relation above there exists a number $q\geq q_{0}$ such
that $z$ lies in $\Sigma A_{q}(n,p)$, $\Sigma D_{q}(n,p)$ or $\Sigma
E_{q}(n,p)$. Consequently, it follows that

(A) if $z\in\Sigma A_{q}(n,p)$, then we have that $X_{q_{0}}=A_{q_{0}}$ and
$z_{k}\in\Sigma A_{q_{0}}(n,p)\subset\Omega(n,p)$ by (i),

(D) if $z\in\Sigma D_{q}(n,p)$, then we have that $X_{q_{0}}=A_{q_{0}}$ or
$D_{q_{0}}$ and $z_{k}\in\Sigma X_{q_{0}}(n,p)\subset\Omega(n,p)$ by (ii),

(E)\ if $z\in\Sigma E_{q}(n,p)$, then we have that $X_{q_{0}}=A_{q_{0}}$,
$D_{q_{0}}$\ or $E_{q_{0}}$ and $z_{k}\in\Sigma X_{q_{0}}(n,p)\subset
\Omega(n,p)$ by (iii).

This is a contradiction.
\end{proof}

\section{Primary obstruction}

In what follows we set $\Sigma D_{j}(N,P)=\emptyset$ for $1\leq j\leq3$ and
$\Sigma E_{j}(N,P)=\emptyset$ unless $6\leq j\leq8$. For a section $s\in
\Gamma_{\Omega}^{tr}(N,P)$, we set $S^{I_{k}}(s)=s^{-1}(\Sigma^{I_{k}}(N,P))$,
$S^{X_{k}}(s)=s^{-1}(\Sigma X_{k}(N,P))$, $S^{\overline{X}_{k}}(s)=s^{-1}%
(\Sigma\overline{X}_{k}(N,P))$, where $X_{k}$ refers to $A_{k}$, $D_{k}$ or
$E_{k}$. We have $S^{\overline{X}_{k}}(s)=\cup_{j=k}^{\infty}S^{X_{j}}(s)$. We
often write $S^{I_{k}}$, $S^{X_{k}}$ and $S^{\overline{X}_{k}}$ by omitting
$(s)$ if there is no confusion. Furthermore, we set $(s|S^{X_{k}})^{\ast
}(\mathbf{K})=K(S^{X_{k}}(s))$, $(s|S^{X_{k}})^{\ast}\mathbf{Q}=Q(S^{X_{k}%
}(s))$ and $(s|S^{X_{k}})^{\ast}(\mathbf{K}_{j}^{I})=K_{j}^{I}(S^{X_{k}}(s))$.

We may assume that the section $s$ given in Theorem 0.1 lies in $\Gamma
_{\Omega}^{tr}(N,P)$. Let $C_{k+1}$ refer to the union $C\cup S^{\overline
{A}_{k+1}}(s)\cup S^{\overline{D}_{k+1}}(s)\cup S^{\overline{E}_{k+1}}(s)$
($k\geq0$). We assume that there exists an $\Omega$-regular map $g_{k+1}$
defined on a neighborhood of $C_{k+1}$, where $j^{\mathbf{\infty}}g_{k+1}=s$ holds.

\begin{theorem}
Let $n\geq p\geq2$ and $k\geq1$. Let $N$, $P$, $\Omega(N,P)$ and $s\in
\Gamma_{\Omega}^{tr}(N,P)$ be given as in Theorem 0.1. Let $C_{k+1}$ and
$g_{k+1}$ be given as above. Then there exists a homotopy $s_{\lambda}%
\in\Gamma_{\Omega}(N,P)$ relative to a neighborhood of $C_{k+1}$ with the
following properties.

$(4.1.1)$ $s_{0}=s$ and $s_{1}\in\Gamma_{\Omega}^{tr}(N,P)$.

$(4.1.2)$ There exists an $\Omega$-regular map $g_{k}$ which is defined on a
neighborhood of $C_{k}$ for $k>1$ and on $N$ for $k=1$, where
$j^{\mathbf{\infty}}g_{k}=s_{1}$ holds.

$(4.1.3)$ In the situation that $\Omega(N,P)$ does not contain $\Sigma
D_{j}(N,P)$ or $\Sigma E_{k}(N,P)$, we have that $s_{\lambda}^{-1}(\Sigma
A_{i}(N,P))=s^{-1}(\Sigma A_{i}(N,P))$ for $k\geq i>1$.
\end{theorem}

The case $k=1$ of Theorem 4.1 follows from [An8, Theorem 0.5]. In fact, let
$N_{0}=N\setminus(S^{\overline{A}_{2}}(s)\cup S^{\overline{D}_{4}}(s)\cup
S^{\overline{E}_{6}}(s))$, $\Omega^{n-p+1,0}(N_{0},P)=\Sigma^{n-p}%
(N_{0},P)\cup\Sigma^{n-p+1,0}(N_{0},P)$ and $U$, $U^{\prime}$ be closed
neighborhoods of $C_{2}$ with $U\subset$Int$U^{\prime}$, where $g_{2}$ is
defined. Since $s\in\Gamma_{\Omega}^{tr}(N_{{}},P)$, $s|N_{0}$ is a section of
$\Gamma_{\Omega^{n-p+1,0}}^{tr}(N_{0},P)$ and $g_{2}|U^{\prime}\cap N_{0}$ is
a fold-map. Hence, we obtain a homotopy $u_{\lambda}\in\Gamma_{\Omega
^{n-p+1,0}}(N_{0},P)$ relative to a neighborhood of $U\cap N_{0}$ and a fold
map $f_{0}:N_{0}\rightarrow P$ such that $s_{0}|N_{0}=u_{0}$ and
$u_{1}=j^{\mathbf{\infty}}f_{0}$. Then we obtain a homotopy $s_{\lambda}$
required in the case $k=1$ of Theorem 4.1 by defining $s_{\lambda}%
|N_{0}=u_{\lambda}$ and $s_{\lambda}|U=j^{\mathbf{\infty}}g_{2}$.

We will prove Theorem 4.1 for $k\geq2$ in Section 6.

Here we give a proof of Theorem 0.1 by using Theorem 4.1.

\begin{proof}
[Proof of Theorem 0.1]We may assume without loss of generality that
$s\in\Gamma_{\Omega}^{tr}(N,P)$. Then we may set $C_{n+1}=C$. By the
assumption there exists the $\Omega$-regular map $g_{n+1}=g$ defined on a
neighborhood of $C_{n+1}$. Then we can prove Theorem 0.1 by the downward
induction on $k$, and by using the above argument for the final step $k=1$.
\end{proof}

\begin{remark}
In Theorem 0.1 we assume, in addition, that $\Omega(N,P)$ does not contain any
$\Sigma D_{j}(N,P)$ or $\Sigma E_{k}(N,P)$ and that $s\in\Gamma_{\Omega}%
^{tr}(N,P)$. Then we have that $s_{\lambda}^{-1}(\Sigma A_{i}(N,P))=s^{-1}%
(\Sigma A_{i}(N,P))=(j^{\mathbf{\infty}}f)^{-1}(\Sigma A_{i}(N,P))$ for each
$i$ with $i>1$.
\end{remark}

In what follows we take a number $k$ with $k>1$. We begin by preparing several
notions and results, which are necessary for the proof of Theorem 4.1. For the
map $g_{k+1}$ and the closed subset $C_{k+1}$, we take an open neighborhood
$U(C_{k+1})^{\prime}$ of $C_{k+1}$, where $j^{\mathbf{\infty}}g_{k+1}=s$.
Without loss of generality we may assume that $N\setminus U(C_{k+1})^{\prime}$
is nonempty. Take a smooth function $h_{C_{k+1}}:N\rightarrow\lbrack0,1]$ such
that
\begin{equation}
\left\{
\begin{array}
[c]{ll}%
h_{C_{k+1}}(x)=1 & \text{for }x\in C_{k+1},\\
h_{C_{k+1}}(x)=0 & \text{for }x\in N\setminus U(C_{k+1})^{\prime},\\
0<h_{C_{k+1}}(x)<1 & \text{for }x\in U(C_{k+1})^{\prime}\setminus C_{k+1}.
\end{array}
\right.
\end{equation}
By the Sard Theorem ([H2]) there is a regular value $r$ of $h_{C_{k+1}}$ with
$0<r<1$. Then $h_{C_{k+1}}^{-1}(r)$ is a submanifold and we set $U(C_{k+1}%
)=h_{C_{k+1}}^{-1}([r,1])$. We decompose $N\setminus\mathrm{Int}U(C_{k+1})$
into the connected components, say $L_{1},\ldots,L_{j},\ldots$. It suffices to
prove Theorem 4.1 for each $L_{j}\cup$Int$U(C_{k+1})$. Since $\partial
N=\emptyset$, we have that $N\setminus U(C_{k+1})$ has empty boundary. If
$L_{j}$ is not compact, then Theorem 4.1 holds for $L_{j}\cup$Int$U(C_{k+1})$
by Gromov's theorem ([G1, Theorem 4.1.1]). Therefore, it suffices to consider
the special case where

(C1) $N\setminus\mathrm{Int}U(C_{k+1})$ is compact, connected and nonempty,

(C2) $\partial U(C_{k+1})$ is a submanifold of dimension $n-1$,

(C3) for the smooth function $h_{C_{k+1}}:N\rightarrow\lbrack0,1]$ satisfying
(4.1) there is a sufficiently small positive real number $\varepsilon$ with
$r-2\varepsilon>0$ such that $r-t\varepsilon$ ($0\leq t\leq2$) are all regular
values of $h_{C_{k+1}}$. We have that $h_{C_{k+1}}^{-1}([r-2\varepsilon,1])$
is contained in $U(C_{k+1})^{\prime}$.

We set $U(C_{k+1})_{t}=h_{C_{k+1}}^{-1}([r-(2-t)\varepsilon,1])$. In
particular, we have $U(C_{k+1})_{2}=U(C_{k+1})$. Furthermore, we may assume that

(C4) $s\in\Gamma_{\Omega}^{tr}(N,P)$ and $S^{X_{j}}(s)$ ($j\leq k$) are
transverse to $\partial U(C_{k+1})_{0}$ and $\partial U(C_{k+1})_{2}$ on a
neighborhood of $S^{X_{k}}(s)$.

In what follows we choose and fix a Riemannian metric of $N$, which satisfies

\textbf{Orthogonality Condition}: Let\textit{ }$I=J$\textit{ or }$\frak{J}%
$\textit{. If }$K_{j}^{I}(S^{X_{k}})/K_{j+1}^{I}(S^{X_{k}})$\textit{ is of
positive dimension, then }$K_{j}^{I}(S^{X_{k}})/K_{j+1}^{I}(S^{X_{k}}%
)$\textit{ is orthogonal to }$S^{I_{j}}(s)$\textit{ in }$S^{I_{j-1}}%
(s)$\textit{ over }$S^{X_{k}}(s)$\textit{\ for }$k\geq j\geq1$\textit{
}$(S^{I_{0}}(s)=N)$\textit{.}

Let $\nu(\overline{X}_{k})$ be the normal bundle $(T(J^{\mathbf{\infty}%
}(N,P))|_{\Sigma{\overline{X}_{k}}})/T(\Sigma{\overline{X}_{k}}(N,P))$ and
$\nu(X_{k})=\nu(\overline{X}_{k})|_{\Sigma X_{k}}$. Let $\mathbf{j}%
_{\mathbf{K}}:\mathbf{K}\rightarrow\nu({{{{\overline{X}_{k})}}}}$ over $\Sigma
X_{k}(N,P)$ be the composition of the inclusion $\mathbf{K}\rightarrow
{{{{T(J}}}}^{\mathbf{\infty}}{{(N,P))}}$ and the projection ${{{{T(J}}}%
}^{\mathbf{\infty}}{{(N,P))\rightarrow}}\nu({{{{\overline{X}_{k})}}}}$. We
have the monomorphism%
\[
\mathbf{j}_{\mathbf{K}}\circ(s|S^{X_{k}})^{\mathbf{K}}:K(S^{X_{k}}%
)\rightarrow\mathbf{K}\rightarrow\nu(X_{k}).
\]

For $s\in\Gamma_{\Omega}^{tr}(N,P)$, let $\frak{n}(s,X_{k})$ or simply
$\frak{n}(X_{k})$ be the orthogonal normal bundle of $S^{X_{k}}(s)$ in $N$.
Furthermore, we have the bundle map
\[
ds|\frak{n}(s,X_{k}):\frak{n}(s,X_{k})\rightarrow\nu(X_{k})
\]
covering $s|S^{X_{k}}(s):S^{X_{k}}(s)\rightarrow\Sigma X_{k}(N,P)$. Let
$\mathbf{i}_{\frak{n}(s,X_{k})}:\frak{n}(s,X_{k})\subset TN|_{S^{X_{k}}}$
denote the inclusion. We define $\Psi(s,X_{k}):K(S^{X_{k}})|_{S^{X_{k}}%
}\rightarrow\frak{n}(s,X_{k})\subset TN|_{S^{X_{k}}}$ to be the composition%
\begin{align*}
&  \mathbf{i}_{\frak{n}(s,X_{k})}\circ(ds|\frak{n}(s,X_{k}))^{-1}%
\circ\mathbf{j}_{\mathbf{K}}\circ(s|S^{X_{k}})^{\mathbf{K}}|(K(S^{X_{k}%
})|_{S^{X_{k}}})\\
\text{ \ \ \ \ }  &  :K(S^{X_{k}})|_{S^{X_{k}}}\rightarrow\nu(X_{k}%
)\rightarrow\frak{n}(s,X_{k})\hookrightarrow TN|_{S^{X_{k}}}.
\end{align*}
Let $i_{K(S^{X_{k}})}:K(S^{X_{k}})\rightarrow TN$ be the inclusion.

\begin{remark}
If $f$ is an $\Omega$-regular map, then it follows from the definition of the
total tangent bundle $\mathbf{D}$ that $i_{K(S^{X_{k}}(j^{\mathbf{\infty}}%
f))}=\Psi(j^{\mathbf{\infty}}f,X_{k})$.
\end{remark}

In what follows let $M=S^{X_{k}}(s)\setminus$Int$(U(C_{k+1}))$. Let
Mono$(K(S^{X_{k}})|_{M},TN|_{M})$ denote the subset of Hom$(K(S^{X_{k}}%
)|_{M},TN|_{M})$, which consists of all monomorphisms $K(S^{X_{k}}%
)_{c}\rightarrow T_{c}N$, $c\in M$. We denote the bundle of the local
coefficients $\mathcal{B}(\pi_{j}(\mathrm{Mono}(K(S^{X_{k}})_{c},T_{c}N))),$
$c\in M,$ by $\mathcal{B}(\pi_{j})$, which is a covering space over $M$ with
fiber $\pi_{j}(\mathrm{Mono}(K(S^{X_{k}})_{c},T_{c}N))$ defined in [Ste,
30.1]. From the obstruction theory due to [Ste, 36.3], it follows that the
obstructions for $i_{K(S^{X_{k}})}|_{M}$ and $\Psi(s,X_{k})|_{M}$ to be
homotopic are the primary differences $d(i_{K(S^{X_{k}})}|_{M},\Psi
(s,X_{k})|_{M})$, which are defined in the cohomology groups with local
coefficients $H^{j}(M,\partial M;\mathcal{B}(\pi_{j}))$. We show that all of
them vanish. In fact, we have that $\dim M=n-(n-p+k)=p-k<p-1$. Since
\textrm{Mono}$(\mathbf{R}^{n-p+1},\mathbf{R}^{n})$ is identified with
$GL(n)/GL(p-1)$, it follows from [Ste, 25.6] that $\pi_{j}($\textrm{Mono}%
$(\mathbf{R}^{n-p+1},\mathbf{R}^{n}))\cong\{\mathbf{0}\}$ for $j<p-1$. Hence,
there exists a homotopy $\psi^{M}(s,X_{k})_{\lambda}:K(S^{X_{k}}%
)|_{M}\rightarrow TN|_{M}$ relative to $M\cap U(C_{k+1})_{1}$ in
\textrm{Mono}$(K(S^{X_{k}})|_{M},TN|_{M})$ such that $\psi^{M}(s,X_{k}%
)_{0}=i_{K(S^{X_{k}})}|_{M}$ and $\psi^{M}(s,X_{k})_{1}$$=\Psi(s,X_{k})|_{M}$
by the definition of the primary difference. Let $\mathrm{Iso}(TN|_{M}%
,TN|_{M})$\ denote the subspace of $\mathrm{Hom}(TN|_{M},TN|_{M})$, which
consists of all isomorphisms of $T_{c}N$, $c\in M$.\ The restriction map
\[
r_{M}:\mathrm{Iso}(TN|_{M},TN|_{M})\rightarrow\mathrm{Mono}(K(S^{X_{k}}%
)|_{M},TN|_{M})
\]
defined by $r_{M}(h)=h|(K(S^{X_{k}})_{c})$, $h\in\mathrm{Iso}(T_{c}N,T_{c}N)$,
induces a structure of a fiber bundle with fiber Iso$(\mathbf{R}%
^{p-1},\mathbf{R}^{p-1})\times{\mathrm{Hom}}(\mathbf{R}^{p-1},\mathbf{R}%
^{n-p+1})$. By applying the covering homotopy property of the fiber bundle
$r_{M}$ to the sections $id_{TN|_{M}}$ and the homotopy $\psi^{M}%
(s,X_{k})_{\lambda},$ we obtain a homotopy $\Psi^{M}(s,X_{k})_{\lambda}:$
$TN|_{M}\rightarrow TN|_{M}$ such that $\Psi^{M}(s,X_{k})_{0}=id_{TN|_{M}}$,
$\Psi^{M}(s,X_{k})_{\lambda}|_{c}=id_{T_{c}N}$ for all $c\in M\cap
U(C_{k+1})_{1}$ and $r_{M}\circ\Psi^{M}(s,X_{k})_{\lambda}=\psi^{M}%
(s,X_{k})_{\lambda}$. We define $\Phi(s,X_{k})_{\lambda}:$ $TN|_{M}\rightarrow
TN|_{M}$ by $\Phi(s,X_{k})_{\lambda}=(\Psi^{M}(s,X_{k})_{\lambda})^{-1}$.

Let us fix a direct sum decomposition%
\begin{equation}%
\begin{array}
[c]{l}%
\nu({{{{A_{k})=\oplus}}}}_{j=1}^{k}\nu(J_{j}\subset J_{j-1}),\\
\nu({{{{D_{k})=(\oplus}}}}_{j=1}^{2}\nu(\frak{J}_{j}\subset\frak{J}%
_{j-1}))\oplus(\oplus_{j=5}^{k}\nu(\overline{D}_{j}\subset\overline{D}%
_{j-1})),\\
\nu({{{{E_{k})=(\oplus}}}}_{j=1}^{3}\nu(\frak{J}_{j}\subset\frak{J}%
_{j-1}))\oplus(\oplus_{j=7}^{k}\nu(\overline{E}_{j}\subset\overline{E}%
_{j-1})),
\end{array}
\end{equation}
over $\Sigma{X_{k}}(N,P)$, which induces the direct sum decomposition
$\mathbf{K=K/K}_{2}^{I}\oplus\mathbf{K}_{2}^{I}\mathbf{/K}_{3}^{I}%
\oplus\mathbf{K}_{3}^{I}$.

Let $\frak{n}(s,\overline{X}_{j}\subset\overline{X}_{j-1})$ be the orthogonal
normal bundle of $S^{\overline{X}_{j}}(s)$ in $S^{\overline{X}_{j-1}}(s)$ and
set $\frak{n}(s,X_{j}\subset\overline{X}_{j-1})=\frak{n}(s,\overline{X}%
_{j}\subset\overline{X}_{j-1})|_{S^{X_{j}}}$ where $\overline{X}_{j}=I_{j}$,
or $X_{j}=A_{j}$, $D_{j}$, $E_{j}$ . Then we have the canonical direct sum
decomposition such as%
\begin{equation}%
\begin{array}
[c]{l}%
\frak{n}(s,{{{{A_{k})=\oplus}}}}_{j=1}^{k}\frak{n}(s,J_{j}\subset J_{j-1}),\\
\frak{n}(s,{{{{D_{k})=(\oplus}}}}_{j=1}^{2}\frak{n}(s,\frak{J}_{j}%
\subset\frak{J}_{j-1}))\oplus(\oplus_{j=5}^{k}\frak{n}(s,\overline{D}%
_{j}\subset\overline{D}_{j-1})),\\
\frak{n}(s,{{{{E_{k})=(\oplus}}}}_{j=1}^{3}\frak{n}(s,\frak{J}_{j}%
\subset\frak{J}_{j-1}))\oplus(\oplus_{j=7}^{k}\frak{n}(s,\overline{E}%
_{j}\subset\overline{E}_{j-1})).
\end{array}
\end{equation}
over $S^{X_{k}}(s)$. We can take the direct sum decompositions in (4.2) to be
compatible with those in (4.3).

\section{Lemmas}

In the proof of the following lemma, $\Phi(s,X_{k})_{\lambda}|_{c}$ ($c\in M$)
is regarded as a linear isomorphism of $T_{c}N$. Let $r_{0}$ be a small
positive real number with $r_{0}<1/10$.

\begin{lemma}
Let $k>1$. Let $X_{k}$ be any of $A_{k}$, $D_{k}$ and $E_{k}$. Let $s\in
\Gamma_{\Omega}^{tr}(N,P)$ be a section satisfying the hypotheses of Theorem
$4.1$. Then there exists a homotopy $s_{\lambda}$ relative to $U(C_{k+1}%
)_{2-3r_{0}}$ in $\Gamma_{\Omega}^{tr}(N,P)$ with $s_{0}=s$ satisfying

$(5.1.1)$ for any $\lambda$, $S^{X_{k}}(s_{\lambda})=S^{X_{k}}(s)$ and
$\pi_{P}^{\mathbf{\infty}}\circ s_{\lambda}|S^{X_{k}}(s_{\lambda})=\pi
_{P}^{\mathbf{\infty}}\circ s|S^{X_{k}}(s),$

$(5.1.2)$ we have that $i_{K(S^{X_{k}}(s_{1}))}|_{S^{X_{k}}}=\Psi(s_{1}%
,X_{k})$ and $K(S^{X_{k}}(s_{1}))_{c}\subset\frak{n}(s,X_{k})_{c}$ for any
point $c\in S^{X_{k}}(s_{1})$.
\end{lemma}

\begin{proof}
Recall the exponential map $\exp_{N,x}:T_{x}N\rightarrow N$ defined near
$\mathbf{0}\in T_{x}N$. We write an element of $\frak{n}(X_{k})_{c}$ as
$\mathbf{v}_{c}$. There exists a small positive number $\delta$ such that the
map
\[
e:D_{\delta}(\frak{n}(X_{k})_{c})|_{M}\rightarrow N
\]
defined by $e(\mathbf{v}_{c})=\exp_{N,c}(\mathbf{v}_{c})$ is an embedding and
that the images of $e$ for $X_{k}=A_{k}$, $D_{k}$ and $E_{k}$ do not mutually
intersect, where $c\in M$ and $\mathbf{v}_{c}\in D_{\delta}(\frak{n}%
(X_{k})_{c})$ (note that $e|M$ is the inclusion). Let $\rho:[0,\infty
)\rightarrow\mathbf{R}$ be a decreasing smooth function such that $0\leq
\rho(t)\leq1$, $\rho(t)=1$ if $t\leq\delta/10$ and $\rho(t)=0$ if $t\geq
\delta$.

If we represent $s(x)\in\Omega(N,P)$ by a jet $j_{x}^{\mathbf{\infty}}%
\sigma_{x}$ for a germ $\sigma_{x}:(N,x)\rightarrow(P,\pi_{P}^{\mathbf{\infty
}}\circ s(x))$, then we define the homotopy $s_{\lambda}$ of $\Gamma_{\Omega
}^{tr}(N,P)$ using $\Phi(s,X_{k})_{\lambda}$ by%
\begin{equation}
\left\{
\begin{array}
[c]{ll}%
\begin{array}
[c]{l}%
s_{\lambda}(e(\mathbf{v}_{c}))\\
=j_{e(\mathbf{v}_{c})}^{\mathbf{\infty}}(\sigma_{e(\mathbf{v}_{c})}\circ
\exp_{N,c}\circ\Phi(s,X_{k})_{\rho(\Vert\mathbf{v}_{c}\Vert)\lambda}|_{c}%
\circ\exp_{N,c}^{-1})
\end{array}
& \text{\textrm{if} $c\in M$ \textrm{and} }\Vert\text{$\mathbf{v}_{c}\Vert
\leq\delta,$}\\%
\begin{array}
[c]{l}%
s_{\lambda}(x)=s(x)
\end{array}
& \text{\textrm{if} $x\notin\mathrm{Im}(e).$}%
\end{array}
\right.
\end{equation}
Here, $\Phi(s,X_{k})_{\rho(\Vert\mathbf{v}_{c}\Vert)\lambda}|_{c}$\ refers to
$\ell(\mathbf{v}_{c})\circ(\Phi(s,X_{k})_{\rho(\Vert\mathbf{v}_{c}%
\Vert)\lambda}|_{c})\circ\ell(-\mathbf{v}_{c})$, where $\ell(\mathbf{v})$ is
the parallel translation defined by $\ell(\mathbf{v})(\mathbf{a}%
)=\mathbf{a}+\mathbf{v}$. If $\Vert\mathbf{v}_{c}\Vert\geq\delta$, then
$\Phi(s,X_{k})_{\rho(\Vert\mathbf{v}_{c}\Vert)\lambda}|_{c}=\Phi(s,X_{k}%
)_{0}|_{c}$, and if $c\in M\cap U(C_{k+1})_{2-3r_{0}}$, then $\Phi
(s,X_{k})_{\lambda}|_{c}=\Phi(s,X_{k})_{0}|_{c}$. Hence, $s_{\lambda}$ is well
defined. Furthermore, we have that

(1) $\pi_{P}^{\mathbf{\infty}}\circ s_{\lambda}(x)=\pi_{P}^{\mathbf{\infty}%
}\circ s(x)$,

(2) $S^{X_{k}}(s_{\lambda})=S^{X_{k}}(s)$,

(3) if $c\in S^{X_{k}}(s)$, then we have that $\frak{n}(s_{1},X_{k}%
)_{c}\supset K(S^{X_{k}}(s_{1}))_{c}$,

(4) $s_{\lambda}$ is transverse to $\Sigma X_{k}(N,P)$.

\noindent The property (5.1.2) is satisfied for $s_{1}$ by (5.1).
\end{proof}

In what follows we set $d_{1}(s,X_{k})=(s|S^{X_{k}})^{\ast}(\mathbf{d}_{1})$.
We also choose and fix a Riemannian metric of $P$ and identify $Q(S^{X_{k}})$
with the orthogonal complement of Im$(d_{1}(s,X_{k}))$ in $(\pi_{P}%
^{\mathbf{\infty}}\circ s|S^{X_{k}})^{\ast}(TP)$.

\begin{lemma}
Let $k>1$. Let $X_{k}$ be any of $A_{k}$, $D_{k}$ and $E_{k}$. Let $s$ be a
section of $\Gamma_{\Omega}^{tr}(N,P)$ satisfying the property $(5.1.2)$ for
$s$ $($in place of $s_{1}$$)$ of Lemma $5.1$. Then there exists a homotopy
$s_{\lambda}$ relative to $U(C_{k+1})_{2-3r_{0}}$ in $\Gamma_{\Omega}%
^{tr}(N,P)$ with $s_{0}=s$ such that

$(5.2.1)$ $S^{X_{k}}(s_{\lambda})=S^{X_{k}}(s)$ for any $\lambda$,

$(5.2.2)$ $\pi_{P}^{\mathbf{\infty}}\circ s_{1}|S^{X_{k}}(s_{1})$ is an
immersion into $P$ such that $d(\pi_{P}^{\mathbf{\infty}}\circ s_{1}|S^{X_{k}%
}(s_{1})):T(S^{X_{k}}(s_{1}))\linebreak \rightarrow TP$ is equal to $(\pi
_{P}^{\mathbf{\infty}}\circ s_{1})^{TP}\circ d_{1}(s_{1},X_{k})|T(S^{X_{k}%
}(s_{1}))$, where $(\pi_{P}^{\mathbf{\infty}}\circ s_{1})^{TP}:(\pi
_{P}^{\mathbf{\infty}}\circ s_{1})^{\ast}(TP)\rightarrow TP$ is the canonical
induced bundle map,

$(5.2.3)$ we have that $i_{K(S^{X_{k}}(s_{1}))}|_{S^{X_{k}}}=\Psi(s_{1}%
,X_{k})$ and $K(S^{X_{k}}(s_{1}))_{c}\subset\frak{n}(s,X_{k})_{c}$ for any
point $c\in S^{X_{k}}(s_{1})$.

\begin{proof}
Since
\[%
\begin{array}
[c]{ll}%
\mathbf{K}\cap T(\Sigma^{J_{k},0}(N,P))=\{\mathbf{0}\} & \text{ for }%
A_{k}\text{ }(k\geq1),\\
\mathbf{K}\cap T(\Sigma^{\frak{J}_{2},0}(N,P))=\{\mathbf{0}\} & \text{ for
}D_{k}\,(k\geq4),\\
\mathbf{K}\cap T(\Sigma^{\frak{J}_{3},0}(N,P))=\{\mathbf{0}\} & \text{ for
}E_{6}\text{ and }E_{7},\\
\mathbf{K}\cap T(\Sigma^{\frak{J}_{4},0}(N,P))=\{\mathbf{0}\} & \text{ for
}E_{8},\\
&
\end{array}
\]
it follows that $(\pi_{P}^{\mathbf{\infty}}\circ s)^{TP}\circ d_{1}%
(s,X_{k})|T(S^{X_{k}})$ is a monomorphism. By the Hirsch Immersion Theorem
([H1, Theorem 5.7]) there exists a homotopy of monomorphisms $m_{\lambda
}^{\prime}:T(S^{X_{k}})\rightarrow TP$ covering a homotopy $m_{\lambda
}:S^{X_{k}}\rightarrow P$ relative to $U(C_{k+1})_{2-4r_{0}}$\ such that
$m_{0}^{\prime}=(\pi_{P}^{\mathbf{\infty}}\circ s)^{TP}\circ d_{1}%
(s,X_{k})|T(S^{X_{k}})$ and that $m_{1}$ is an immersion with $d(m_{1}%
)=m_{1}^{\prime}$. Then we can extend $m_{\lambda}^{\prime}$ to a homotopy
$\widetilde{m_{\lambda}^{\prime}}:TN|_{S^{X_{k}}}\rightarrow TP$ of
homomorphisms of constant rank $p-1$ relative to $U(C_{k+1})_{2-3r_{0}}$ so
that $\widetilde{m_{0}^{\prime}}=(\pi_{P}^{\mathbf{\infty}}\circ s)^{TP}\circ
d_{1}(s,X_{k})$. In fact, let $m:S^{X_{k}}\times\lbrack0,1]\rightarrow
P\times\lbrack0,1]$ and $m^{\prime}:T(S^{X_{k}})\times\lbrack0,1]\rightarrow
TP\times\lbrack0,1]$ be the maps defined by $m(c,\lambda)=(m_{\lambda
}(c),\lambda)$ and $m^{\prime}(\mathbf{v},\lambda)=(m_{\lambda}^{\prime
}(\mathbf{v}),\lambda)$ respectively. Let $m^{\ast}(m^{\prime}):T(S^{X_{k}%
})\times\lbrack0,1]\rightarrow m^{\ast}(TP\times\lbrack0,1])$ be the canonical
monomorphism induced from $m^{\prime}$ by $m$.\ Let $\mathcal{F}%
_{1}=\mathrm{\operatorname{Im}}(m^{\ast}(m^{\prime}))$ and $\mathcal{F}_{2}$
be the orthogonal complement of $\mathcal{F}_{1}$ in $m^{\ast}(TP\times
\lbrack0,1])$. Since $\mathcal{F}_{2}$ is isomorphic to $(\mathcal{F}%
_{2}|_{S^{X_{k}}\times0})\times\lbrack0,1]$, we obtain a monomorphism of rank
$k-1$%
\[
j_{\mathcal{F}}:\mathrm{\operatorname{Im}}(d_{1}(s,X_{k})|\frak{n}%
(X_{k}))\times\lbrack0,1]\rightarrow\mathcal{F}_{2}\text{ \ \ \ \ \ over
}S^{X_{k}}\times\lbrack0,1]\text{.}%
\]
Since $d_{1}(s,X_{k})|(TN|_{S^{X_{k}}})$ is of constant rank $p-1$ and induces
the homomorphism of kernel rank $n-p+1$%
\[
d:\frak{n}(X_{k})\times\lbrack0,1]\rightarrow\mathrm{\operatorname{Im}}%
(d_{1}(s,X_{k})|\frak{n}(X_{k}))\times\lbrack0,1]\overset{j_{\mathcal{F}}%
}{\rightarrow}\mathcal{F}_{2}\text{,}%
\]
we define $\widetilde{m^{\prime}}$ to be the composition%
\begin{align*}
TN|_{S^{X_{k}}}\times\lbrack0,1]  &  \cong(T(S^{X_{k}})\oplus\frak{n}%
(X_{k}))\times\lbrack0,1]\overset{\underrightarrow{\text{ }m^{\ast}(m^{\prime
})\oplus d\text{ \ }}}{}\mathcal{F}_{1}\oplus\mathcal{F}_{2}\\
&  \rightarrow\mathrm{\operatorname{Im}}(m^{\ast}(m^{\prime}))\oplus
\mathrm{Cok}(m^{\ast}(m^{\prime}))\cong m^{\ast}(TP\times\lbrack
0,1])\overset{\underrightarrow{\text{ }m^{TP\times\lbrack0,1]}\text{\ }}}%
{}TP\times\lbrack0,1],
\end{align*}
where $m^{TP\times\lbrack0,1]}:m^{\ast}(TP\times\lbrack0,1])\rightarrow
TP\times\lbrack0,1]$ is the canonical bundle map. We define $\widetilde
{m_{\lambda}^{\prime}}$ to be $(\widetilde{m_{\lambda}^{\prime}}%
(\mathbf{v}),\lambda)=\widetilde{m^{\prime}}(\mathbf{v},\lambda)$.

Next we construct a homotopy $s_{\lambda}:N\rightarrow\Omega(N,P)$ from
$\widetilde{m_{\lambda}^{\prime}}.$ We write, by $\Sigma^{n-p+1}(N,P)^{\prime
}$, the submanifold of $J^{1}(N,P)=J^{1}(TN,TP)$, which corresponds to
$\Sigma^{n-p+1}(N,P)$ to distinguish them. Then $\pi_{1}^{\mathbf{\infty}%
}|\Sigma X_{k}(N,P):\Sigma X_{k}(N,P)\rightarrow\Sigma^{n-p+1}(N,P)^{\prime}$
becomes a fiber bundle. We regard $\widetilde{m_{\lambda}^{\prime}}$ as a
homotopy $S^{X_{k}}\rightarrow\Sigma^{n-p+1}(N,P)^{\prime}.$\ By the covering
homotopy property to $s|S^{X_{k}}$and $\widetilde{m_{\lambda}^{\prime}}$, we
obtain a homotopy $s_{\lambda}^{\prime}:S^{X_{k}}\rightarrow\Sigma X_{k}(N,P)$
covering $\widetilde{m_{\lambda}^{\prime}}$ relative to $U(C_{k+1})_{2-3r_{0}%
}$ such that $s_{0}^{\prime}=s|S^{X_{k}}$.

By the transversality of $s$, we regard small tubular neighborhoods of
$S^{X_{k}}(s)$ and $\Sigma X_{k}(N,P)$ as vector bundles and that $s$ induces
a bundle map between them when restricted. By applying the homotopy extension
property to $s$ and $s_{\lambda}^{\prime}$, we first extend $s_{\lambda
}^{\prime}$ to a homotopy defined on this tubular neighborhood of $S^{X_{k}}$
and then extend it to a required homotopy $s_{\lambda}\in\Gamma_{\Omega}%
^{tr}(N,P)$, which satisfies $s_{0}=s$, $s_{\lambda}|S^{X_{k}}=s_{\lambda
}^{\prime}$ and $\ s_{\lambda}|U(C_{k+1})_{2-3r_{0}}=s|U(C_{k+1})_{2-3r_{0}}$.
This is a standard argument in topology and the details are left to the reader.
\end{proof}
\end{lemma}

Here we give two lemmas necessary for the proof of Theorem 4.1. Let
$\pi:E\rightarrow S$ be a smooth $(n-p+k)$-dimensional vector bundle with a
metric over a $(p-k)$-dimensional manifold, where $S$ is identified with the
zero-section. Then we can identify $\exp_{E}|D_{\varepsilon}(E)$ with
$id_{D_{\varepsilon}(E)}$.

\begin{lemma}
Let $\pi:E\rightarrow S$ be given as above. Let $f_{i}:E\rightarrow P$
$(i=1,2)$ be $\Omega$-regular maps which have only singularities of types
$A_{j}$, $D_{j}$ and $E_{j}$ with $j\leq k$ and those of type $X_{k}$, where
$X_{k}$ is one of $A_{k}$, $D_{k}$ and $E_{k}$, exactly on $S$\ such that, for
every point $c\in S,$

$(5.3.1)$ $f_{1}|S=f_{2}|S$, which are immersions,

$(5.3.2)$ $S=S^{X_{k}}(j^{\mathbf{\infty}}f_{1})=S^{X_{k}}(j^{\mathbf{\infty}%
}f_{2})$,

$(5.3.3)$ $K(S^{X_{k}}(j^{\mathbf{\infty}}f_{1}))_{c}=K(S^{X_{k}%
}(j^{\mathbf{\infty}}f_{2}))_{c}$ are tangent to $\pi^{-1}(c)$,

$(5.3.4)$ $T_{c}(S^{I_{j-1}}(j^{\mathbf{\infty}}f_{1}))=T_{c}(S^{I_{j-1}%
}(j^{\mathbf{\infty}}f_{2}))$, $((j^{\mathbf{\infty}}f_{1}|S)^{\ast}%
\mathbf{P}_{j})_{c}=((j^{\mathbf{\infty}}f_{2}|S)^{\ast}\mathbf{P}_{j})_{c}%
$\ and
\[
(j^{\mathbf{\infty}}f_{1}|S)^{\ast}(\mathbf{d}_{j+1}^{I}\circ
d(j^{\mathbf{\infty}}f_{1}))_{c}=(j^{\mathbf{\infty}}f_{2}|S)^{\ast
}(\mathbf{d}_{j+1}^{I}\circ d(j^{\mathbf{\infty}}f_{2}))_{c}%
\]
for each number $j$ and $I=J$ or $\frak{J}$,

$(5.3.5)$ if $X_{k}=D_{k}$, then for each number $j\geq3,$%
\[
(j^{\mathbf{\infty}}f_{1}|S)^{\ast}(\mathbf{r}_{j})_{c}=(j^{\mathbf{\infty}%
}f_{2}|S)^{\ast}(\mathbf{r}_{j})_{c}\text{ \ and \ }(j^{\mathbf{\infty}}%
f_{1}|S)^{\ast}(d(\mathbf{r}_{j}))_{c}=(j^{\mathbf{\infty}}f_{2}|S)^{\ast
}(d(\mathbf{r}_{j}))_{c},
\]

$(5.3.6)$ if $X_{k}=E_{7}$, then%
\[
(j^{\mathbf{\infty}}f_{1}|S)^{\ast}(d(\widetilde{\mathbf{d}}_{4}^{\frak{J}%
}|\bigcirc^{4}\mathbf{K}_{3}^{\frak{J}}\mathbf{))}_{c}=(j^{\mathbf{\infty}%
}f_{2}|S)^{\ast}(d(\widetilde{\mathbf{d}}_{4}^{\frak{J}}|\bigcirc
^{4}\mathbf{K}_{3}^{\frak{J}}\mathbf{))}_{c},
\]

$(5.3.7)$ if $X_{k}=E_{8}$, then%
\[%
\begin{array}
[c]{l}%
(j^{\mathbf{\infty}}f_{1}|S)^{\ast}(d(\widetilde{\mathbf{d}}_{4}^{\frak{J}%
}|\bigcirc^{3}\mathbf{K}_{4}^{\frak{J}}\bigcirc(\mathbf{K}_{2}^{\frak{J}%
}/\mathbf{K}_{4}^{\frak{J}}))\mathbf{)}_{c}=(j^{\mathbf{\infty}}f_{2}%
|S)^{\ast}(d(\widetilde{\mathbf{d}}_{4}^{\frak{J}}|\bigcirc^{3}\mathbf{K}%
_{4}^{\frak{J}}\bigcirc(\mathbf{K}_{2}^{\frak{J}}/\mathbf{K}_{4}^{\frak{J}%
})\mathbf{))}_{c},\\
(j^{\mathbf{\infty}}f_{1}|S)^{\ast}(\widetilde{\mathbf{d}}_{5}^{\frak{J}%
}|\bigcirc^{5}\mathbf{K}_{4}^{\frak{J}}\mathbf{)}_{c}=(j^{\mathbf{\infty}%
}f_{2}|S)^{\ast}(\widetilde{\mathbf{d}}_{5}^{\frak{J}}|\bigcirc^{5}%
\mathbf{K}_{4}^{\frak{J}}\mathbf{)}_{c}%
\end{array}
\]

\noindent Let $\eta:S\rightarrow\lbrack0,1]$ be any smooth function. Let
$\varepsilon:S\rightarrow\mathbf{R}$ be a sufficiently small positive smooth
function. For any $c\in S$, $\mathbf{v}_{c}\in\pi^{-1}(c)$ with $\Vert
\mathbf{v}_{c}\Vert\leq\varepsilon(c)$, let $\mathbf{f}^{\eta}(\mathbf{v}%
_{c})$ denote $\exp_{P,f_{1}(c)}((1-\eta(c))\exp_{P,f_{1}(c)}^{-1}%
(f_{1}(\mathbf{v}_{c}))+\eta(c)\exp_{P,f_{2}(c)}^{-1}(f_{2}(\mathbf{v}_{c})))$.

Then the map $\mathbf{f}^{\eta}:D_{\varepsilon}(E)\rightarrow P$ is a
well-defined $\Omega$-regular map such that

$(1)$ $\mathbf{f}^{\eta}|S=f_{1}|S=f_{2}|S,$

$(2)$ $S=S^{X_{k}}(j^{\mathbf{\infty}}\mathbf{f}^{\eta})$,

$(3)$ $K(S^{X_{k}}(j^{\mathbf{\infty}}\mathbf{f}^{\eta})_{c}=K(S^{X_{k}%
}(j^{\mathbf{\infty}}f_{1}))_{c}=K(S^{X_{k}}(j^{\mathbf{\infty}}f_{2}))_{c}$
are tangent to $\pi^{-1}(c)$,

$(4)$ $T_{c}(S^{I_{j-1}}(j^{\mathbf{\infty}}\mathbf{f}^{\eta}))=T_{c}%
(S^{I_{j-1}}(j^{\mathbf{\infty}}f_{1}))$, $((j^{\mathbf{\infty}}%
\mathbf{f}^{\eta}|S)^{\ast}\mathbf{P}_{j})_{c}=((j^{\mathbf{\infty}}%
f_{1}|S)^{\ast}\mathbf{P}_{j})_{c}$\ and
\[
(j^{\mathbf{\infty}}\mathbf{f}^{\eta})^{\ast}(\mathbf{d}_{j+1}^{I}\circ
d(j^{\mathbf{\infty}}\mathbf{f}^{\eta}))_{c}=(j^{\mathbf{\infty}}f_{1})^{\ast
}(\mathbf{d}_{j+1}^{I}\circ d(j^{\mathbf{\infty}}f_{1}))_{c}%
=(j^{\mathbf{\infty}}f_{2})^{\ast}(\mathbf{d}_{j+1}^{I}\circ
d(j^{\mathbf{\infty}}f_{2}))_{c},
\]
for each number $j$ and $I=J$ or $\frak{J}$,

$(5)$ if $X_{k}=D_{k}$, then, for each number $j\geq3,$%
\[%
\begin{array}
[c]{l}%
(j^{\mathbf{\infty}}\mathbf{f}^{\eta})^{\ast}(\mathbf{r}_{j})_{c}%
=(j^{\mathbf{\infty}}f_{1})^{\ast}(\mathbf{r}_{j})_{c}=(j^{\mathbf{\infty}%
}f_{2})^{\ast}(\mathbf{r}_{j})_{c},\\
(j^{\mathbf{\infty}}\mathbf{f}^{\eta})^{\ast}(d(\mathbf{r}_{j}\mathbf{))}%
_{c}=(j^{\mathbf{\infty}}f_{1})^{\ast}(d(\mathbf{r}_{j}\mathbf{))}%
_{c}=(j^{\mathbf{\infty}}f_{2})^{\ast}(d(\mathbf{r}_{j}\mathbf{))}_{c},
\end{array}
\]

$(6)$ if $X_{k}=E_{7},$ then%
\[
(j^{\mathbf{\infty}}\mathbf{f}^{\eta})^{\ast}(d(\widetilde{\mathbf{d}}%
_{4}^{\frak{J}}|\bigcirc^{4}\mathbf{K}_{3}^{\frak{J}}\mathbf{))}%
_{c}=(j^{\mathbf{\infty}}f_{1})^{\ast}(d(\widetilde{\mathbf{d}}_{4}^{\frak{J}%
}|\bigcirc^{4}\mathbf{K}_{3}^{\frak{J}}\mathbf{))}_{c}=(j^{\mathbf{\infty}%
}f_{2})^{\ast}(d(\widetilde{\mathbf{d}}_{4}^{\frak{J}}|\bigcirc^{4}%
\mathbf{K}_{3}^{\frak{J}}\mathbf{))}_{c},
\]

$(7)$ if $X_{k}=E_{8},$ then
\begin{align*}
&  (j^{\mathbf{\infty}}\mathbf{f}^{\eta})^{\ast}(d(\widetilde{\mathbf{d}}%
_{4}^{\frak{J}}|\bigcirc^{3}\mathbf{K}_{4}^{\frak{J}}\bigcirc(\mathbf{K}%
_{2}^{\frak{J}}/\mathbf{K}_{4}^{\frak{J}}))\mathbf{)}_{c}\\
&  =(j^{\mathbf{\infty}}f_{1})^{\ast}(d(\widetilde{\mathbf{d}}_{4}^{\frak{J}%
}|\bigcirc^{3}\mathbf{K}_{4}^{\frak{J}}\bigcirc(\mathbf{K}_{2}^{\frak{J}%
}/\mathbf{K}_{4}^{\frak{J}}))\mathbf{)}_{c}=(j^{\mathbf{\infty}}f_{2})^{\ast
}(d(\widetilde{\mathbf{d}}_{4}^{\frak{J}}|\bigcirc^{3}\mathbf{K}_{4}%
^{\frak{J}}\bigcirc(\mathbf{K}_{2}^{\frak{J}}/\mathbf{K}_{4}^{\frak{J}%
})\mathbf{))}_{c},\\
&  (j^{\mathbf{\infty}}\mathbf{f}^{\eta})^{\ast}(\widetilde{\mathbf{d}}%
_{5}^{\frak{J}}|\bigcirc^{5}\mathbf{K}_{4}^{\frak{J}})_{c}=(j^{\mathbf{\infty
}}f_{1})^{\ast}(\widetilde{\mathbf{d}}_{5}^{\frak{J}}|\bigcirc^{5}%
\mathbf{K}_{4}^{\frak{J}}\mathbf{)}_{c}=(j^{\mathbf{\infty}}f_{2})^{\ast
}(\widetilde{\mathbf{d}}_{5}^{\frak{J}}|\bigcirc^{5}\mathbf{K}_{4}^{\frak{J}%
}\mathbf{)}_{c}.
\end{align*}

\begin{proof}
Let us take a Riemannian metric on $E$ which is compatible with the metric of
the vector bundle $E$\ over $S$. In particular, $S$\ is a Riemannian
submanifold\ of $E$. Furthermore, take a Riemannian metric on $P$\ such that
$f(S)\cap P$\ is a Riemannian submanifold\ of $P$\ around $f(c)$. Then the
local coordinates of $\exp_{N,c}(\mathbf{K}_{c})$\ and $\exp_{P,f(c)}%
(\mathbf{Q}_{c})$\ are independent of the coordinates of $S$, where
$\mathbf{Q}_{c}$\ is regarded as a line subspace of $T_{f(c)}P$.\ 

Since $j^{\mathbf{\infty}}f_{i}$ are transverse to $\Sigma X_{k}(E,P)$, we
have $\frak{n}(j^{\mathbf{\infty}}f_{1},I_{j}\subset I_{j-1})_{c}%
=\frak{n}(j^{\mathbf{\infty}}f_{2},I_{j}\subset I_{j-1})_{c}$ for $I=J$ or
$\frak{J}$. Furthermore, it follows similarly that $\frak{n}(j^{\mathbf{\infty
}}f_{1},\overline{D}_{j+1}\subset\overline{D}_{j})_{c}=\frak{n}%
(j^{\mathbf{\infty}}f_{2},\overline{D}_{j+1}\subset\overline{D}_{j})_{c}$ for
$X_{k}=D_{k}$ and $\frak{n}(j^{\mathbf{\infty}}f_{1},\overline{E}_{j+1}%
\subset\overline{E}_{j})_{c}=\frak{n}(j^{\mathbf{\infty}}f_{2},\overline
{E}_{j+1}\subset\overline{E}_{j})_{c}$ for $X_{k}=E_{k}$.

We may consider $\eta(c)$ as a constant when dealing with higher intrinsic
derivatives in the lemma by the identification (1.2) and the property of the
total tangent bundle $\mathbf{D}$ given in the beginning\ of Section 2. Then
the assertions (1)-(7) follow from the assumptions and the properties of
$\Sigma X_{k}(N,P)$.

Let $\varepsilon$ be sufficiently small. Then since $\Omega(n,p)$ is open and
$j^{\mathbf{\infty}}\mathbf{f}^{\eta}(S)\subset\Sigma X_{k}(N,P)$,
$\mathbf{f}^{\eta}$ is an $\Omega$-regular map. Since $\mathbf{f}^{\eta}$ is
transverse to $\Sigma X_{k}(N,P)$, we have $S^{X_{k}}(j^{\mathbf{\infty}%
}\mathbf{f}^{\eta})=S$.
\end{proof}
\end{lemma}

The proof of the following lemma is elementary, and so is left to the reader.

\begin{lemma}
Let $\pi:E\rightarrow S$ be given as above, and let $(\Omega,\Sigma)$ be a
pair of a smooth manifold and its submanifold of codimension $n-p+k$. Let
$\varepsilon:S\rightarrow\mathbf{R}$ be a sufficiently small positive smooth
function. Let $h:D_{\varepsilon}(E)\rightarrow(\Omega,\Sigma)$ be a smooth map
such that $S=h^{-1}(\Sigma)$ and that $h$ is transverse to $\Sigma$. Then
there exists a smooth homotopy $h_{\lambda}:(D_{\varepsilon}(E),S)\rightarrow
(\Omega,\Sigma)$ between $h$ and $\exp_{\Omega}\circ dh|D_{\varepsilon}(E)$
such that

$(1)$ $h_{\lambda}|S=h_{0}|S$, $S=h_{\lambda}^{-1}(\Sigma)=h_{0}^{-1}(\Sigma)$
for any $\lambda$,

$(2)$ $h_{\lambda}$ is smooth and is transverse to $\Sigma$ for any $\lambda$,

$(3)$ $h_{0}=h$ and $h_{1}(\mathbf{v}_{c})=\exp_{\Omega,h(c)}\circ
dh(\mathbf{v}_{c})$ for $c\in S$ and $\mathbf{v}_{c}\in D_{\varepsilon}(E_{c}).$
\end{lemma}

\section{Proof of Theorem 4.1}

For the normal bundles $\frak{n}(X_{k})$ ($=\frak{n}(s,X_{k})$) and $Q$
$(=Q(S^{X_{k}}))$\ over $S^{X_{k}}(s)$, we recall that Hom$(\Sigma_{j=1}%
^{\ell}\bigcirc^{j}\frak{n}(X_{k}),Q)|_{S^{X_{k}}}$ is identified with the set
of polynomials of degree $\leq\ell$ having the constant $0$ with coefficients
depending on a point of $S^{X_{k}}(s)$ (see [Mats, Ch. III, Section 2]).\ For
a point $c\in S^{X_{k}}(s)$, take an open neighborhood $U$ around $c$ such
that $\frak{n}(X_{k})|_{U}$ and $Q|_{U}$ are the trivial bundles, say
$U\times\mathbf{R}^{n-p+k}$ and $U\times\mathbf{R}$ respectively, where
$\mathbf{R}^{n-p+k}$ has coordinates $(x_{1},\ldots,x_{n-p+k})$ and
$\mathbf{R}$ has $y$. Then an element of Hom$(\bigcirc^{j}\frak{n}%
(X_{k}),Q)|_{U}$ is identified with a polynomial $y(c)=\Sigma_{|\omega
|=j}a^{\omega}(c)x_{1}^{\omega_{1}}x_{2}^{\omega_{2}}\cdots x_{n-p+k}%
^{\omega_{n-p+k}}$, $c\in U$, where $\omega=(\omega_{1},\omega_{2}%
,\cdots,\omega_{n-p+k})$, $\omega_{i}\geq0$ ($i=1,\cdots,n-p+k$),
$|\omega|=\omega_{1}+\cdots+\omega_{n-p+k}$ and $a^{\omega}(c)$ is a real
number. If $a^{\omega}(c)$ are smooth functions of $c,$ then $\{a^{\omega
}(c)\}$ defines a smooth section of Hom$(\bigcirc^{j}\frak{n}(X_{k}%
),Q)|_{S^{X_{k}}}$ over $U$.

We first introduce several homomorphisms between vector bundles over
$S^{X_{k}}(s)$, which are used for the construction of the required $\Omega
$-regular map in Theorem 4.1.

Let $s\in\Gamma_{\Omega}^{tr}(N,P)$. By deforming $s$ if necessary, we may
assume without loss of generality that $s$ satisfies (5.1.2) of Lemma 5.1 and
(5.2.2) of Lemma 5.2, where $s_{1}$ is replaced by $s$.

In the following, let $K=K(S^{X_{k}})$, $Q=Q(S^{X_{k}})$, $K_{j}$ ($j\geq2$)
refer to $K_{j}^{J}(S^{X_{k}})$ for $X_{k}=A_{k}$ and $K_{j}^{\frak{J}%
}(S^{X_{k}})$ for $X_{k}=D_{k}$ or $E_{k}$, and let $L=(s|S^{X_{k}})^{\ast
}(\mathbf{L})$. We now describe the isomorphisms%
\begin{equation}%
\begin{array}
[c]{ll}%
(s|S^{I_{j}})^{\ast}(\mathbf{d}_{j+1}^{I}\circ ds|\frak{n}(s,I_{j}\subset
I_{j-1})) & \text{for }I=J\text{ or }\frak{J,}\\
(s|S^{\overline{D}_{j+1}})^{\ast}(d(\mathbf{r}_{j-1})\circ ds|\frak{n}%
(s,\overline{D}_{j+1}\subset\overline{D}_{j})), & \\
(s|S^{\overline{E}_{7}})^{\ast}(d(\widetilde{\mathbf{d}}_{4}^{\frak{J}%
}|\bigcirc^{4}\mathbf{K}_{3}^{\frak{J}})\circ ds|\frak{n}(s,\overline{E}%
_{7}\subset\overline{E}_{6})), & \\
(s|S^{E_{8}})^{\ast}(d(\widetilde{\mathbf{d}}_{4}^{\frak{J}}|\bigcirc
^{3}\mathbf{K}_{3}^{\frak{J}}\bigcirc(\mathbf{K}_{2}^{\frak{J}}/\mathbf{K}%
_{4}^{\frak{J}}))\circ ds|\frak{n}(s,E_{8}\subset\overline{E}_{7})), &
\end{array}
\end{equation}
which are induced by $s|S^{I_{j}}$, $s|S^{\overline{D}_{j+1}}$,
$s|S^{\overline{E}_{7}}$ and $s|S^{E_{8}}$ respectively. They yield the
decompositions of the target bundles derived from $\mathbf{K\supset K}_{2}%
^{J}\supset\mathbf{L}$ and $\mathbf{K}_{2}^{\frak{J}}\supset\mathbf{K}%
_{3}^{\frak{J}}$, and hence the decomposition of the normal bundles as follows.

For $A_{k}$,%
\begin{equation}%
\begin{array}
[c]{ll}%
\frak{n}(s,J_{1}\subset J_{0})=K/K_{2}\oplus T_{2}^{A}\rightarrow
\text{Hom}(K/K_{2}\oplus K_{2},Q)\text{ \ by }(2.1), & \\
\frak{n}(s,J_{j}\subset J_{j-1})=T_{j+1}^{A}\rightarrow\text{Hom}(\bigcirc
^{j}K_{2},Q)\text{ }(2\leq j\leq k-1)\text{ \ by }(2.2), & \text{ }\\
\frak{n}(s,J_{k}\subset J_{k-1})=K_{2}\rightarrow\text{Hom}(\bigcirc^{k}%
K_{2},Q)\text{ \ by (5), (8) in Section 2 and }(2.2), &
\end{array}
\end{equation}
over $S^{A_{k}}$, where (i) $K/K_{2}$\ refers to the orthogonal complements of
$K_{2}$ in $K,$ and (ii) $T_{2}^{A}$ refers to the $1$-subbundle of
$\frak{n}(s,J_{1}\subset J_{0})$, which corresponds to the direct summand in
the right-hand side.

For $D_{k}$,%
\begin{equation}%
\begin{array}
[c]{l}%
\frak{n}(s,\frak{J}_{1}\subset\frak{J}_{0})=K/K_{2}\oplus T_{2}^{D_{4}%
}\rightarrow\text{Hom}(K/K_{2}\oplus K_{2},Q)\text{ (}k=4\text{) by }(2.1),\\
\frak{n}(s,\frak{J}_{1}\subset\frak{J}_{0})=K/K_{2}\oplus T_{2,1}^{D_{k}%
}\oplus T_{2,2}^{D_{k}}\rightarrow\text{Hom}(K/K_{2}\oplus K_{2}/L\oplus
L,Q)\text{ (}k\geq5\text{) by }(2.1),\\
\frak{n}(s,\frak{J}_{2}\subset\frak{J}_{1})=K_{2}\oplus T_{3}^{D_{4}%
}\rightarrow\text{Hom}(\mathcal{P}^{D_{4}},Q)\text{ (}k=4\text{) \ by
}(2.3)\text{, (D-i), (D-iv-}3),\\
\frak{n}(s,\frak{J}_{2}\subset\frak{J}_{1})=K_{2}/L\oplus L\oplus T_{3}%
^{D_{k}}\rightarrow\text{Hom}(\mathcal{P}^{D},Q)\text{ (}k\geq5\text{) by
}(2.3)\text{, (D-ii), (D-iv-}3\text{),}\\
\frak{n}(s,\overline{D}_{j}\subset\overline{D}_{j-1})=T_{j-1}^{D_{k}%
}\rightarrow\text{Hom}(\bigcirc^{j-2}L,Q)\text{ (}5\leq j\leq k\text{)\ by
}(3.2),
\end{array}
\end{equation}
over $S^{D_{k}}$, where

(i) $K/K_{2}$ and $K_{2}/L$\ refer to the orthogonal complements of $K_{2}$ in
$K$, $L$ in $K_{2}$ respectively,

(ii) $T_{2}^{D_{4}}$, $T_{2,1}^{D_{k}}$, $T_{2,2}^{D_{k}}$ and $T_{3}^{D_{k}}$
refer to the subbundles of the normal bundles concerned which corresponds to
the direct summands in the right-hand side respectively (dim$T_{2}^{D_{4}}=2$
and dim$T_{3}^{D_{4}}=1$),

(iii) $\mathcal{P}^{D_{4}}=\bigcirc^{2}K_{2}=V_{1}\oplus V_{2}$ and
$\mathcal{P}^{D}=L\bigcirc K_{2}/L\oplus\bigcirc^{2}K_{2}/L\oplus\bigcirc
^{2}L$. Here, $V_{2}$ is the $1$-subbundle of $\bigcirc^{2}K_{2}$ which
consists of all elements $v$ with $(s|S_{3}^{D_{4}})^{\ast}(\widetilde
{\mathbf{d}}_{3}^{\frak{J}})|(K_{2}\bigcirc\{v\})=0$, and $V_{1}$ is the
orthogonal complement of $V_{2}$ in $\bigcirc^{2}K_{2}$.

For $E_{k}$,%
\begin{equation}%
\begin{array}
[c]{l}%
\frak{n}(s,\frak{J}_{1}\subset\frak{J}_{0})=K/K_{2}\oplus T_{2,1}^{E}\oplus
T_{2,2}^{E}\rightarrow\text{Hom}(K/K_{2}\oplus K_{2}/K_{3}\oplus
K_{3},Q)\text{ \ \ by }(2.1),\\
\frak{n}(s,\frak{J}_{2}\subset\frak{J}_{1})=K_{2}/K_{3}\oplus T_{3,1}%
^{E}\oplus T_{3,2}^{E}\rightarrow\text{Hom}(\bigcirc^{2}K_{2}/K_{3}%
\oplus\bigcirc^{2}K_{3}\oplus K_{3}\bigcirc K_{2}/K_{3},Q)\\
\text{ }%
\ \ \ \ \ \ \ \ \ \ \ \ \ \ \ \ \ \ \ \ \ \ \ \ \ \ \ \ \ \ \ \ \ \ \ \ \ \ \ \ \ \ \ \ \ \ \ \ \ \ \ \ \ \ \ \ \ \ \ \ \ \ \ \ \ \ \ \text{by
\ }(2.3),\text{ }(2.4),\\
\frak{n}(s,E_{6}\subset\frak{J}_{2})=K_{3}\oplus T_{4}^{E_{6}}\rightarrow
\text{Hom}(\mathcal{P}^{E},Q)\text{ \ \ for }E_{6}\text{ by }(2.5)\text{,
(E-i), (E-iii-a)},\\
\frak{n}(s,\frak{J}_{3}\subset\frak{J}_{2})=T_{4}^{E_{7}}\oplus K_{3}%
\rightarrow\text{Hom}(\mathcal{P}^{E},Q)\text{ \ \ for }E_{7}\text{ by
}(2.5)\text{, (E-i)},\text{ (E-iii-b)},\\
\frak{n}(s,\frak{J}_{3}\subset\frak{J}_{2})=T_{4,1}^{E_{8}}\oplus
T_{4,2}^{E_{8}}\rightarrow\text{Hom}(\mathcal{P}^{E},Q)\text{ \ \ for }%
E_{8}\text{ by }(2.5)\text{, (E-iii)},\\
\frak{n}(s,E_{7}\subset\overline{E}_{6})=T_{5}^{E_{7}}\rightarrow
\text{Hom}(\bigcirc^{4}K_{3},Q)\text{ \ \ for }E_{7}\text{ by }(3.3)\text{,
(E-iii-b)},\\
\frak{n}(s,\overline{E}_{7}\subset\overline{E}_{6})|_{S^{E_{8}}}%
=K_{4}|_{S^{E_{8}}}\rightarrow\text{Hom}(\bigcirc^{4}K_{4},Q)|_{S^{E_{8}}%
}\text{ \ \ for }E_{8}\text{\ by (E-iii-c),}.\\
\frak{n}(s,E_{8}\subset\overline{E}_{7})=T_{5}^{E_{8}}\rightarrow
\text{Hom}(\bigcirc^{3}K_{4}\bigcirc K_{2}/K_{4},Q)\text{ \ \ for }E_{8}\text{
by }(\frak{J}E\text{-}2),\text{ (E-iii)},
\end{array}
\end{equation}
over $S^{E_{k}}$, where

(i) $\frak{n}(s,E_{6}\subset\frak{J}_{2})=\frak{n}(s,\frak{J}_{3}%
\subset\frak{J}_{2})|_{S^{E_{6}}}$ and $\mathcal{P}^{E}=\bigcirc^{3}%
K_{3}\oplus\bigcirc^{2}K_{3}\bigcirc K_{2}/K_{3},$

(ii) $K/K_{2}$, $K_{2}/L$, and $K_{2}/K_{3}$ or $K_{2}/K_{4}$\ refer to the
orthogonal complements of $K_{2}$ in $K$, $L$ in $K_{2}$, and $K_{3}$ or
$K_{4}$\ in $K_{2}$ respectively,

(iii) $T_{j,1}^{E}$, $T_{j,2}^{E}$, $T_{4}^{E_{6}}$, $T_{4}^{E_{7}}$,
$T_{4,j}^{E_{8}}$, $T_{5}^{E_{7}}$ and $T_{5}^{E_{8}}$ refer to $1$-subbundle
of the normal bundles concerned, which correspond to the direct summands in
the right-hand side respectively.

The isomorphisms in (6.2), (6.3) and (6.4) canonically induce the following
homomorphisms over $S^{X_{k}}(s)$ respectively.

For $A_{k}$,%
\[%
\begin{array}
[c]{ll}%
\widetilde{d}_{2}^{A}(s):\bigcirc^{2}K/K_{2}\oplus K_{2}\bigcirc T_{2}%
^{A}\rightarrow Q, & \\
\widetilde{d}_{j+1}^{A}(s):\bigcirc^{j}K_{2}\bigcirc T_{j+1}^{A}\rightarrow
Q\text{ }(2\leq j\leq k-1), & \\
\widetilde{d}_{k+1}^{A}(s):\bigcirc^{k+1}K_{2}\rightarrow Q. &
\end{array}
\]

For $D_{k}$,%
\[%
\begin{array}
[c]{ll}%
\widetilde{d}_{2}^{D_{4}}(s):\bigcirc^{2}K/K_{2}\oplus K_{2}\bigcirc
T_{2}^{D_{4}}\rightarrow Q & (k=4),\\
\widetilde{d}_{2}^{D_{k}}(s):\bigcirc^{2}K/K_{2}\oplus K_{2}/L\bigcirc
T_{2,1}^{D_{k}}\oplus L\bigcirc T_{2,2}^{D_{k}}\rightarrow Q & (k\geq5),\\
\widetilde{d}_{3}^{D_{4}}(s):V_{1}\bigcirc K_{2}\oplus V_{2}\bigcirc
T_{3}^{D_{4}}\rightarrow Q & (k=4),\\
\widetilde{d}_{3}^{D_{k}}(s):\bigcirc^{2}K_{2}/L\bigcirc L\oplus\bigcirc
^{2}L\bigcirc T_{3}^{D_{k}}\rightarrow Q & \text{(}k\geq5\text{),}\\
\widetilde{d}(r_{j-2})(s):\bigcirc^{j-2}L\bigcirc T_{j-1}^{D_{k}}\rightarrow
Q\text{ }(5\leq j\leq k) & \text{(}k\geq5\text{).}%
\end{array}
\]
Furthermore, the homomorphism $\mathbf{r}_{j-1}$ induces the homomorphism by
$s|S^{D_{k}}$%
\[
r_{j-1}(s):\bigcirc^{j-1}L\rightarrow Q\text{ }(5\leq j\leq k)\text{\ \ \ over
}S^{D_{k}}(s)\text{ (}k\geq5\text{).}%
\]

For $E_{k}$,%
\[%
\begin{array}
[c]{ll}%
\widetilde{d}_{2}^{E}(s):\bigcirc^{2}K/K_{2}\oplus K_{2}/K_{3}\bigcirc
T_{2,1}^{E}\oplus K_{3}\bigcirc T_{2,2}^{E}\rightarrow Q, & \\
\widetilde{d}_{3}^{E}(s):\bigcirc^{3}K_{2}/K_{3}\oplus\bigcirc^{2}%
K_{3}\bigcirc T_{3,1}^{E}\oplus K_{3}\bigcirc K_{2}/K_{3}\bigcirc T_{3,2}%
^{E}\rightarrow Q\text{,} & \\
\widetilde{d}_{4}^{E_{6}}(s):\bigcirc^{4}K_{3}\oplus\bigcirc^{2}K_{3}\bigcirc
K_{2}/K_{3}\bigcirc T_{4}^{E_{6}}\rightarrow Q & \text{for }E_{6}\text{,}\\
\widetilde{d}_{4}^{E_{7}}(s):\bigcirc^{3}K_{3}\bigcirc T_{4}^{E_{7}}%
\oplus\bigcirc^{3}K_{3}\bigcirc K_{2}/K_{3}\rightarrow Q & \text{for }%
E_{7}\text{,}\\
\widetilde{d}_{4}^{E_{8}}(s):\bigcirc^{3}K_{3}\bigcirc T_{4,1}^{E_{8}}%
\oplus\bigcirc^{2}K_{3}\bigcirc K_{2}/K_{3}\bigcirc T_{4,2}^{E_{8}}\rightarrow
Q & \text{for }E_{8}\text{,}\\
\widetilde{d}_{5}^{E_{7}}(s):\bigcirc^{4}K_{3}\bigcirc T_{5}^{E_{7}%
}\rightarrow Q & \text{for }E_{7}\text{,}\\
\widetilde{d}_{5}^{E_{8}}(s):\bigcirc^{5}K_{4}\oplus\bigcirc^{3}K_{4}\bigcirc
K_{2}/K_{4}\bigcirc T_{5}^{E_{8}}\rightarrow Q & \text{for }E_{8}\text{,}%
\end{array}
\]
where $\widetilde{d}_{5}^{E_{8}}(s)|\bigcirc^{5}K_{4}$ comes from
$(s|S^{E_{8}})^{\ast}(\widetilde{\mathbf{d}}_{5}^{\frak{J}}|\bigcirc
^{5}\mathbf{K}_{4}^{\frak{J}})$.

We define the sections of Hom$(\Sigma_{j=1}^{k+1}\bigcirc^{j}\frak{n}%
(X_{k}),Q)$%
\[%
\begin{array}
[c]{ll}%
q^{A_{k}}(s)=\widetilde{d}_{2}^{A}(s)+\Sigma_{j=3}^{k}\widetilde{d}_{j}%
^{A}(s)+\widetilde{d}_{k+1}^{A}(s) & \text{over }S^{A_{k}}(s),\\
q^{D_{4}}(s)=\widetilde{d}_{2}^{D_{4}}(s)+\widetilde{d}_{3}^{D_{4}}(s) &
\text{over }S^{D_{4}}(s),\\
q^{D_{k}}(s)=\widetilde{d}_{2}^{D_{k}}(s)+\widetilde{d}_{3}^{D_{k}}%
(s)+\Sigma_{j=5}^{k}\widetilde{d}(r_{j-2})(s)+r_{k-1}(s)\text{ }(k\geq5) &
\text{over }S^{D_{k}}(s),\\
q^{E_{6}}(s)=\widetilde{d}_{2}^{E}(s)+\widetilde{d}_{3}^{E}(s)+\widetilde
{d}_{4}^{E_{6}}(s) & \text{over }S^{E_{6}}(s),\\
q^{E_{7}}(s)=\widetilde{d}_{2}^{E}(s)+\widetilde{d}_{3}^{E}(s)+\widetilde
{d}_{4}^{E_{7}}(s)+\widetilde{d}_{5}^{E_{7}}(s) & \text{over }S^{S_{7}}(s),\\
q^{E_{8}}(s)=\widetilde{d}_{2}^{E}(s)+\widetilde{d}_{3}^{E}(s)+\widetilde
{d}_{4}^{E_{8}}(s)+\widetilde{d}_{5}^{E_{8}}(s) & \text{over }S^{E_{8}}(s).
\end{array}
\]

Then we obtain the smooth fiber map
\begin{align}
(\pi_{P}^{\mathbf{\infty}}\circ s|S^{X_{k}})^{TP}\circ(d_{1}(s,X_{k}%
)|\frak{n}(s,X_{k})+q^{X_{k}}(s))  &  :\nonumber\\
\frak{n}(s,X_{k})  &  \rightarrow(\pi_{P}^{\mathbf{\infty}}\circ s|S^{X_{k}%
})^{\ast}(TP)\rightarrow TP
\end{align}
covering the immersion $\pi_{P}^{\mathbf{\infty}}\circ s|S^{X_{k}}:S^{X_{k}%
}(s)\rightarrow P.$

\begin{remark}
We explain what $V_{2}$ is and how the normal form for $D_{4}$ in Introduction
is induced from $\widetilde{d}_{3}^{D_{4}}(s)$. For $c\in S^{D_{4}}(s)$ we set
$\widetilde{d}_{3}=\widetilde{d}_{3}^{D_{4}}(s)|\bigcirc^{3}K_{2,c}=a_{0}%
u^{3}+a_{1}u^{2}\ell+a_{2}u\ell^{2}+a_{3}\ell^{3}$. A generator of
$\mathrm{Hom}((V_{2})_{c},Q_{c})$ is a nonsingular quadratic form of
$\mathrm{Hom}(\bigcirc^{2}K_{2,c},Q_{c})$. Suppose to the contrary. Then there
are coordinates $(u,\ell)$ around $c$ such that $(V_{2})_{c}$ is generated by
$\partial/\partial u\bigcirc\partial/\partial u=\partial^{2}/\partial u^{2}$.
It follows from $(s|S^{D_{4}})^{\ast}(\widetilde{\mathbf{d}}_{3}^{\frak{J}%
})|(K_{2}\bigcirc V_{2})=0$ that the coefficients of $u^{3}$ and $u^{2}\ell
\ $in the polynomial $\widetilde{d}_{3}$ vanish. Namely, $\widetilde{d}_{3}$
is written as $a_{2}u\ell^{2}+a_{3}\ell^{3}=(a_{2}u+a_{3}\ell)\ell^{2}$. This
is of type $S_{5}$ or $S_{E}$. This is a contradiction.

Suppose that $\partial^{2}/\partial u^{2}\pm\partial^{2}/\partial\ell^{2}$ is
the generator of $(V_{2})_{c}$. Then $(s|S_{3}^{D_{4}})^{\ast}(\widetilde
{\mathbf{d}}_{3}^{\frak{J}})|(K_{2}\bigcirc V_{2})=0$ implies that
\[
\left(  \frac{\partial^{3}}{\partial u^{3}}\pm\frac{\partial^{3}}{\partial
u\partial\ell^{2}}\right)  \widetilde{d}_{3}=0\text{ \ \ and \ \ }\left(
\frac{\partial^{3}}{\partial u^{2}\partial\ell}\pm\frac{\partial^{3}}%
{\partial\ell^{3}}\right)  \widetilde{d}_{3}=0.
\]
Then we have $3a_{0}\pm a_{2}=0$ and $a_{1}\pm3a_{3}=0$. Hence, if $a_{0}=0$,
then we have $\widetilde{d}_{3}=\mp3a_{3}u^{2}\ell+a_{3}\ell^{3}$ $(a_{3}%
\neq0)$. If $a_{0}\neq0$, then we may assume $a_{0}=1$ and%
\[
\widetilde{d}_{3}=u^{3}\mp3a_{3}u^{2}\ell\mp3u\ell^{2}+a_{3}\ell^{3}.
\]
Let $k$ be a real number such that $k^{3}\mp3a_{3}k^{2}\mp3k+a_{3}=0$, and let
$u=v+k\ell$. Then we have%
\[
u^{3}\mp3a_{3}u^{2}\ell\mp3u\ell^{2}+a_{3}\ell^{3}=v(v^{2}+3(k\mp a_{3}%
)v\ell+3(k^{2}\mp2a_{3}k\mp1)\ell^{2}).
\]
If we set $A=(k^{2}\mp2a_{3}k\mp1)$, then $A\neq0$ and $4A\neq3(k\mp
a_{3})^{2}$. In fact, if $A\neq0$, then we have%
\[
v^{2}+3(k\mp a_{3})v\ell+3A\ell^{2}=(1-\frac{3(k\mp a_{3})^{2}}{4A}%
)v^{2}+3A(\ell+\frac{k\mp a_{3}}{2A}v)^{2}.
\]
If $4A=3(k\mp a_{3})^{2}$, then $\widetilde{d}_{3}$ is of type $S_{5}$,
Otherwise this is of type $S_{4}^{\pm}$. If $A=0$, then $a_{3}k^{2}+2k\mp
a_{3}=0$, and hence $k=0$ or $a_{3}^{2}=1$. Since $A=0$, the case $k=0$ does
not occur. If $a_{3}^{2}=1$, then $\widetilde{d}_{3}=(u\mp\ell)^{3}$, which is
of type $S_{E}$.
\end{remark}

\begin{proof}
[Proof of Theorem 4.1]By Lemmas 5.1 and 5.2 we may assume that $s$ satisfies
(5.1.2) and (5.2.2), where $s_{1}$ is replaced by $s$. We define
$\mathbf{E}(k)$ to be the union of all $\exp_{N}(D_{\delta\circ s}%
(\frak{n}(X_{k})))$, where $\delta$ is a sufficiently small positive smooth
function defined on $\Sigma X_{k}(N,P)$ for $X_{k}=A_{k}$, $D_{k}$ and $E_{k}$
such that

(i) $\delta\circ s|(S^{X_{k}}\setminus$Int$U(C_{k+1}){_{2}})$ is constant,

(ii) $\exp_{N}(D_{\delta\circ s}(\frak{n}(X_{k})))\setminus$Int$U(C_{k+1}%
){_{2}}$ for $X_{k}=A_{k}$, $D_{k}$ and $E_{k}$ do not intersect mutually.

\noindent This is a tubular neighborhood of $S^{X_{k}\text{'}}$s.

It is enough for the proof of Theorem 4.1 except for (4.1.3) to prove the
following assertion:

(\textbf{A}) {There exists a homotopy $H_{\lambda}$ relative to $U(C_{k+1}%
)_{2-r_{0}}$ in $\Gamma_{\Omega}^{tr}(N,P)$ with $H_{0}=s$ satisfying the
following for each }$X_{k}=A_{k}$, $D_{k}$ and $E_{k}${. }

{$(1)$ $S{^{X_{k}}}(H_{\lambda})=S^{X_{k}}$ for any $\lambda$. }

{$(2)$ We have an }$\Omega$-regular{ map $G$ defined on a neighborhood of
${U(C_{k+1})}_{2-r_{0}}\cup\mathbf{E}(k)$ to $P$ such that $j^{\mathbf{\infty
}}G=H_{1}$ on }${U(C_{k+1})}_{2-r_{0}}\cup\mathbf{E}(k).$

By the Riemannian metric on $P$, we identify $Q$ with the orthogonal line
bundle of Im$(d_{1}(s,X_{k}))$ in $(\pi_{P}^{\mathbf{\infty}}\circ s|S^{X_{k}%
})^{\ast}(TP)$. Then $\exp_{P}\circ(\pi_{P}^{\mathbf{\infty}}\circ s|S^{X_{k}%
})^{TP}|D_{\gamma}(Q)$ is an immersion for some small positive function
$\gamma$. In the proof we express a point of $\mathbf{E}(k)$ as $\mathbf{v}%
_{c}$, where $c\in S^{X_{k}},$ $\mathbf{v}_{c}\in\frak{n}(X_{k})_{c}$ and
$\Vert\mathbf{v}_{c}\Vert\leq\delta(s(c))$. In the proof we say that a smooth
homotopy
\[
k_{\lambda}:(\mathbf{E}(k),\partial\mathbf{E}(k))\rightarrow(\Omega
(N,P),\Omega(N,P)\setminus\Sigma X_{k}(N,P))
\]
has the property (C) if it satisfies that for any $\lambda$

(C-1) $k_{\lambda}^{-1}(\Sigma X_{k}(N,P))=S^{X_{k}}$, and $\pi_{P}%
^{\mathbf{\infty}}\circ k_{\lambda}|S^{X_{k}}=\pi_{P}^{\mathbf{\infty}}\circ
k_{0}|S^{X_{k}\text{ }}$and,

(C-2) $k_{\lambda}$ is smooth and transverse to $\Sigma X_{k}(N,P)$.

Recall the fiber map $d_{1}(s,X_{k})|\frak{n}(X_{k})+q^{X_{k}}(s)$ over
$S^{X_{k}}$ in (6.5). If we choose $\delta$ sufficiently small compared with
$\gamma$, then we can define the $\Omega$-regular map $g_{0}:\mathbf{E}%
(k)\rightarrow P$ by
\begin{equation}
g_{0}(\mathbf{v}_{c})=\exp_{P,\pi_{P}^{\mathbf{\infty}}\circ s(c)}\circ
(\pi_{P}^{\mathbf{\infty}}\circ s|S^{X_{k}})^{TP}\circ(d_{1}(s,X_{k}%
)|\frak{n}(X_{k})+q^{X_{k}}(s))_{c}\circ\exp_{N,c}^{-1}(\mathbf{v}_{c})
\end{equation}
{for each }$X_{k}=A_{k}$, $D_{k}$ and $E_{k}${. }It follows from Remark 3.2
that $g_{0}$ has each point $c\in S^{X_{k}}$ as the singularity of type
$X_{k}$\ and vice versa. Now we need to modify $g_{0}$ by using Lemma 5.3 so
that $g_{0}$ is compatible with $g_{k+1}$. Let $\eta:S^{X_{k}}\rightarrow
\mathbf{R}$ be a smooth function such that

(i) $0\leq\eta(c)\leq1$ for $c\in S^{X_{k}},$

(ii) $\eta(c)=0$ for $c\in S^{X_{k}}\cap{U(C_{k+1})}_{2-3r_{0}}$,

(iii) $\eta(c)=1$ for $c\in S^{X_{k}}\setminus{U(C_{k+1})}_{2-4r_{0}}$.

\noindent Then consider the map $G:{U(C_{k+1})}_{2-3r_{0}}\cup\mathbf{E}%
(k)\rightarrow P$ defined by
\[
\left\{
\begin{array}
[c]{lll}%
G(x)=g_{k+1}(x) & \text{if $x\in{U(C_{k+1})}_{2-3r_{0}}$}, & \\
G(\mathbf{v}_{c})=(1-\eta(c))g_{k+1}(\mathbf{v}_{c})+\eta(c)g_{0}%
(\mathbf{v}_{c}) & \text{if }\mathbf{v}_{c}\text{$\in\mathbf{E}(k)$}. &
\end{array}
\right.
\]
It follows from Lemma 5.3 that $G$ is an $\Omega$-regular map defined on
${U(C_{k+1})}_{2-3r_{0}}\cup\mathbf{E}(k)$, that $G|\mathbf{E}(k)$ has the
singularities of type $X_{k}$ exactly on $S^{X_{k}}$, and that for any \ $c\in
S^{X_{k}}$, the properties (1)-(7) of Lemma 5.3 hold for $G$. Furthermore, we
note that if $c\in\mathbf{E}(k)\cap{U(C_{k+1})}_{2-3r_{0}},$ then
$G(\mathbf{v}_{c})=g_{k+1}(\mathbf{v}_{c})$.

Set $\exp_{\Omega}=\exp_{\Omega(N,P)}$ for short. Let $h_{1}^{1}$, $h_{0}%
^{3}:(\mathbf{E}(k),S^{X_{k}})\rightarrow(\Omega(N,P),\Sigma X_{k}(N,P))$ be
the maps defined by
\begin{align*}
h_{1}^{1}(\mathbf{v}_{c})  &  =\exp_{\Omega,s(c)}\circ d_{c}s\circ(\exp
_{N,c})^{-1}(\mathbf{v}_{c}),\\
h_{0}^{3}(\mathbf{v}_{c})  &  =\exp_{\Omega,j^{\mathbf{\infty}}G(c)}\circ
d_{c}(j^{\mathbf{\infty}}G)\circ(\exp_{N,c})^{-1}(\mathbf{v}_{c}).
\end{align*}
By applying Lemma 5.4 to the section $s$ and $h_{1}^{1}$, we first obtain a
homotopy $h_{\lambda}^{1}\in\Gamma_{\Omega}^{tr}(\mathbf{E}(k),P)$ between
$h_{0}^{1}=s$ and $h_{1}^{1}$ on $\mathbf{E}(k)$ satisfying the properties
(1), (2) and (3) of Lemma 5.4. Similarly we obtain a homotopy $h_{\lambda}%
^{3}\in\Gamma_{\Omega}^{tr}(\mathbf{E}(k),P)$ between $h_{0}^{3}$ and
$h_{1}^{3}=j^{\mathbf{\infty}}G$ on $\mathbf{E}(k)$ satisfying the properties
(1), (2) and (3) of Lemma 5.4.

Next we construct a homotopy of bundle maps $\frak{n}(X_{k})\rightarrow
\nu(X_{k})$ covering a homotopy $S^{X_{k}}\rightarrow\Sigma X_{k}(N,P)$
between $ds|\frak{n}(X_{k})$ and $d(j^{\mathbf{\infty}}G)|\frak{n}(X_{k})$.
Recall the additive structure of $J^{\mathbf{\infty}}(N,P)$ in (1.2). Then we
have the homotopy $\kappa_{\lambda}:\mathbf{E}(k)\rightarrow J^{\mathbf{\infty
}}(N,P)$ defined by
\[
\kappa_{\lambda}(c)=(1-\lambda)s(c)+\lambda j^{\mathbf{\infty}}G(c)\quad
\text{covering }\pi_{P}^{\mathbf{\infty}}\circ s|S^{X_{k}}:S^{X_{k}%
}\rightarrow P,
\]
where $\pi_{P}^{\mathbf{\infty}}\circ s|S^{X_{k}}$ is the immersion as in
$(5.2.2)$.

We show that $\kappa_{\lambda}|S^{X_{k}}$ is actually a homotopy of $S^{X_{k}%
}$ into $\Sigma X_{k}(N,P)$. Recall the identification $(s|S^{X_{k}})^{\ast
}\mathbf{P}\cong(\pi_{P}^{\mathbf{\infty}}\circ s)^{\ast}(TP)$ and $s^{\ast
}\mathbf{D}\cong TN$.

It follows from the decomposition of $\frak{n}(X_{k})$ in (4.3) that%
\begin{equation}%
\begin{array}
[c]{l}%
(s|S^{I_{j}})^{\ast}(\mathbf{d}_{j+1}^{I}\circ ds|\frak{n}(s,I_{j}\subset
I_{j-1}))=(j^{\mathbf{\infty}}G|S^{I_{j}})^{\ast}(\mathbf{d}_{j+1}^{I}\circ
d(j^{\mathbf{\infty}}G)|\frak{n}(s,I_{j}\subset I_{j-1})),\\
(s|S^{\overline{D}_{j+1}})^{\ast}(d(\mathbf{r}_{j-1})\circ ds|\frak{n}%
(s,\overline{D}_{j+1}\subset\overline{D}_{j}))\\
\text{ \ \ \ \ \ \ \ \ \ \ \ \ \ \ \ \ \ \ \ \ \ \ \ \ \ \ }%
=(j^{\mathbf{\infty}}G|S^{\overline{D}_{j+1}})^{\ast}(d(\mathbf{r}_{j-1})\circ
d(j^{\mathbf{\infty}}G)|\frak{n}(s,\overline{D}_{j+1}\subset\overline{D}%
_{j})),\\
(s|S^{D_{k}})^{\ast}(\mathbf{r}_{k-1})=(j^{\mathbf{\infty}}G|S^{D_{k}})^{\ast
}(\mathbf{r}_{k-1}),\\
(s|S^{\overline{E}_{7}})^{\ast}(d(\widetilde{\mathbf{d}}_{4}^{\frak{J}%
}|\bigcirc^{4}\mathbf{K}_{3}^{\frak{J}})\circ ds|\frak{n}(s,\overline{E}%
_{7}\subset\overline{E}_{6}))\\
\text{ \ \ \ \ \ \ \ \ \ \ \ \ \ \ \ \ \ \ \ \ \ \ \ \ \ \ \ }%
=(j^{\mathbf{\infty}}G|S^{\overline{E}_{7}})^{\ast}(d(\widetilde{\mathbf{d}%
}_{4}^{\frak{J}}|\bigcirc^{4}\mathbf{K}_{3}^{\frak{J}})\circ
d(j^{\mathbf{\infty}}G)|\frak{n}(s,\overline{E}_{7}\subset\overline{E}%
_{6})),\\
(s|S^{E_{8}})^{\ast}(d(\widetilde{\mathbf{d}}_{4}^{\frak{J}}|\bigcirc
^{3}\mathbf{K}_{4}^{\frak{J}}\bigcirc(\mathbf{K}_{2}^{\frak{J}}/\mathbf{K}%
_{4}^{\frak{J}}))\circ ds|\frak{n}(s,E_{8}\subset\overline{E}_{7}))\\
\text{ \ \ \ \ \ \ \ \ \ \ \ \ \ \ \ \ \ \ \ \ \ \ \ \ \ \ \ }%
=(j^{\mathbf{\infty}}G|S^{E_{8}})^{\ast}(d(\widetilde{\mathbf{d}}%
_{4}^{\frak{J}}|\bigcirc^{3}\mathbf{K}_{4}^{\frak{J}}\bigcirc(\mathbf{K}%
_{2}^{\frak{J}}/\mathbf{K}_{4}^{\frak{J}}))\circ d(j^{\mathbf{\infty}%
}G)|\frak{n}(s,E_{8}\subset\overline{E}_{7}))\\
(s|S^{E_{8}})^{\ast}(\widetilde{\mathbf{d}}_{5}^{\frak{J}}|\bigcirc
^{5}\mathbf{K}_{4}^{\frak{J}})=(j^{\mathbf{\infty}}G|S^{E_{8}})^{\ast
}(\widetilde{\mathbf{d}}_{5}^{\frak{J}}|\bigcirc^{5}\mathbf{K}_{4}^{\frak{J}%
}),
\end{array}
\end{equation}
over respective $S^{X_{k}},$\ where $I=J$ or $\frak{J}$. These formulas are
the direct consequence of the construction of $d_{1}(s,X_{k})|\frak{n}%
(X_{k})+q^{X_{k}}(s)$ appearing in the definition of $G$ and the definitions
of the intrinsic derivatives and $\mathbf{r}_{k-1}$ in Sections 2 and 3.
Hence, it follows from (6.7) that for any $c\in S^{X_{k}}$ we have
$\frak{n}(\kappa_{\lambda},X_{k})_{c}=\frak{n}(X_{k})_{c}$ and $Q(\kappa
_{\lambda})_{c}=Q_{c}$. Hence, the equalities of the homomorphisms in (6.7)
also hold when $s$ is replaced by $\kappa_{\lambda}$ $(0\leq\lambda\leq1)$.
Therefore, $\kappa_{\lambda}|S^{X_{k}}$ gives a homotopy of $S^{X_{k}}$ into
$\Sigma X_{k}(N,P)$.

Let us consider the commutative diagram%
\[%
\begin{array}
[c]{llllll}%
\frak{n}(X_{k}) & \text{ \ \ \ }\overset{\underrightarrow{d(\kappa_{\lambda
})|\frak{n}(X_{k})\text{ }}}{} & \nu(X_{k}) & \text{ \ \ \ }\overset
{\underrightarrow{\mathbf{d}(X_{k})\text{ }}}{} & \mathbf{H(}X_{k}) & \\
\text{ \ }\parallel &  & \text{ \ \ \ }\uparrow\scriptstyle{(\kappa_{\lambda
}|S^{X_{k}})^{\nu(X_{k})}} &  & \text{ \ }\uparrow\scriptstyle{(\kappa
_{\lambda}|S^{X_{k}})^{\mathbf{H(}X_{k})}} & \\
\frak{n}(X_{k}) & \overset{\underrightarrow{(\kappa_{\lambda}|S^{X_{k}}%
)^{\ast}(d(\kappa_{\lambda})|\frak{n}(X_{k}))\text{ }}}{} & (\kappa_{\lambda
}|S^{X_{k}})^{\ast}(\nu(X_{k})) & \overset{\underrightarrow{(\kappa_{\lambda
}|S^{X_{k}})^{\ast}(\mathbf{d}(X_{k}))\text{ }}}{} & H(X_{k}). &
\end{array}
\]
Here, $\mathbf{H(}X_{k})$ (resp. $\mathbf{d}(X_{k})$) denotes the direct sum
of the target bundles of the homomorphisms (resp. the sum of the
homomorphisms) $\mathbf{d}_{j+1}^{I}|\nu(I_{j}\subset I_{j-1})$ for $I=J$ or
$\frak{J}$, $d(\mathbf{r}_{j-1})|\nu(\overline{D}_{j+1}\subset\overline{D}%
_{j})$, $d(\widetilde{\mathbf{d}}_{4}^{\frak{J}}|\bigcirc^{4}\mathbf{K}%
_{3}^{\frak{J}})|\nu(\overline{E}_{7}\subset\overline{E}_{6})$ and
$d(\widetilde{\mathbf{d}}_{4}^{\frak{J}}|\bigcirc^{3}\mathbf{K}_{4}^{\frak{J}%
}\bigcirc(\mathbf{K}_{2}^{\frak{J}}/\mathbf{K}_{4}^{\frak{J}}))|\nu
(E_{8}\subset\overline{E}_{7})$ such that $\nu(I_{j}\subset I_{j-1})$,
$\nu(\overline{D}_{j+1}\subset\overline{D}_{j})$, $\nu(\overline{E}_{7}%
\subset\overline{E}_{6})$ and $\nu(E_{8}\subset\overline{E}_{7})$, which
appear in the direct sum decompositions of $\nu(X_{k})$ in (4.2), and let
$H(X_{k})=(s|S^{X_{k}})^{\ast}(\mathbf{H}(X_{k}))$. Namely, $H(X_{k})$ is the
direct sum of the target bundles of the isomorphisms in (6.2), (6.3) and (6.4)
whose source normal bundles appear in the direct sum decomposition of
$\frak{n}(X_{k})$ in (4.3). Then it follows from (6.7) that%
\[
(\kappa_{\lambda}|S^{X_{k}})^{\ast}(\mathbf{d}(X_{k})\circ d(\kappa_{\lambda
})|\frak{n}(X_{k}))=(s|S^{X_{k}})^{\ast}(\mathbf{d}(X_{k})\circ ds|\frak{n}%
(X_{k}))
\]
for any $\lambda$. This implies that $\kappa_{\lambda}$ is transverse to
$\Sigma X_{k}(N,P)$ for any $\lambda$.

We define $h_{\lambda}^{2}:(\mathbf{E}(k),S^{X_{k}})\rightarrow(\Omega
(N,P),\Sigma X_{k}(N,P))$ by%
\[
h_{\lambda}^{2}(\mathbf{v}_{c})=\exp_{\Omega,\kappa_{\lambda}(c)}\circ
d_{c}(\kappa_{\lambda})\circ(\exp_{N,c})^{-1}(\mathbf{v}_{c}).
\]
Then we have that $h_{0}^{2}(\mathbf{v}_{c})=h_{1}^{1}(\mathbf{v}_{c}%
)=\exp_{\Omega,s(c)}\circ d_{c}s\circ(\exp_{N,c})^{-1}(\mathbf{v}_{c})$ and
$h_{0}^{3}(\mathbf{v}_{c})=h_{1}^{2}(\mathbf{v}_{c})=\exp_{\Omega
,j^{\mathbf{\infty}}G(c)}\circ d_{c}(j^{\mathbf{\infty}}G)\circ(\exp
_{N,c})^{-1}(\mathbf{v}_{c})$ on $\mathbf{E}(k)$.

Let $\overline{h}_{\lambda}\in\Gamma_{\Omega}^{tr}(\mathbf{E}(k)\cup
{U(C_{k+1})}_{2-3r_{0}},P)$ be the homotopy which is obtained by pasting
$h_{\lambda}^{1}$, $h_{\lambda}^{2}$ and $h_{\lambda}^{3}$. Since $h_{\lambda
}^{1}$ and $h_{\lambda}^{3}$ do not keep $s|(\mathbf{E}(k)\cap({U(C_{k+1}%
)}_{2-3r_{0}}\setminus$Int${U(C_{k+1})}_{2}))$ in general, we need to modify
$h_{\lambda}$ as follows.

Since $h_{0}^{1}(\mathbf{v}_{c})=h_{1}^{3}(\mathbf{v}_{c})=s(\mathbf{v}_{c})$
for $\mathbf{v}_{c}\in\mathbf{E}(k)\cap{U(}C_{k+1}{)}_{2-3r_{0}}$, we may
assume in the construction of $h_{\lambda}^{1}$, $h_{\lambda}^{2}$ and
$h_{\lambda}^{3}$ that if $\mathbf{v}_{c}\in\mathbf{E}(k)\cap{U(}C_{k+1}%
{)}_{2-3r_{0}}$, then $h_{\lambda}^{2}(\mathbf{v}_{c})=h_{0}^{2}%
(\mathbf{v}_{c})=h_{1}^{2}(\mathbf{v}_{c})$ and $h_{\lambda}^{1}%
(\mathbf{v}_{c})=h_{1-\lambda}^{3}(\mathbf{v}_{c})$ for any $\lambda$. By
using these properties of $h_{\lambda}^{1}$, $h_{\lambda}^{2}$ and
$h_{\lambda}^{3}$, we can modify $\overline{h}_{\lambda}$ to a homotopy
$h_{\lambda}\in\Gamma_{\Omega}^{tr}(\mathbf{E}(k),P)$\ satisfying the property
(C) such that

(1) $h_{\lambda}(\mathbf{v}_{c})=h_{0}(\mathbf{v}_{c})=s(\mathbf{v}_{c})$ for
any $\lambda$\ and any $\mathbf{v}_{c}\in\mathbf{E}(k)\cap{U(}C_{k+1}%
{)}_{2-2r_{0}}$,

(2) $h_{0}(\mathbf{v}_{c})=s(\mathbf{v}_{c})$ for any $\mathbf{v}_{c}%
\in\mathbf{E}(k),$

(3) $h_{1}(\mathbf{v}_{c})=j^{\mathbf{\infty}}G(\mathbf{v}_{c})$ for any
$\mathbf{v}_{c}\in\mathbf{E}(k)$.

By (1), we can extend $h_{\lambda}$ to the homotopy $H_{\lambda}^{\prime}%
\in\Gamma_{\Omega}^{tr}(\mathbf{E}(k)\cup{U(}C_{k+1}{)}_{2-2r_{0}}%
,P)$\ defined by $H_{\lambda}^{\prime}|\mathbf{E}(k)=h_{\lambda}$ and
$H_{\lambda}^{\prime}|{U(}C_{k+1}{)}_{2-2r_{0}}=s|{U(}C_{k+1}{)}_{2-2r_{0}}$.

By applying the homotopy extension property to $s$ and $H_{\lambda}^{\prime}$
for each $X_{k}=A_{k}$, $D_{k}$ and $E_{k}$, we obtain an extended homotopy
\[
H_{\lambda}:(N,S^{X_{k}})\rightarrow(\Omega(N,P),\Sigma X_{k}(N,P))
\]
relative to ${U(C_{k+1})}_{2-r_{0}}$\ with $H_{0}=s$. Furthermore, we replace
$\delta$ and $\mathbf{E}(k)$ by smaller ones. Then $H_{\lambda}$ is a required
homotopy in $\Gamma_{\Omega}^{tr}(N,P)$ in the assertion (\textbf{A}).

We need further argument for (4.1.3) using the Thom's first Isotopy Lemma (see
[Math3, (8.1)]). We suppose that $\Omega(N,P)$ does not contain $\Sigma
D_{j}(N,P)$ or $\Sigma E_{k}(N,P)$. Let $\widehat{\delta}$ and $\delta
^{\prime}$ be smooth functions with $0<\widehat{\delta}<\delta^{\prime}%
<\delta$. Let $U(A_{k})=$Int$(\mathbf{E}(k))\cup{U(}C_{k+1}{)}_{2-2r_{0}}$ and
$U(k)_{\varsigma}=\exp_{N}(D_{\varsigma\circ s}(\frak{n}(A_{k})))\cup
{U(}C_{k+1}{)}_{2-(3/2)r_{0}}$ for $\zeta=\widehat{\delta},$ $\delta^{\prime}%
$. Let $p^{I}:U(A_{k})\times I\rightarrow I$ be the projection onto the second
factor. Let $S^{A_{j}}(\{H_{\lambda}^{\prime}\})$ be the subset of
$U(A_{k})\times I$ which consists of all points $(x,\lambda)$ such that
$H_{\lambda}^{\prime}(x)\in\Sigma A_{j}(N,P)$ for each $j\leq k$. Since
$H_{\lambda}^{\prime}$ is transverse to $\Sigma A_{k}(N,P)$, we may assume
that $H_{\lambda}^{\prime}$ is transverse to all $\Sigma A_{j}(N,P)$ on
$U(A_{k})$ for any $\lambda$. Since $p^{I}$ is a submersion onto $I$, we can
construct a smooth vector field $V(x,\lambda)$ on $U(A_{k})\times I$ such that

(i) if $(x,\lambda)\in S^{A_{j}}(H_{\lambda}^{\prime})$, then $V(x,\lambda)$
is tangent to $S^{A_{j}}(H_{\lambda}^{\prime})$ at $(x,\lambda),$

(ii) $dp^{I}(V(x,\lambda))=\partial/\partial\lambda,$

(iii) if $(x,\lambda)\in(S^{A_{k}}\cup{U(}C_{k+1}{)}_{2-(3/2)r_{0}})\times I$,
then $V(x,\lambda)=\partial/\partial\lambda$.

\noindent Then by considering the integral curve of $V(x,\lambda)$, we can
choose $\widehat{\delta}$, $\delta^{\prime}$ and a smooth isotopy of
embeddings $\tau_{\lambda}:U(k)_{\delta^{\prime}}\rightarrow U(A_{k})$ such that

(i) $\tau_{0}(x)=x,$

(ii) each path $(\tau_{\lambda}(x),\lambda)$, $x\in U(k)_{\delta^{\prime}}$
and $\lambda\in I$ lies in $U(A_{k})\times I$ and $(d/d\lambda)(\tau_{\lambda
}(x),\lambda)=V(x,\lambda),$

(iii) if $x\in S^{A_{k}}\cup{U(}C_{k+1}{)}_{2-(3/2)r_{0}}$, then
$\tau_{\lambda}(x)=x$ for any $\lambda.$

\noindent Define $H_{\lambda}^{\prime}\bullet j^{\mathbf{\infty}}\tau
_{\lambda}\in\Gamma_{\Omega}(U(k)_{\delta^{\prime}},P)$\ by $H_{\lambda
}^{\prime}\bullet j^{\mathbf{\infty}}\tau_{\lambda}(x)=H_{\lambda}^{\prime
}(x)\circ j^{\mathbf{\infty}}\tau_{\lambda}(x)$. Then we have

(1) $S^{A_{j}}(H_{\lambda}^{\prime}\bullet j^{\mathbf{\infty}}\tau_{\lambda
})\cap U(k)_{\delta^{\prime}}=S^{A_{j}}(H_{0}^{\prime})\cap U(k)_{\delta
^{\prime}}$ for $j\leq k$ and any $\lambda,$

(2) $H_{\lambda}^{\prime}\bullet j^{\mathbf{\infty}}\tau_{\lambda}$ is
transverse to all $\Sigma A_{j}(N,P),$

(3) $H_{1}^{\prime}\bullet j^{\mathbf{\infty}}\tau_{1}=j^{\mathbf{\infty}%
}(G\circ\tau_{1})$, where $G\circ\tau_{1}$ is an $\Omega$-regular map.

\noindent Set $U^{k}=U(k)_{\delta^{\prime}}$. Suppose that there exist a
neighborhood $U^{j}$ ($j\leq k$) with $U(k)_{\widehat{\delta}}\subset
$Int$U^{j}$, $U^{i}\subset$Int$U^{i+1}$ ($j\leq i<k$), and a smooth homotopy%
\[
H_{\lambda}^{j}:(U^{j}\cup S^{\overline{A}_{j}}(s),S^{\overline{A}_{j}%
}(s))\rightarrow(\Omega(N,P),\Sigma\overline{A}_{j}(N,P))
\]
relative to $U(k)_{\widehat{\delta}}$\ such that $S^{A_{\ell}}(H_{\lambda}%
^{j})=S^{A_{\ell}}(s)$ $(j\leq\ell\leq k)$, $H_{0}^{j}|U^{j}=s|U^{j}$ and
$H_{\lambda}^{j}|U^{j}=H_{\lambda}^{\prime}\bullet j^{\mathbf{\infty}}%
\tau_{\lambda}|U^{j}$. Then there exists a small tubular neighborhood
$T(j-1)\ $of $S^{\overline{A}_{j}}(s)$ in $S^{\overline{A}_{j-1}}(s)$\ with
projection $p^{T(j-1)}:T(j-1)\rightarrow S^{\overline{A}_{j}}(s)$ and a
neighborhood $U^{j-1}$\ with $U(k)_{\widehat{\delta}}\subset$Int$U^{j-1}$ and
$U^{j-1}\subset$Int$U^{j}$ such that $U^{j-1}\cap T(j-1)=(p^{T(j-1)}%
)^{-1}(S^{\overline{A}_{j}}(s)\cap U^{j-1})$. Since $s$ is transverse to
$\Sigma A_{j}(N,P)$, we have a smooth homotopy%
\[
H_{\lambda}^{T(j-1)}:(U^{j-1}\cup T(j-1),S^{\overline{A}_{j}}(s))\rightarrow
(\Omega(N,P),\Sigma\overline{A}_{j}(N,P))
\]
such that $S^{A_{\ell}}(H_{\lambda}^{T(j-1)})=S^{A_{\ell}}(s)$ $(j-1\leq
\ell\leq k)$, $H_{\lambda}^{T(j-1)}|U^{j-1}\cup S^{\overline{A}_{j}%
}(s)=H_{\lambda}^{j}|U^{j-1}\cup S^{\overline{A}_{j}}(s)$ and $H_{0}%
^{T(j-1)}|U^{j-1}\cup T(j-1)=s|U^{j-1}\cup T(j-1)$ by using (1), (2) and (3)
and applying the homotopy extension property to$\ s|U^{j-1}\cup T(j-1)$ and
$H_{\lambda}^{j}|U^{j-1}\cup S^{\overline{A}_{j}}(s)$. Next we can easily
extend $H_{\lambda}^{T(j-1)}$ to a smooth homotopy%
\[
H_{\lambda}^{j-1}:(U^{j-1}\cup S^{\overline{A}_{j-1}}(s),S^{\overline{A}%
_{j-1}}(s),S^{A_{j-1}}(s))\rightarrow(\Omega(N,P),\Sigma\overline{A}%
_{j-1}(N,P),\Sigma A_{j-1}(N,P))
\]
by applying the homotopy extension property to$\ s|U^{j-1}\cup S^{A_{j-1}}(s)$
and $H_{\lambda}^{T(j-1)}$ so that $S^{A_{j-1}}(H_{\lambda}^{j-1})=S^{A_{j-1}%
}(s)$.

By the downward induction on $j$, we obtain an extended smooth homotopy
\[
H_{\lambda}:(N,S^{A_{k}}(s))\rightarrow(\Omega(N,P),\Sigma A_{k}(N,P))
\]
relative to ${U(C_{k+1})}_{2-r_{0}}$\ such that $H_{0}=s$, $S^{A_{j}%
}(H_{\lambda})=S^{A_{j}}(s)$ for $1\leq j\leq k$. This completes the proof.
\end{proof}

The author does not know whether the Whitney condition (b) for the Thom's
first Isotopy Lemma holds or not among $\Sigma A_{i}(N,P)$, $\Sigma
D_{j}(N,P)$ and $\Sigma E_{k}(N,P).$

\section{Proof of Theorem 0.2}

In this section we prove Theorem 0.2 by applying Theorem 0.1 in the
elimination of higher $D_{k}$ and $E_{k}$ singularities.

Let $\Omega^{\overline{D}}(N,P)$ (resp. $\Omega^{\overline{E}}(N,P)$) denote
the open subbundle of $J^{\infty}(N,P)$ which consists of all regular jets and
$\Sigma A_{i}(N,P)$ ($i\geq1$) and $\Sigma D_{j}(N,P)$ ($j\geq4$)\ (resp.
$\Sigma A_{i}(N,P)$ ($i\geq1$), $\Sigma D_{j}(N,P)$ ($j\geq4$)\ and $\Sigma
E_{k}(N,P)$ ($k\geq6$)).\ Let $\Omega^{D_{k}}(N,P)$ (resp. $\Omega^{D_{5}%
E_{6}}(N,P)$) denote the open subbundle of $J^{\infty}(N,P)$ which consists of
all regular jets and $\Sigma A_{i}(N,P)$ ($i\geq1$) and $\Sigma D_{j}(N,P)$
($k\geq j\geq4$)\ (resp. $\Sigma A_{i}(N,P)$ ($i\geq1$), $\Sigma D_{j}(N,P)$
($j=4$, $5$) and $\Sigma E_{6}(N,P)$).

\begin{proof}
[Proof of Theorem 0.2]By the assumption, $j^{\infty}f$ is a section
$N\rightarrow\Omega^{\overline{D}}(N,P)$ (resp. $N\rightarrow\Omega
^{\overline{E}}(N,P)$). By Proposition 7.1 below, we have the section
$s^{D}:N\rightarrow\Omega^{D_{5}}(N,P)$ (resp. $s^{E}:N\rightarrow
\Omega^{D_{5}E_{6}}(N,P)$) such that $\pi_{P}^{\infty}\circ j^{\infty}%
f=\pi_{P}^{\infty}\circ s^{D}$ (resp. $\pi_{P}^{\infty}\circ j^{\infty}%
f=\pi_{P}^{\infty}\circ s^{E}$). By Theorem 0.1 we obtain a required smooth
map $g$ such that $j^{\infty}g$ and $s^{D}$ (resp. $s^{E}$) are homotopic.
This proves the assertion.
\end{proof}

\begin{proposition}
Let $n>p\geq2$ and $n-p$ be even. Then we have the following.

$(1)$ Given a section $s:N\rightarrow\Omega^{\overline{D}}(N,P)$, there exists
a section $s^{D}:N\rightarrow\Omega^{D_{5}}(N,P)$ such that $\pi_{P}^{\infty
}\circ s=\pi_{P}^{\infty}\circ s^{D}$.

$(2)$ Given a section $s:N\rightarrow\Omega^{\overline{E}}(N,P)$, there exists
a section $s^{E}:N\rightarrow\Omega^{D_{5}E_{6}}(N,P)$ such that $\pi
_{P}^{\infty}\circ s=\pi_{P}^{\infty}\circ s^{E}$.
\end{proposition}

We need the following lemma for the proof of Proposition 7.1.

\begin{lemma}
Let $n>p\geq2$ and $n-p$ be even. Then we have that, for $\frak{J}$,
$\mathbf{Q|}_{\Sigma^{\frak{J}_{2}}(N,P)}$, $\mathbf{L|}_{\Sigma\overline
{D}_{5}(N,P)}$ and $(\mathbf{K}_{2}^{\frak{J}}/\mathbf{K}_{3}^{\frak{J}%
})|_{\Sigma^{\frak{J}_{3}}(N,P)}$ are trivial line bundles equipped with the
canonical orientations respectively.
\end{lemma}

\begin{proof}
By Section 2 (5), $\mathbf{d}_{2}|\mathbf{K}:\mathbf{K}\rightarrow
\mathrm{Hom}(\mathbf{K},\mathbf{Q})$ induces the isomorphism%
\[
\mathbf{K/K}_{2}^{\frak{J}}\rightarrow\text{Hom}(\mathbf{K/K}_{2}^{\frak{J}%
},\mathbf{Q})\text{ \ \ \ over }\Sigma^{\frak{J}_{2}}(N,P),
\]
which yields the nonsingular quadratic form $\mathbf{q:K/K}_{2}^{\frak{J}%
}\bigcirc\mathbf{K/K}_{2}^{\frak{J}}\rightarrow\mathbf{Q}$\ over
$\Sigma^{\frak{J}_{2}}(N,P)$. Since dim$\mathbf{K/K}_{2}^{\frak{J}}=n-p-1$ is
odd, we choose the unique orientation of $\mathbf{Q}_{z}$, expressed by the
unit vector $\mathbf{e}_{P}$, so that the index (the number of the negative
eigen values) of $\mathbf{q}_{z}$, $z\in\Sigma^{\frak{J}_{2}}(N,P)$ is less
than $(n-p-1)/2$.

By ($\frak{J}$-1) and (D-ii) we have the following direct sum decompositions
over $\Sigma\overline{D}_{5}(N,P):$%
\[%
\begin{array}
[c]{l}%
\nu(\frak{J}_{2}\subset\frak{J}_{1})\mathbf{|}_{\Sigma\overline{D}_{5}%
}=(\mathbf{K}_{2}^{\frak{J}}/\mathbf{L}\oplus\mathbf{L}\oplus\mathbf{T}%
^{1})\mathbf{|}_{\Sigma\overline{D}_{5}},\\
\text{Hom(}\bigcirc^{2}\mathbf{K}_{2}^{\frak{J}},\mathbf{Q)|}_{\Sigma
\overline{D}_{5}}\cong\text{Hom(}\mathbf{K}_{2}^{\frak{J}}/\mathbf{L\bigcirc
L\oplus}\bigcirc^{2}(\mathbf{K}_{2}^{\frak{J}}/\mathbf{L)\oplus}\bigcirc
^{2}\mathbf{L},\mathbf{Q)|}_{\Sigma\overline{D}_{5}},
\end{array}
\]
where $\mathbf{T}^{1}$ is the orthogonal complement of $\mathbf{K}%
_{2}^{\frak{J}}$ in $\nu(\frak{J}_{2}\subset\frak{J}_{1})$ over $\Sigma
\overline{D}_{5}(N,P)$.

By (2.3), (D-ii) and (D-iii), ($\mathbf{d}_{3}^{\frak{J}}|\nu(\frak{J}%
_{2}\subset\frak{J}_{1})\mathbf{)|}_{\Sigma\overline{D}_{5}}$ induces the
isomorphisms
\[%
\begin{array}
[c]{l}%
(\mathbf{K}_{2}^{\frak{J}}/\mathbf{L)|}_{\Sigma\overline{D}_{5}}%
\mathbf{\rightarrow}\text{Hom(}\mathbf{K}_{2}^{\frak{J}}/\mathbf{L\bigcirc
L},\mathbf{Q)|}_{\Sigma\overline{D}_{5}},\\
\mathbf{L|}_{\Sigma\overline{D}_{5}}\mathbf{\rightarrow}\text{Hom(}%
\bigcirc^{2}(\mathbf{K}_{2}^{\frak{J}}/\mathbf{L)},\mathbf{Q)|}_{\Sigma
\overline{D}_{5}},\\
\mathbf{T}^{1}|_{\Sigma\overline{D}_{5}}\mathbf{\rightarrow}\text{Hom(}%
\bigcirc^{2}\mathbf{L},\mathbf{Q)|}_{\Sigma\overline{D}_{5}}.
\end{array}
\]
Since $\bigcirc^{2}(\mathbf{K}_{2}^{\frak{J}}/\mathbf{L)}$ has the canonical
orientation, $\mathbf{L}$ has the canonical orientation, expressed by the unit
vector $\mathbf{e}(\mathbf{L})$ over $\Sigma\overline{D}_{5}(N,P)$, by the
second isomorphism.

Let us provide $\mathbf{K}_{2}^{\frak{J}}/\mathbf{K}_{3}^{\frak{J}}$ with the
orientation. By ($\frak{J}$-1), $\mathbf{d}_{3}^{\frak{J}}|\mathbf{K}%
_{2}^{\frak{J}}$ induces the isomorphism $\mathbf{K}_{2}^{\frak{J}}%
/\mathbf{K}_{3}^{\frak{J}}\rightarrow\mathrm{Hom}(\bigcirc^{2}(\mathbf{K}%
_{2}^{\frak{J}}/\mathbf{K}_{3}^{\frak{J}}),\mathbf{Q})$ over $\Sigma
^{\frak{J}_{3}}(N,P)$. Then it comes from the orientation of $\mathrm{Hom}%
(\bigcirc^{2}(\mathbf{K}_{2}^{\frak{J}}/\mathbf{K}_{3}^{\frak{J}}%
),\mathbf{Q})$.
\end{proof}

\begin{lemma}
Let $n>p\geq2$. Let $z$ be a point of $\Sigma E_{6}(N,P)$, $\Sigma E_{7}(N,P)$
or $\Sigma E_{8}(N,P)$. Let $\{z_{m}\}$ be a sequence of $\Sigma\overline
{D}_{5}(N,P)$ which converges to $z$. Then we have

$(1)$ $\{\mathbf{K}_{2,z_{m}}^{\frak{J}}\}$ converges to $\mathbf{K}%
_{2,z}^{\frak{J}}$,

$(2)$ $\{\mathbf{L}_{z_{m}}\}$ converges to $\mathbf{K}_{3,z}^{\frak{J}}$.
\end{lemma}

\begin{proof}
(1) Since $\mathrm{Ker}(\mathbf{d}_{2}^{\frak{J}}|\mathbf{K)=K}_{2}^{\frak{J}%
}$ over $\Sigma^{\frak{J}_{2}}(N,P)$ and since $\Sigma\overline{D}%
_{5}(N,P)\subset\Sigma^{\frak{J}_{2}}(N,P)$ and $\Sigma E_{k}(N,P)\subset
\Sigma^{\frak{J}_{2}}(N,P)$, the assertion (1) follows from the continuity of
$\mathbf{d}_{2}^{\frak{J}}|\mathbf{K}$ on $\Sigma^{\frak{J}_{2}}(N,P)$.

(2) By (1), we take a subsequence of $\{z_{m}\}$ for which $\{\mathbf{L}%
_{z_{m}}\}$ converges to a $1$-dimensional subspace of $\mathbf{K}%
_{2,z}^{\frak{J}}$, say $\mathbf{V}_{z}$. By (2.4) it is enough for the proof
of the assertion (2) to show that $\mathbf{d}_{3,z}^{\frak{J}}|\mathbf{V}_{z}$
vanishes. Suppose that $\mathbf{V}_{z}\mathbf{\neq K}_{3,z}^{\frak{J}}$. Since
$\Sigma E_{k}(N,P)\subset\Sigma^{\frak{J}_{3}}(N,P)$, it follows that
$\mathrm{Ker}(\mathbf{d}_{3,z}^{\frak{J}}|\mathbf{K}_{2,z}^{\frak{J}%
})=\mathbf{K}_{3,z}^{\frak{J}}$. Hence, $\widetilde{\mathbf{d}}_{3,z}%
^{\frak{J}}$ induces the isomorphism $\bigcirc^{3}(\mathbf{K}_{2,z}^{\frak{J}%
}/\mathbf{K}_{3,z}^{\frak{J}})\rightarrow\mathbf{Q}_{z}$, which yields the
isomorphism $\bigcirc^{3}\mathbf{V}_{z}\mathbf{\rightarrow Q}_{z}$ (on the
contrary, we have that $\mathbf{d}_{3,z}^{\frak{J}}|\mathbf{K}_{3,z}%
^{\frak{J}}$ and $\widetilde{\mathbf{d}}_{3,z}^{\frak{J}}|\bigcirc
^{3}\mathbf{K}_{3,z}^{\frak{J}}$ vanish). Since $\lim_{m\rightarrow\infty
}\mathbf{L}_{z_{m}}=\mathbf{V}_{z}$, there is a number $m_{0}$ such that if
$m>m_{0}$, then $\widetilde{\mathbf{d}}_{3,z}^{\frak{J}}|\bigcirc
^{3}\mathbf{L}_{z_{m}}$ does not vanish. This contradicts to the definition of
$\Sigma\overline{D}_{5}(N,P)$ in (D-iii). Hence, we obtain the assertion (2).
\end{proof}

\begin{remark}
Under the same assumption of Lemma 7.2

$(1)$ $\widetilde{\mathbf{L}}$ has the canonical orientation if $U(\Sigma
\overline{D}_{5})$ is chosen as a tubular neighborhood,

$(2)$ the normal bundles for the respective inclusions $U(\Sigma\overline
{D}_{5})\supset\Sigma\overline{D}_{5}(N,P)\supset\cdots\supset\Sigma
\overline{D}_{j}(N,P)$ are all trivial by $(3.2)$.
\end{remark}

\begin{proof}
[Proof of Proposition 7.1]We give the proof only for (2), since the proof for
(1) is parallel by setting $S^{\overline{E}_{6}}(s)=\varnothing$.

In the proof we identify $J^{k}(N,P)$ with $J^{k}(TN,TP)$ by (1.2). By Remark
3.4 , there exists the open subbundle $\Omega^{D_{4}}(N,P)^{\prime}$ of
$J^{3}(N,P)$ (resp. $\Omega^{D_{5}E_{6}}(N,P)^{\prime}$ of $J^{4}(N,P)$) such
that $(\pi_{3}^{\infty})^{-1}(\Omega^{D_{4}}(N,P)^{\prime})=\Omega^{D_{4}%
}(N,P)$ (resp. $(\pi_{4}^{\infty})^{-1}(\Omega^{D_{5}E_{6}}(N,P)^{\prime
})=\Omega^{D_{5}E_{6}}(N,P)$). We may assume $s\in\Gamma_{\Omega^{\overline
{E}}}^{tr}(N,P)$. It follows that $(\pi_{3}^{\infty}\circ s)(N\setminus
(S^{\overline{D}_{5}}(s)\cup S^{\overline{E}_{6}}(s)))\subset\Omega^{D_{4}%
}(N\setminus(S^{\overline{D}_{5}}(s)\cup S^{\overline{E}_{6}}(s)),P)^{\prime}%
$. Since $s(S^{\overline{D}_{5}}(s))\subset\Sigma\overline{D}_{5}(N,P)$ and
$s(S^{\overline{E}_{6}}(s))\subset\Sigma^{\frak{J}_{3}}(N,P)$, we consider
$L=(s|S^{\overline{D}_{5}}(s))^{\ast}\mathbf{L}$ and $K_{3}=(s|S^{\overline
{E}_{6}}(s))^{\ast}\mathbf{K}_{3}^{\frak{J}}$.

We now construct a new section%
\[
\widetilde{u}:N\rightarrow\Omega^{D_{5}E_{6}}(N,P)^{\prime}%
\]
as follows.

We have that $r_{3}(s):\bigcirc^{3}L\rightarrow Q$ over $S^{D_{4}}(s)$ is an
isomorphism and is the null-homomorphism over $S^{\overline{D}_{5}}(s)$.
Furthermore, $r_{4}(s):\bigcirc^{4}L\rightarrow Q$ over $S^{D_{5}}(s)$ is an
isomorphism and is the null-homomorphism over $S^{\overline{D}_{6}}(s)$. Let
$e_{P}$, $\mathbf{e}(L)$ and $\mathbf{e}(K_{3}\bigcirc K_{3})$ be the unit
vectors induced from $\mathbf{e}_{P}$, $\mathbf{e(L)}$ and $\mathbf{e}%
(\mathbf{K}_{3}^{\frak{J}}\bigcirc\mathbf{K}_{3}^{\frak{J}})$, which
represents the canonical orientation of $\mathbf{K}_{3}^{\frak{J}}%
\bigcirc\mathbf{K}_{3}^{\frak{J}}$,\ by $s$ respectively. Then by using Lemma
7.3 we define the smooth isomorphisms $\phi^{D}:\bigcirc^{4}L\rightarrow Q$
over $S^{\overline{D}_{5}}(s)$ and $\phi^{E}:\bigcirc^{4}K_{3}\rightarrow Q$
over $S^{\overline{E}_{6}}(s)$ by $\phi^{D}(\bigcirc^{4}\mathbf{e}(L))=e_{P}$
and $\phi^{E}(\bigcirc^{2}\mathbf{e}(K_{3}\bigcirc K_{3}))=e_{P}$
respectively. Then we can find a section $u_{\phi}:S^{\overline{D}_{5}}(s)\cup
S^{\overline{E}_{6}}(s)\rightarrow\mathrm{Hom}(S^{4}((\pi_{N}^{\infty}\circ
s)^{\ast}(TN)),(\pi_{P}^{\infty}\circ s)^{\ast}(TP))$ such that $u_{\phi
}(x)|\bigcirc^{4}L_{x}=\phi^{D}|_{x}$ for $x\in S^{\overline{D}_{5}}(s)$ and
$u_{\phi}(x)|\bigcirc^{4}K_{3,x}=\phi^{E}|_{x}$ for $x\in S^{\overline{E}_{6}%
}(s)$. Since $S^{\overline{D}_{5}}(s)\cup S^{\overline{E}_{6}}(s)$ is a closed
subset and since Hom$(S^{4}((\pi_{N}^{\infty}\circ s)^{\ast}(TN)),(\pi
_{P}^{\infty}\circ s)^{\ast}(TP))$ is a vector bundle, we can extend $u_{\phi
}$ arbitrarily to the section $\widetilde{u_{\phi}}:N\rightarrow
\mathrm{Hom}(S^{4}((\pi_{N}^{\infty}\circ s)^{\ast}(TN)),(\pi_{P}^{\infty
}\circ s)^{\ast}(TP))$. Then we define $\widetilde{u}$ by $\widetilde{u}%
=\pi_{3}^{\infty}\circ s\oplus\widetilde{u_{\phi}}$ as the section of
$J^{4}(N,P)=$ $J^{4}(TN,TP)$. We regard $\widetilde{u}$ as the section of
$J^{\infty}(N,P)$ over $N$. We prove that $\widetilde{u}\in\Omega^{D_{5}E_{6}%
}(N,P)$. By the construction we have that $r_{4}(\widetilde{u})_{x}=u_{\phi
}(x)|\bigcirc^{4}L_{x}=\phi^{D}|_{x}$ for $x\in S^{\overline{D}_{5}}(s)$ and
$\widetilde{d}_{4,x}^{\frak{J}}|\bigcirc^{4}K_{3,x}=u_{\phi}(x)|\bigcirc
^{4}K_{3,x}=\phi^{E}|_{x}$ for $x\in S^{\overline{E}_{6}}(s)$. For any point
$x\in S^{\overline{D}_{5}}(s)\cup S^{\overline{E}_{6}}(s)$, let $U_{x}$ be a
convex neighborhood of $x$ and let $\ell$ and $k$\ be the coordinates of
$\exp_{N,x}(L_{x})$ and $\exp_{N,x}(K_{3,x})$ respectively. As in the proof of
Lemma 3.1, it follows from the definition of $\mathbf{D}$ that%
\[%
\begin{array}
[c]{ll}%
(\bigcirc^{4}\delta_{\ell})y_{p}|_{\widetilde{u}(x)}=\partial^{4}%
y_{p}/\partial\ell^{4}(x)\neq0 & \text{for }x\in S^{\overline{D}_{5}}(s),\\
(\bigcirc^{4}\delta_{k})y_{p}|_{\widetilde{u}(x)}=\partial^{4}y_{p}/\partial
k^{4}(x)\neq0 & \text{for }x\in S^{\overline{E}_{6}}(s).
\end{array}
\]
Hence, we have that $\widetilde{u}(S^{\overline{D}_{5}}(s))\subset\Sigma
D_{5}(N,P)$ and $\widetilde{u}(S^{\overline{E}_{6}}(s))\subset\Sigma
E_{6}(N,P)$. It is clear that $\widetilde{u}(N\setminus(S^{\overline{D}_{5}%
}(s)\cup S^{\overline{E}_{6}}(s)))\subset\Omega^{D_{4}}(N,P)$. This completes
the proof.
\end{proof}

Let $X=D$ or $E$ and $f:N\rightarrow P$ be an $\Omega^{\overline{X}}$-regular
map. In this paper the Thom polynomial $P(X_{k+1},f)$ of $S^{\overline
{X}_{k+1}}(j^{\infty}f)$ refers to the Poincar\'{e} dual in $N$ of the
fundamental class of $S^{\overline{X}_{k+1}}(j^{\infty}f)$. Let $W=1+W_{1}%
+\cdots+W_{j}+\cdots$ be the total Stiefel-Whitney class and $\overline{W}$ be
the formal inverse of $W$. Let $[a]$ denote the greatest integer not greater
than $a$.

If $n-p$ is odd, then Theorem 0.2 does not hold in general. We have the
following theorem ([An4, Remark 4.3 and Section 8]).

\begin{theorem}
Let $n=p+1$. Let $W_{j}=W_{j}(TN-f^{\ast}(TP))$. Then

$(1)$ if $f:N\rightarrow P$ is an $\Omega^{\overline{E}}$-regular map, then
the Thom polynomial $P(E_{k+1},f)$ of $S^{\overline{E}_{k+1}}(j^{\infty}f)$ is
equal to $W_{k}W_{2}+W_{k-1}(W_{3}+W_{1}W_{2}),$

$(2)$ if $f:N\rightarrow P$ is an $\Omega^{\overline{D}}$-regular map, then
the Thom polynomial $P(D_{k+1},f)$ of $S^{\overline{D}_{k+1}}(j^{\infty}f)$ is
equal to the part of degree $k+2$ of the polynomial%
\[
W(\overline{W}_{1}+\overline{W}_{2})\left\{  \sum_{j=0}^{[\frac{k}{2}%
-2]}\binom{[\frac{k}{2}-2]}{j}\overline{W}_{j}\right\}  +W\left\{  \sum
_{j=0}^{[\frac{k}{2}-1]}\binom{[\frac{k}{2}-1]}{j}\overline{W}_{j+1}\right\}
.
\]
\end{theorem}

Recall that codim$\Sigma^{2,2,2}(N,P)=9$ ([B, Theorem 6.1]) and codim$\Sigma
^{3}(N,P)=6$. By [H1, Theorem 7.1], there is an immersion $\mathbf{RP}%
^{6}\rightarrow\mathbf{R}^{7}$. Hence, any map $f:\mathbf{RP}^{6}%
\times\mathbf{RP}^{2}\rightarrow\mathbf{R}^{7}$ is homotopic to an
$\Omega^{\overline{E}}$-regular map. We have that $W(\mathbf{RP}^{6}%
\times\mathbf{RP}^{2})=(1+a)^{7}\otimes(1+b)^{3}$, where $a$ and $b$ are the
generators of $H^{1}(\mathbf{RP}^{6};\mathbf{Z}/(2))$ and $H^{1}%
(\mathbf{RP}^{2};\mathbf{Z}/(2))$ respectively. By Theorem 7.5 (1), we have
that $P(E_{7},f)=a^{6}\otimes b^{2}\neq0$ in this case. This implies that we
cannot eliminate the singularity of type $E_{7}$ from $f:\mathbf{RP}^{6}%
\times\mathbf{RP}^{2}\rightarrow\mathbf{R}^{7}$.

For $\Omega^{\overline{D}}$-regular maps it is not easy to find this kind of
examples. By Theorem 7.5 (2), we have that $P(D_{4},f)=W_{1}W_{4}$,
$P(D_{5},f)=P(D_{6},f)=0$ and%
\[
P(D_{k+1},f)=W_{k}W_{2}+W_{k-1}(W_{3}+W_{1}W_{2})\text{ \ \ for }k=6,7.
\]
Since an $\Omega^{\overline{D}}$-regular map does not have the singularity of
type $E_{k}$, we have by Theorem 7.5 (1) that $P(E_{k+1},f)=W_{k}W_{2}%
+W_{k-1}(W_{3}+W_{1}W_{2})=0$\ for $k=5,6,7$. Hence, we have that
$P(D_{k+1},f)=0$ for $k=5,6,7$.

\bigskip

Department of Mathematical Sciences, Faculty of Science,

Yamaguchi University, Yamaguchi 753-8512, Japan

e-mail: andoy@po.cc.yamaguchi-u.ac.jp
\end{document}